# THE SPREAD OF A RUMOR OR INFECTION IN A MOVING POPULATION

By Harry Kesten[1] and Vladas Sidoravicius[2]

*Cornell University and IMPA*

We consider the following interacting particle system: There is a "gas" of particles, each of which performs a continuous-time simple random walk on $\mathbb{Z}^d$, with jump rate $D_A$. These particles are called $A$-particles and move independently of each other. They are regarded as individuals who are ignorant of a rumor or are healthy. We assume that we start the system with $N_A(x,0-)$ $A$-particles at $x$, and that the $N_A(x,0-), x \in \mathbb{Z}^d$, are i.i.d., mean-$\mu_A$ Poisson random variables. In addition, there are $B$-particles which perform continuous-time simple random walks with jump rate $D_B$. We start with a finite number of $B$-particles in the system at time 0. $B$-particles are interpreted as individuals who have heard a certain rumor or who are infected. The $B$-particles move independently of each other. The only interaction is that when a $B$-particle and an $A$-particle coincide, the latter instantaneously turns into a $B$-particle.

We investigate how fast the rumor, or infection, spreads. Specifically, if $\widetilde{B}(t) := \{x \in \mathbb{Z}^d : \text{a } B\text{-particle visits } x \text{ during } [0,t]\}$ and $B(t) = \widetilde{B}(t) + [-1/2, 1/2]^d$, then we investigate the asymptotic behavior of $B(t)$. Our principal result states that if $D_A = D_B$ (so that the $A$- and $B$-particles perform the same random walk), then there exist constants $0 < C_i < \infty$ such that almost surely $\mathcal{C}(C_2 t) \subset B(t) \subset \mathcal{C}(C_1 t)$ for all large $t$, where $\mathcal{C}(r) = [-r, r]^d$. In a further paper we shall use the results presented here to prove a full "shape theorem," saying that $t^{-1}B(t)$ converges almost surely to a nonrandom set $B_0$, with the origin as an interior point, so that the true growth rate for $B(t)$ is linear in $t$.

If $D_A \neq D_B$, then we can only prove the upper bound $B(t) \subset \mathcal{C}(C_1 t)$ eventually.

Received December 2003; revised October 2004.

[1]Supported by NSF Grant DMS-99-70943 and by a Tage Erlander Professorship.

[2]Supported by FAPERJ Grant E-26/151.905/2001, CNPq (Pronex).

**1. Introduction.** We study the interacting particle system described in the first paragraph of the abstract. A construction of such a process will be discussed in the beginning of the next section.

In addition to the possible interpretations of such systems mentioned in the abstract, the $B$-particles have been interpreted as "packets of energy" which together with $A$-particles produce more energy, according to the reaction $B + A \to 2B$ (see [11]). If memory serves us well, the study of these systems was suggested by Frank Spitzer to the first author around 1980. At that time only the case when the $A$- and $B$-particles perform the same random walks (i.e., $D_A = D_B$) seems to have been considered. Recently, the so-called frog model—which has $D_A = 0$, that is, the $A$-particles do not move—has been treated by [1] and [11]. In this special case, in which the $A$-particles stand still, the model has subadditivity properties which were used to prove a full shape theorem. More specifically, it is proven in these references that there exists a nonrandom set $B_0$ such that almost surely (abbreviated to a.s. in the sequel) for all $\varepsilon > 0$

$$(1.1) \qquad (1-\varepsilon)B_0 \subset \frac{1}{t}B(t) \subset (1+\varepsilon)B_0 \qquad \text{eventually.}$$

In this paper we mainly deal with the case $D_A = D_B$. However, the upper bound for $B(t)$ (see Theorem 1 below) is relatively easy and is proven even for $D_A \neq D_B$. Probably this bound was known to several people already. It turns out that a lower bound for $B(t)$ in Theorem 2, in the case $D_A = D_B$, can be obtained by the methods of [7]. It is still an open problem whether $B(t)$ grows linearly with $t$ when $D_A > 0$, but $D_A \neq D_B$. In this case we can only prove that $B(t) \supset \mathcal{C}(K_1 t/(\log t)^p)$ eventually, for some constants $K_1, p > 0$. (We do not give the proof here.)

Throughout we shall use $N_A(x,t)$ ($N_B(x,t)$) to denote the number of $A$-particles (resp. $B$-particles) at position $x$ at time $t$. $N_B$ denotes the total number of $B$-particles at time 0. We always take $0 < N_B < \infty$ and consider $N_B$, as well as the positions of the initial $B$-particles, as nonrandom. At a site $x$ with a $B$-particle at time 0 all particles immediately turn to $B$-particles. We write $N_A(x,0-)$ for the number of $A$-particles at $x$ "just before" the $B$-particles are added to the system, and $N_B(x,0-)$ for the number of $B$-particles added at $x$. In accordance with these rules we take $N_A(x,0) = 0, N_B(x,0) = N_A(x,0-) + N_B(x,0-)$ at a site $x$ to which a $B$-particle is added at time 0. If no $B$-particle is added at $x$ at time 0, then $N_A(x,0) = N_A(x,0-)$ and $N_B(x,0) = 0$. We further define

$$\widetilde{B}(t) = \{x \in \mathbb{Z}^d \colon \text{ a } B\text{-particle visits } x \text{ during } [0,t]\},$$

$$B(t) = \widetilde{B}(t) + [-\tfrac{1}{2}, \tfrac{1}{2}]^d,$$

and the cubes

$$(1.2) \qquad \mathcal{C}(r) = [-r,r]^d.$$

Our first theorem states that the rumor/infection cannot spread from the origin faster than linearly in time.

THEOREM 1. *For some constant $C_1 < \infty$, and all sufficiently large $t$,*

$$E\{\text{number of } B\text{-particles with a position outside } \mathcal{C}(C_1 t) \text{ at time } t\} \le 2N_B e^{-t}. \tag{1.3}$$

*Consequently it is a.s. the case that*

$$B(t) \subset \mathcal{C}(2C_1 t) \qquad \text{eventually.} \tag{1.4}$$

This result holds for any $D_A, D_B \ge 0$ and probably is even valid if one allows the $A$- and $B$-particles to perform any random walk with bounded jumps of mean zero. The next theorem shows that the rumor/infection spreads at least linearly in time, but we can only prove this if both the $A$- and $B$-particles perform simple random walks with the same jump rate.

THEOREM 2. *If $D_A = D_B$, then there exists a constant $C_2 > 0$ such that for each constant $K > 0$*

$$P\{\mathcal{C}(C_2 t) \not\subset B(t)\} \le \frac{1}{t^K} \qquad \text{for all large } t. \tag{1.5}$$

*Consequently, a.s.*

$$\mathcal{C}(\tfrac{1}{2}C_2 t) \subset B(t) \qquad \text{eventually.} \tag{1.6}$$

For proving a shape theorem we will need a form of Theorem 2 which also gives some information about the possible occurrence of $A$-particles amid the spreading $B$-particles. More specifically, the same proof as for Theorem 2 can be used to prove the next theorem. This answers a question raised after a lecture on this material; unfortunately we do not remember who the questioner was.

THEOREM 3. *If $D_A = D_B$, then for all $K$ there exists a constant $C_3 = C_3(K)$ such that*

$$P\{\text{there is a vertex } z \text{ and an } A\text{-particle at the space–time point } (z,t) \text{ while there also was a } B\text{-particle at } z \text{ at some time} \le t - C_3[t \log t]^{1/2}\} \le \frac{1}{t^K} \quad \text{for all sufficiently large } t. \tag{1.7}$$

*Consequently, for large t,*

$$P\{\text{at time } t \text{ there is a site in } \mathcal{C}(C_2 t/2) \text{ which is occupied by an } A\text{-particle}\} \le \frac{2}{t^K}. \tag{1.8}$$

REMARK 1. It can be checked that the constants $C_1, C_2$ do not depend on the number or positions of the initial $B$-particles. However, the lower bounds for the times for which (1.3)–(1.6) are valid do depend on these initial data.

*Some heuristics.* The proof of Theorem 1 is basically a Peierls argument. This proof relies in part on the construction of the process given in Section 2. It associates to each $B$-particle, $\rho$ say, present at time $t$, a so-called genealogical path which describes the sequence of $B$-particles which "transmitted the rumor/infection" from the initial $B$-particles to $\rho$ at time $t$, and also describes the relevant pieces of the paths of these intermediate particles. One proves (1.3) by taking the expectation of the number of genealogical paths which lead to a $B$-particle outside $\mathcal{C}(C_1 t)$ at time $t$.

By far the most involved proof here is that of Theorem 2, which gives a lower bound on the spread of the rumor/infection. To help the intuition, it is best to think of the one-dimensional case, started with one $B$-particle at the origin and no other $B$-particle. All the major difficulties appear already in this special case. Until the last two paragraphs of these heuristic remarks we therefore take $d = 1$.

In this one-dimensional case, there is for each $t$ a rightmost $B$-particle, at position $\mathcal{R}(t)$ say, and a leftmost $B$-particle at position $-\mathcal{L}(t)$. At time $t$ all particles in $[-\mathcal{L}(t), \mathcal{R}(t)]$ are $B$-particles and all particles outside $[-\mathcal{L}(t), \mathcal{R}(t)]$ are $A$-particles. Basically we want to show that $\liminf_{t\to\infty} \mathcal{R}(t)/t > 0$ and similarly for $\mathcal{L}(t)$. If there is exactly one particle at $\mathcal{R}(t)$ at time $t$, then $\mathcal{R}(\cdot)$ behaves like a simple random walk, that is, $P\{\mathcal{R}(t+dt) = \mathcal{R}(t) \pm 1\} = Dt/2 + O(dt^2)$, with $D$ standing for the common value of $D_A$ and $D_B$. However, if there is more than one particle at $\mathcal{R}(t)$ at time $t$, then the rightmost particle moves one step to the right as soon as one of the particles at $\mathcal{R}(t)$ makes a jump to the right, whereas the rightmost position moves a step to the left only when all particles at $\mathcal{R}(t)$ move to the left. Thus, the rightmost $B$-particle has a drift to the right at all times when there is more than one particle at $\mathcal{R}(t)$. When there is at least one other particle (of either type) "close to" the rightmost $B$-particle, then there is a positive probability that in the next time unit another particle will coincide with the rightmost $B$-particle. This will still provide $\mathcal{R}(\cdot)$ with an upwards drift. By using large deviation estimates for martingales one can see that the only way for $\mathcal{R}(t)/t$ to become small (with a nonnegligible probability) is if the particle at $\mathcal{R}(s)$ has for most $s \in [0,t]$ no particle (of any type) nearby. We therefore want to show that the probability of this event goes to 0. One is tempted to try and prove this by studying the environment as seen from the position $\mathcal{R}(t)$. However, this approach seems difficult because the dependence between $\mathcal{R}(t)$ and the particles near $\mathcal{R}(t)$ is very complicated. We have been unable to use this approach. Instead, it turns out to be easier to prove a much stronger property, which uses almost no property of the path $s \mapsto \mathcal{R}(s)$. Roughly speaking we prove that *every* space–time path $s \mapsto \widehat{\pi}(s)$ with not too many jumps during $[0,t]$ has some particle "near $\widehat{\pi}(s)$ most of the time."

To make this more specific, we introduce some notation. A *path* $\pi = (x_0, \dots, x_m)$ is a sequence of integers with $x_{j+1} - x_j = \pm 1, 1 \le i \le m$. We regard the $x_j$ as the successive positions of a space–time path $\widehat{\pi}$. There are many space–time paths which traverse the same positions in the same order. A *space–time path* $\widehat{\pi}$ is specified by giving its successive positions $x_i$ and jump times $s_i$. For $s_1 < s_2 < \cdots$ we shall sometimes denote the path which jumps to $x_i$ at time $s_i$ by $\widehat{\pi}(\{s_i, x_i\})$. We make the convention that $s_0 = 0$, and unless stated otherwise, $x_0 = \mathbf{0}$. In addition we are here only discussing space–time paths over the time interval $[0,t]$, so we tacitly take $s_m \le t$. $\widehat{\pi}(\{s_i, x_i\})$ is then the path which is at position $x_i$ during $[s_i, s_{i+1})$ for $0 \le i < m$, and at position $x_m$ during $[s_m, t]$. If it is important that the path has exactly $m$ jump times, then we shall write $\widehat{\pi}(\{s_i, x_i\}_{i \le m})$. Throughout this proof we shall only consider paths which are contained in

$$\mathcal{C}(t \log t) = [-t \log t, t \log t].$$

Of particular interest for us is the following class of paths with exactly $\ell$ jumps:

$$\Xi(\ell, t) = \{\widehat{\pi}(\{s_i, x_i\}_{0 \le i \le \ell}) \quad \text{with } 0 = s_0 < s_1 < \cdots < s_\ell < t \text{ and } x_i \in \mathcal{C}(t \log t)\}. \tag{1.9}$$

Instead of using the path followed by $\mathcal{R}(\cdot)$, we shall construct special paths $\widehat{\pi}$ with the property that there is a $B$-particle at $(\widehat{\pi}(s), s)$ for all $s \le t$ [so that automatically $\mathcal{R}(t) \ge \widehat{\pi}(t)$], and such that these paths are with high probability in $\Xi(\ell, t)$ for some $\ell \le 2Dt$, and also have a drift to the right at any time $s$ when there are at least two particles at $\widehat{\pi}(s)$. Thus, it will be sufficient to show that every space–time path $\widehat{\pi} \in \bigcup_{\ell \le 2Dt} \Xi(\ell, t)$ has some particle "near $\widehat{\pi}(s)$ most of the time."

To this end we choose a large integer $C_0$ and partition space–time $\mathbb{Z} \times [0, \infty)$ into the following blocks of size $\Delta_r := C_0^{6r}$:

$$\mathcal{B}_r(i,k) = [i\Delta_r, (i+1)\Delta_r) \times [k\Delta_r, (k+1)\Delta_r).$$

We call these intervals $r$-blocks. We shall soon define "good" and "bad" $r$-blocks. There is a standard percolation argument which also partitions space into large blocks which can be good or bad, and then shows that on the one hand the bad blocks do not percolate, and on the other hand that no percolation of bad blocks implies a desired property. In our case the desired property would be that any space–time path $\widehat{\pi} \in \Xi(\ell, t)$ intersects at most $\varepsilon t$ bad $r$-blocks for a suitable $\ell$ and for a small $\varepsilon$ [see (1.17) below]. This is indeed the desired property we are after, but we have not succeeded in simply working with $r$-blocks for one fixed $r$, because of the complicated dependence of the configurations in different $r$-blocks. Instead we work with $r$-blocks for all $r$. This is why we say that our proof is based on a multiscale argument. We also need the following sets (see Figure 1):

$$\begin{aligned}\widetilde{\mathcal{B}}_r(i,k) &:= [(i-3)\Delta_r, (i+4)\Delta_r) \times [(k-1)\Delta_r, (k+1)\Delta_r),\\ \mathcal{Q}_r(x) &:= [x, x + C_0^r),\\ V_r(i) &:= [(i-3)\Delta_r, (i+4)\Delta_r)\end{aligned}$$

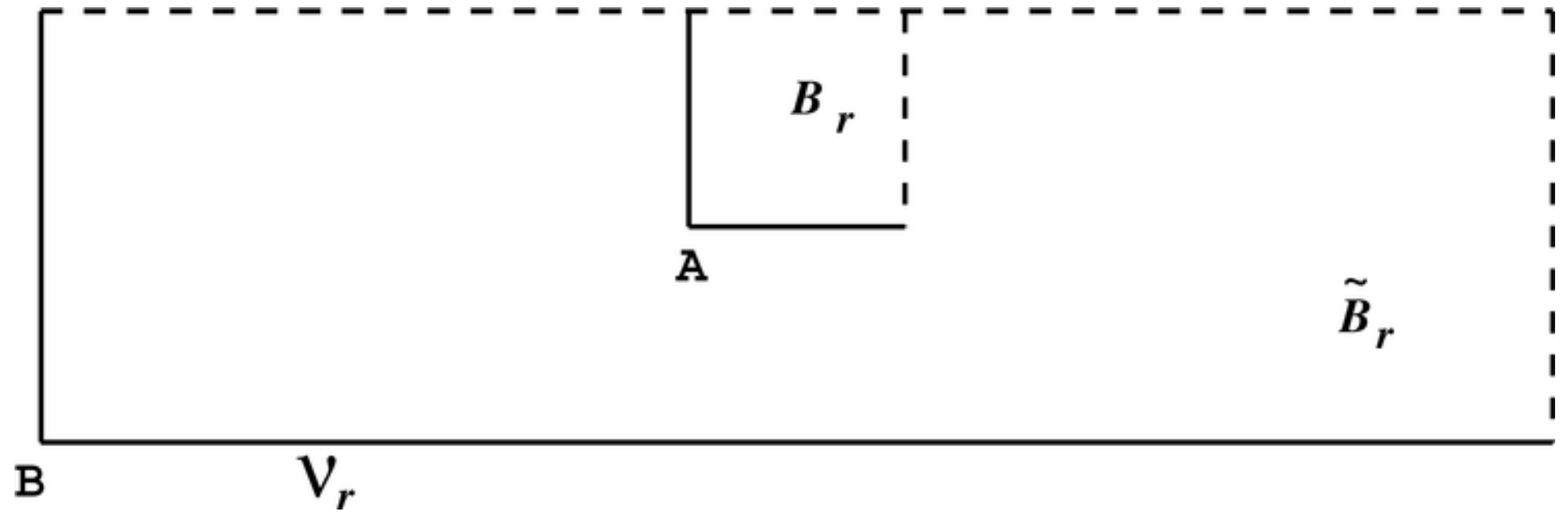

FIG. 1. *Relative location of the sets $\mathcal{B}_r(i,k), \widetilde{\mathcal{B}}_r(i,k)$ and $\mathcal{V}_r(i,k)$ for $d=1$. These sets are "left closed, right open," that is, the solid segments are in the sets, but the dashed segments are not. The space and time directions are along the horizontal and vertical axes, respectively. The points $A$ and $B$ are the space–time points $(i\Delta_r, k\Delta_r)$ and $((i-3)\Delta_r, (k-1)\Delta_r)$, respectively. $\mathcal{V}_r$ is the line segment which constitutes the bottom of $\widetilde{\mathcal{B}}_r$.*

(these are just intervals of length $C_0^r$ and $7\Delta_r$, resp.), and the *pedestal* of $\mathcal{B}_r(i,k)$:

$$\mathcal{V}_r(i,k) = V_r(i) \times \{(k-1)\Delta_r\}.$$

We also need to count numbers of particles in certain sets. We define $N^*(x,t)$ as the number of particles at the space–time point $(x,t)$ in the system which evolves freely, without any $B$-particles. In this system—which we sometimes denote by $\mathcal{P}^*$—we start off with $N_A(x,0) = N_A(x,0-)$ particles at $x$ at time 0 and let all these particles perform independent random walks without any interaction. Note that $N^*(x,t) \le N_A(x,t) + N_B(x,t)$. The important counts are

$$U_r(x,v) = \sum_{y \in \mathcal{Q}_r(x)} N^*(y,v) = \sum_{y : x \le y < x + C_0^r} N^*(y,v). \tag{1.10}$$

(We shall need this only for integer times $v$.) We call the $r$-block $\mathcal{B}_r(i,k)$ *bad* if

$$U_r(x,v) < \gamma_r \mu_A C_0^r \qquad \text{for some } (x,v) \text{ with integer } v \text{ for which}$$
$$\mathcal{Q}_r(x) \times \{v\} \text{ is contained in } \widetilde{\mathcal{B}}_r(i,k).$$

The $\gamma_r$ in this definition are given by (4.3) below. For the time being the only important properties are that the $\gamma_r$ are strictly increasing (but slowly) and satisfy

$$0 < \gamma_0 < \gamma_r < \gamma_\infty \leq \tfrac{1}{2}, \qquad r > 0. \tag{1.11}$$

Roughly speaking, the bad blocks are blocks in which the number of $A$-particles in some spacelike cube of specified size and which is nearby in space–time, is less than half the expected amount. Indeed, it is well known that in our setup each $U_r(x,v)$ has a Poisson distribution of mean $\mu_A C_0^r$. The pedestal $\mathcal{V}_r(i,k)$ is called bad if

$$U_r(x,(k-1)\Delta_r) < \gamma_r \mu_A C_0^r \qquad \text{for some } x \text{ with } \mathcal{Q}_r(x) \subset V_r(i).$$

A block (resp. pedestal) is called *good* if it is not bad.

If a space–time path $\widehat{\pi}$ is in a good $r$-block at a given time $s$, then there are a reasonable number of particles within distance $C_0^r$ of $\widehat{\pi}(s)$ at time $s$, by definition of a good block. We therefore would like to show that "most" space–time paths intersect "few" bad blocks during $[0,t]$. To quantify this statement we define

$$\phi_r(\widehat{\pi}) = \text{number of bad } r\text{-blocks which intersect the space–time path } \widehat{\pi}, \tag{1.12}$$

$$\Phi_r(\ell) = \sup_{\widehat{\pi} \in \Xi(\ell,t)} \phi_r(\widehat{\pi}), \tag{1.13}$$

$$\psi_{r+1}(\widehat{\pi}) = \text{number of } (r+1)\text{-blocks which intersect the space–time path } \widehat{\pi} \text{ and which have a good pedestal but contain a bad } r\text{-block} \tag{1.14}$$

and

$$\Psi_r(\ell) = \sup_{\widehat{\pi} \in \Xi(\ell,t)} \psi_r(\widehat{\pi}) \tag{1.15}$$

(we suppress the dependence on $t$ in these quantities). The principal part of the proof is to show that for any choice of $K > 0$ and $\varepsilon_0 > 0$ there exists an $r_0$ such that

$$P\{\Phi_r(\ell) \geq \varepsilon_0 C_0^{-6r}(t+\ell) \text{ for some } r \geq r_0, \ell \geq 0\} \leq \frac{2}{t^K} \tag{1.16}$$

for all large $t$ (see Proposition 8). This result has the desired form, because any path $\widehat{\pi}$ spends at most $C_0^{6r}$ time units in a given $r$-block, and therefore at most $C_0^{6r}\phi_r(\widehat{\pi})$ time units in bad blocks during $[0,t]$. Moreover, as we stated before, we only need to consider space–time paths in $\bigcup_{\ell\le 2Dt}\Xi(\ell,t)$. Thus if the property in braces in (1.16) holds, then any $\widehat{\pi}\in\bigcup_{\ell\le 2Dt}\Xi(\ell,t)$ satisfies

$$C_0^{6r}\phi_r(\widehat{\pi})\le C_0^{6r}\sup_{\ell\le zDt}\Phi_r(\ell)\le\varepsilon_0(1+2D)t, \tag{1.17}$$

and spends at most $\varepsilon_0(1+2D)t$ time units in bad blocks (for $r\ge r_0$). For $\varepsilon=\varepsilon_0(1+2D)<1/2$ this shows that the paths of interest to us have a drift to the right for at least $t/2$ time units.

This leaves us with the problem of proving (1.16). This is done by means of a recurrence relation (with random terms) for the $\Phi_r$. Note that each bad $r$-block has to lie either in a good $(r+1)$-block or in a bad $(r+1)$-block. Since any $(r+1)$-block contains exactly $C_0^{12}$ $r$-blocks, the number of bad $r$-blocks which intersect a path $\widehat{\pi}$, and which are contained in a bad $(r+1)$-block (and which necessarily intersects $\widehat{\pi}$) is at most $C_0^{12}\phi_{r+1}(\widehat{\pi})\le C_0^{12}\Phi_{r+1}$. A similar estimate holds for the number of bad $r$-blocks which intersect $\widehat{\pi}$ and which are contained in a good $(r+1)$-block. If one also takes into account that any good block has a good pedestal, by definition, then it is not hard to see that

$$\phi_r(\widehat{\pi})\le C_0^{12}\Phi_{r+1}(\ell)+C_0^{12}\psi_{r+1}(\widehat{\pi}).$$

In turn, by taking the sup over $\widehat{\pi}\in\Xi(\ell,t)$ this gives

$$\Phi_r(\ell)\le C_0^{12}\Phi_{r+1}(\ell)+C_0^{12}\Psi_{r+1}(\ell). \tag{1.18}$$

In addition, it follows from simple estimates for Poisson variable that outside a set of probability $t^{-K}$ there are no spacelike intervals $Q_r(x)$ which intersect $[-t\log t,t\log t]$ and with $U_r(x,s)<\gamma_r\mu_A C_0^r$ for any $r\ge R(t)$, where $R(t)$ is the unique integer with $C_0^{R(t)}\ge K_4\log t>C_0^{R(t)-1}$ (for a suitable constant $K_4$). Thus with high probability $\Phi_r(\ell)=0$ for all $\ell$ and $r\ge R(t)$. Thus we can start with the "boundary condition" $\Phi_{R(t)}(\ell)=0$ and then work our way downward to conclude that also $\Phi_{r_0}(\ell)$ is $o(t)$ for some fixed $r_0$, provided we can show that the $\Psi_r(\ell)$ are suitably small (with high probability). This last fact is shown by using the following lower bound for the $U_r(x,v)$: Let $\mathcal{B}_{r+1}(i,k)$ be the unique $(r+1)$-block which contains $(x,v)$. Then define

$$\begin{aligned}W_r(x,v)=&\text{ number of particles in the system }\mathcal{P}^*\text{ in }\mathcal{Q}_r(x)\times\{v\}\\&\text{ which were in }V_{r+1}(i)\text{ at time }(k-1)\Delta_{r+1}.\end{aligned}$$

We call the $r$-block $\mathcal{B}_r(i,k)$ *inferior* if

$$W_r(x,v)<\gamma_r\mu_A C_0^r\qquad\text{for some }(x,v)\text{ with integer }v\text{ for which }\mathcal{Q}_r(x)\times\{v\}\text{ is contained in }\widetilde{\mathcal{B}}_r(i,k).$$

It is apparent from the definitions that $W_r(x,v) \le U_r(x,v)$. Therefore, a bad block is also inferior and it suffices to show that $\widetilde{\Psi}_r(\ell)$ is small, where $\widetilde{\Psi}$ is defined by changing "bad" in the definition (1.14) to "inferior." Now let

$$\mathcal{A}(i,k) = \mathcal{A}(i,k,r) := \{\mathcal{B}_{r+1}(i,k) \text{ contains some inferior } r\text{-block } \mathcal{B}_r(j,q)\}.$$

The advantage of the $W_r$ over the $U_r$ is that they lead to much better independence properties of the $\mathcal{A}(i,k)$ than if we had defined $\mathcal{A}(i,k)$ with "some bad $r$-block" instead of "some inferior $r$-block." In fact, once we know which particles are in the pedestal $\mathcal{V}_{r+1}(i,k)$ of $\mathcal{B}_{r+1}(i,k)$, whether or not $\mathcal{A}(i,k)$ occurs depends only on the particles in $\mathcal{V}_{r+1}(i,k)$ and not on particles in any pedestal $\mathcal{V}_{r+1}(j,v)$ with $v \le k$ and disjoint from $\mathcal{V}_{r+1}(i,k)$. With a little more work one shows that for fixed $a \in \{0,1,\dots,11\}$ and $b \in \{0,1\}$ the collection of pairs $(i,k)$ with $i \equiv a \bmod 12, k \equiv b \bmod 2$ for which $\mathcal{V}_{r+1}(i,k)$ is good, but $\mathcal{A}(i,k)$ occurs, is stochastically smaller than an independent percolation system in which each site $(i,k), i \equiv a \bmod 12, k \equiv b \bmod 2$ has probability $\rho_{r+1}$ of being open. Here $\rho_{r+1}$ is an upper bound for the probability that an $(r+1)$-block with a good pedestal contains an inferior $r$-block. We shall show in Lemma 6, by straightforward large deviation estimates, that in dimension 1 we can take

$$\rho_{r+1} = 9C_0^{12(r+1)} \exp[-\tfrac{1}{2}\gamma_r \mu_A C_0^{r/4}].$$

It is for this estimate that the $\gamma_r$ are chosen strictly increasing. Roughly speaking, a good $(r+1)$-block has density at least $\gamma_{r+1}\mu_A$ of particles in its pedestal. It is then possible to bound the probability that such an $(r+1)$-block contains an $r$-block with density $\le \gamma_r \mu_A$ for a suitable $\gamma_r < \gamma_{r+1}$.

From here on one can follow known arguments from percolation and large deviations to obtain an estimate for the tail of the distribution of $\Psi_{r+1}(\ell)$ (see Lemma 7). Finally, the recurrence relation (1.18) then gives (1.16) in the one-dimensional case. As pointed out before, (1.16) guarantees that with high probability every relevant space–time path has drift to the right for at least half the time and this is enough to obtain $\liminf_{t\to\infty} \mathcal{R}(t)/t > 0$ a.s.

At this stage it may be useful to say a few words about the case of dimension greater than 1. There is no clear analogue of $\mathcal{R}(t)$, or at least none that is helpful. Instead of constructing paths which have a drift to the right at least half the time one now fixes an $x \in \mathcal{C}(C_2 t) \cap \mathbb{Z}^d$ and tries to construct a space–time path $\lambda(\cdot) = \lambda(\cdot, x)$ which has a $B$-particle at $\lambda(s)$ for all $s$, and which has a tendency to move toward $x$. In fact, our $\lambda(s)$ behaves like a ($d$-dimensional) simple random walk at times $s$ when there is only one particle at $\widehat{\pi}(s)$, but if there are at least two particles at $\lambda(s)$ and a particle jumps away from $\lambda(s-)$ at time $s$, then the conditional expectation of $\|\lambda(s) - x\|_2$ is smaller than $\|\lambda(s-) - x\|_2$. This will give us a path which with high probability reaches $x$ during $[0,t]$, provided the path has at least two particles "near $\lambda(s)$" at least a positive fraction of the time. In this way all points $x \in \mathcal{C}(C_2 t) \cap \mathbb{Z}^d$ are reached by the infection during $[0, t]$. From there on there are only minor differences between the cases $d = 1$ and $d > 1$ for Theorem 2.

Theorem 3 is very similar to Theorem 2. Roughly speaking, if there is a $B$-particle at a given site $x$ at some time $s \le t - C_3[t \log t]^{1/2}$, then by the estimates for Theorem 2 each site $z \in x + \mathcal{C}(C_2(t - s))$ will be reached by some $B$-particle during $[s, t]$. This is proven by constructing some random path from $x$ to $z$ in the same manner as in the last pargaraph. However, for Theorem 3 there is a small difference. We do not need a $B$-particle which reaches a given site $z$, but any particle which is at $z$ at time $t$ should have coincided with a $B$-particle during $[s, t]$ (or already have type $B$ itself at time $s$). To show that this is the case we construct a random path which has a tendency to move toward such a moving particle, rather than toward the fixed site $z$. Only trivial changes in the construction of useful paths are required, but no real changes in the estimates are needed.

In Section 2 we describe a possible construction of our process, but we do not give a proof here that this construction results in a strong Markov process. A proof of this fact can be found in an earlier version of this paper (see [8]). The proof of the upper bound for the spread of the rumor/infection, that is, of Theorem 1, is given in Section 3. The rather involved proof of Theorem 2 is given in Section 4. Finally, the proof of Theorem 3 is similar to that of Theorem 2. This proof is given in Section 5.

## 2. Construction of a strong Markov process.

Throughout this paper we make the following convention about constants. $K_i$ will denote a strictly positive, finite constant, whose precise value is unimportant for our purposes. The value of the same $K_i$ may be different in different formulas. We use $C_i$ for constants whose value remains fixed throughout the paper. They will again have values in $(0, \infty)$. If necessary, we indicate on what other quantities a constant depends at the time when it is first introduced. Throughout $\|x\|$ denotes the $\ell^\infty$ norm of the vector $x = (x(1), \dots, x(d)) \in \mathbb{R}^d$, that is,

$$\|x\| = \max_{1 \le i \le d} |x(i)|. \tag{2.1}$$

$\mathbf{0}$ will denote the origin (in $\mathbb{Z}^d$ or $\mathbb{R}^d$).

In this section we shall indicate how to construct our process on a suitable probability space as a strong Markov process. We shall skip most proofs. Even though the main results in this paper are for the case $D_A = D_B$, we do not make this assumption yet, so that we can prove Theorem 1 also if $D_A \ne D_B$. A complete proof of the strong Markov property, even for the case $D_A \ne D_B$, can be found in an earlier version of this paper (see [8]). The usual way to prove that an interacting particle system can be represented by a strong Markov process is to construct the process as a function from a probability space into a state space with a suitable topology in which the process is right-continuous and has the Feller property. We did not succeed in finding such a topology. We were only able to construct our process as a right-continuous process $\{Y_t\}$ for which $t \mapsto Q(Y_t, \mathcal{E})$ is almost surely right-continuous for sufficiently many $\mathcal{E}$, where $Q(y, \cdot)$ are the transition probabilities from $y$. This suffices for the strong Markov property (see [5], Theorem 5.10 and the remark following it, or [3], Theorem I.8.11 and its proof). However, we need a somewhat involved definition of $Y_t$ as a function on a probability space.

We want to construct our process as a Markov process with a given initial state which contains only finitely many $B$-particles. Our first task is to choose the state space $\Sigma_0$ for our process. We shall assume that there are countably many particles in our system, which are labeled $\rho_1, \rho_2, \dots$. A particle keeps the same label throughout. The state of our system is described by specifying the location and type of each particle. We shall also add an artificial cemetery point $\partial$ for each particle to its coordinate space. Thus, the state space will be taken as a subset of $\Sigma := \prod_{k\geq 1}((\mathbb{Z}^d \cup \partial_k) \times \{A, B\})$. If $\sigma = (\sigma'(k), \sigma''(k))$ is a generic point of $\Sigma$, then $\sigma'(k)$ represents the position of $\rho_k$ and $\sigma''(k)$ represents the type of $\rho_k$. Occasionally it will be more convenient to use the notation $\sigma'(\rho), \sigma''(\rho), \partial(\rho)$ for the position, type and cemetery point of a particle $\rho$, without specifying which of the particles $\rho_k$ equals $\rho$. To describe the state space $\Sigma_0$ we introduce a process $\{Y_t\}_{t\geq 0}$. A priori, each $Y_t$ takes values in $\Sigma$. Later we add conditions to make sure that $Y_t$ takes values in $\Sigma_0$. We need some definitions. $\{S_t^\eta\}_{t\geq 0}$ will be a random walk with the same distribution as the random walks performed by the particles of type $\eta$ (with $S_0^\eta = \mathbf{0}$, $\eta = A$ or $B$). We further attach to each particle $\rho$ present at time 0 two random walk paths $t \mapsto \pi_A(t, \rho)$ and $t \mapsto \pi_B(t, \rho)$. Each $\{\pi_\eta(t, \rho)\}_{t\geq 0}$ has the same distribution as $\{S_t^\eta\}_{t\geq 0}$. For the case $D_A \neq D_B$ all these paths are chosen independently. For the case $D_A = D_B$ it is more convenient to take $\pi_A = \pi_B$, so that for each particle only the paths $\pi_A$ have to be chosen, with the $\pi_A(\cdot, \rho)$ for different $\rho$ completely independent. Note that we take all these paths right-continuous. We write $\pi(t, \rho)$ and $\eta(t, \rho)$ for the position and type of $\rho$ at time $t$, respectively.

We want to let an $A$-particle $\rho$ which starts at $z$ move along the path $t \mapsto z + \pi_A(t, \rho)$ until the time $\theta(\rho)$, say, at which it changes to a $B$-particle, or until $\rho$ is moved to its cemetery point, if this time comes at or before $\theta(\rho)$ (see below). If $\theta(\rho)$ comes before $\rho$ is moved to its cemetery point, then from $\theta(\rho)$ on, $\rho$ follows the path $t \mapsto z + \pi_A(\theta(\rho), \rho) + \pi_B(t, \rho) - \pi_B(\theta(\rho), \rho)$. For all $\rho$ which have type $B$ at time 0, we take $\theta(\rho) = 0$ and let $\rho$ move along the path $t \mapsto z + \pi_B(t, \rho)$ for all $t \geq 0$ ($z$ again denotes the initial position of $\rho$). Also for the case $D_A = D_B$ each particle $\rho$ moves along the single path $t \mapsto z + \pi_A(t, \rho)$ until the time at which $\rho$ is moved to its cemetery point (this time may be infinite). Formally, we proceed as follows. We assume that initially there are in total only finitely many $B$-particles, and that none of these sits at a cemetery point. We set $\tau_0 = 0$. Now let $k = 0$, or let $k \geq 1$ and assume that we have already found the first $k$ times $\tau_1 \leq \tau_2 \leq \cdots \leq \tau_k$ at which a $B$-particle has coincided with an $A$-particle. We also assume that at each of these times only finitely many $A$-particles turned into $B$-particles, so that at time $\tau_k$ there are still only finitely many $B$-particles in the system. Assume further that we have determined the paths of all particles during the interval $[0, \tau_k]$. Then we know at time $\tau_k$ which particles are $B$-particles and also the positions of all particles. We then assign to each particle $\rho$ the tentative continuation of its path on $[\tau_k, \infty)$, which it would follow if it never changed type anymore. The tentative continuation of the particle paths is given by

$$\begin{aligned}\widetilde{\pi}_k(\tau_k + t, \rho) \\ = \begin{cases} \pi(\tau_k, \rho) + [\pi_A(\tau_k + t, \rho) - \pi_A(\tau_k, \rho)], & \text{if } \eta(\tau_k, \rho) = A, \\ \pi(\tau_k, \rho) + [\pi_B(\tau_k + t, \rho) - \pi_B(\tau_k, \rho)], & \text{if } \eta(\tau_k, \rho) = B. \end{cases}\end{aligned} \tag{2.2}$$

We have to allow that some particles sit at their cemetery point. We therefore interpret the right-hand side of (2.2) as $\partial(\rho)$ if $\pi(\tau_k, \rho) = \partial(\rho)$. As the reader can check in the definitions below, this has the effect that any particle stays at its cemetery point once it reaches this cemetery point. After that such a particle no longer interacts with the other particles and plays no further role in the construction of the paths of the other particles. We now use these $\widetilde{\pi}_k$ to define

$$\begin{aligned}\tau_{k+1} &= \inf\{t > \tau_k : \text{ a } B\text{-particle coincides with an } A\text{-particle at time } t \\ &\qquad\qquad \text{if the particles move according to the } \widetilde{\pi}_k\} \\ &= \inf\{t > \tau_k : \widetilde{\pi}_k(t, \rho') = \widetilde{\pi}_k(t, \rho'') \text{ for some } \rho', \rho'' \\ &\qquad\qquad \text{with } \eta(\tau_k, \rho') = B, \eta(\tau_k, \rho'') = A\}.\end{aligned} \tag{2.3}$$

We then take

$$\pi(s, \rho) = \begin{cases} \widetilde{\pi}_k(s, \rho), & \text{for } \tau_k \leq s \leq \tau_{k+1} \text{ if } \eta(\tau_k, \rho) = A, \\ \widetilde{\pi}_k(s, \rho), & \text{for } s \geq \tau_k \text{ if } \eta(\tau_k, \rho) = B. \end{cases} \tag{2.4}$$

Moreover,

$$\eta(s, \rho) = \begin{cases} A, & \text{for } \tau_k \leq s < \tau_{k+1} \text{ if } \eta(\tau_k, \rho) = A, \\ B, & \text{for } s \geq \tau_k \text{ if } \eta(\tau_k, \rho) = B. \end{cases} \tag{2.5}$$

In addition we take $\eta(\tau_{k+1}, \rho) = B$ for those $\rho$ which have $\eta(\tau_k, \rho) = A$ and which coincide at time $\tau_{k+1}$ with a $\rho'$ which has $\eta(\tau_k, \rho') = B$. For this special set of particles $\rho$ we take $\theta(\rho) = \tau_{k+1}$ and call $\theta(\rho)$ the *switching time* of $\rho$. For all other particles their type remains unchanged at $\tau_{k+1}$. If $\rho$ is already of type $B$ at time 0, then we define its switching time to be 0. Note that if there are sites with both $A$- and $B$-particles in $\sigma$, then $\tau_1 = 0$ according to (2.3), and all $A$-particles which are at the same location as a $B$-particle in $\sigma$ immediately change their type to $B$. These particles have switching time equal to 0.

These definitions give us $Y_t$ through time $\tau_{k+1}$ and we can repeat the procedure to go till time $\tau_{k+2}$, and so on. We stop the process at

$$\begin{aligned}\widehat{\tau} := \inf\{\tau_k : &\text{infinitely many } A\text{-particles} \\ &\text{turn into a } B\text{-particle at time } \tau_k, \text{ or } \tau_{k+1} = \tau_k\}.\end{aligned}$$

Note that a.s. $\tau_{k+1} = \tau_k$ can occur only if there are coincidences of $B$- and $A$-particles immediately after $\tau_k$, so that there must be infinitely many $B$-particles at $\tau_k + \varepsilon$ for any $\varepsilon > 0$. (E.g., such a situation would arise if at some time there are infinitely many particles at a site $x$ and a $B$-particle adjacent to $x$.) We shall actually choose $\Sigma_0$ such that this possibility has probability 0. We also cannot continue beyond $\tau_\infty := \lim_{k\to\infty} \tau_k$. We define for $t < \min\{\widehat{\tau}, \tau_\infty\}$,

$$\nu(t) = \text{total number of } B\text{-particles at time } t$$

and

$$Y_t'(\rho) = \pi(t,\rho), \qquad Y_t''(\rho) = \eta(t,\rho).$$

If $\min\{\widehat{\tau}, \tau_\infty\} > 0$ and $t \geq \min\{\widehat{\tau}, \tau_\infty\}$, then we take

$$Y_t(\rho) = \begin{cases} (\partial(\rho), A), & \text{if } \eta(s,\rho) = A \text{ for all } s < \min\{\widehat{\tau}, \tau_\infty\}, \\ (\pi(\theta(\rho),\rho) + \pi_B(t,\rho) - \pi_B(\theta(\rho),\rho), B), & \\ \qquad\qquad \text{if } \theta(\rho) < \min\{\widehat{\tau}, \tau_\infty\}. \end{cases} \tag{2.6}$$

If $\min\{\widehat{\tau}, \tau_\infty\} = 0$, then we take for $t \geq 0$

$$Y_t(\rho) = \begin{cases} (\partial(\rho), A), & \text{if } \eta(0,\rho) = A, \\ (\pi(0,\rho) + \pi_B(t,\rho), B), & \text{if } \eta(0,\rho) = B. \end{cases} \tag{2.7}$$

We further take

$$\nu(t) = \infty \qquad \text{for } t \geq \min\{\widehat{\tau}, \tau_\infty\}.$$

Thus, at $\min\{\widehat{\tau}, \tau_\infty\}$ all particles which still have type $A$ are moved to their cemetery, while the $B$-particles continue as $B$-particles along the appropriate path prescribed by their $\pi_B$. Since we start off with no $B$-particles at any cemetery point, the relations (2.2), (2.6) and (2.7) guarantee that there never are $B$-particles at the cemetery points. Thus $\nu(t)$ is actually the number of $B$-particles in $\mathbb{Z}^d$ at time $t$.

The preceding defines for each initial state $\sigma$ with finitely many $B$-particles a process $\{Y_t\}_{t\geq 0}$. We write $P^\sigma$ for the probability measure governing this process.

We easily see that our definitions give us the following three properties which agree with the intuitive description of our system:

$$\begin{aligned}&\text{if } \rho \text{ is already of type } B \text{ at time } \tau_k,\\ &\text{then it will stay of type } B \text{ for all } t \geq \tau_k;\end{aligned} \tag{2.8}$$

(note that if $\eta(\tau_k, \rho) = B$, then we have two possible prescriptions for $\pi(s, \rho)$ and $\eta(s, \rho)$ on $[\tau_{k+1}, \infty)$, one using (2.4) and (2.5) as written, and the other using (2.4) and (2.5) with $k$ replaced by $k + 1$, but these two prescriptions agree)

$$\begin{aligned}&\text{if } \rho \text{ has type } A \text{ at time } \tau_k, \text{ then it must have been}\\ &\text{of type } A \text{ during the whole interval } [0, \tau_k] \text{ and}\\ &\pi(s, \rho) = \pi(0, \rho) + \pi_A(s, \rho) \text{ for } s \in [0, \tau_k];\end{aligned} \tag{2.9}$$

$$\begin{aligned}&\text{once } \rho \text{ has become of type } B, \text{ then its position changes according}\\ &\text{to } \pi_B(\cdot, \rho), \text{ that is, } \pi(s'', \rho) - \pi(s', \rho) = \pi_B(s'', \rho) - \pi_B(s', \rho)\\ &\text{for } s'' \geq s' \geq \theta(\rho).\end{aligned} \tag{2.10}$$

We also point out that $\nu(t) < \infty$ for all $t < \min\{\widehat{\tau}, \tau_\infty\}$, directly from the definitions. Finally, we define

$$\begin{aligned}\Sigma_0 = \{\sigma \in (\mathbb{Z}^d \times \{A, B\})^{\mathbb{Z}_+} :\\ 1 \leq \text{ (number of } B\text{-particles in } \sigma) < \infty,\\ \text{and } P^\sigma\{\min\{\widehat{\tau}, \tau_\infty\} = \infty | Y_0 = \sigma\} = 1\}.\end{aligned} \tag{2.11}$$

Note that $\sigma \in \Sigma_0$ requires that none of the particles in $\sigma$ are at their cemetery point.

The next two lemmas and Proposition 3 state that $\Sigma_0$ is a good state space for the process $\{Y_t\}$ and that $\{Y_t\}$ restricted to $\Sigma_0$ has the strong Markov property. We expect that most readers will be content to accept this without proof. We therefore do not give their proofs here, but refer the interested reader to [8] for the proofs. Proposition 3 shows that under a product measure of mean-$\mu_A$ Poisson variables for the numbers of $A$-particles on the sites of $\mathbb{Z}^d$, almost all choices lead to an initial point in $\Sigma_0$. In particular $\Sigma_0 \neq \varnothing$. Note, however, that in Lemmas 1, 2 and 3 the numbers of initial $A$-particles at the various sites are not random. The initial state there is any point of $\Sigma$ or $\Sigma_0$, respectively. The basic $\sigma$-fields which we shall use are

$$\mathcal{F}_t^0 := \sigma\text{-field generated by } \{Y_s : s \leq t\}. \tag{2.12}$$

The elements of these $\sigma$-fields are subsets of $\Sigma^{[0,\infty)}$, the path space for $\{Y_t\}_{t\geq 0}$. The coordinate spaces of $\Sigma$, that is, the spaces $(\mathbb{Z}^d \cup \partial_k) \times \{A, B\}$, are countable. We endow them with the discrete topology and use the product of these topologies on $\Sigma$.

Unfortunately the description of $\Sigma_0$ is not very explicit, and it may seem useless to go through such length to find such a state space. Instead one might choose to work only with the process starting with independent Poisson numbers of particles at the sites of $\mathbb{Z}^d$. However, we know of no way to prove that such a process has the strong Markov property without describing the state space $\Sigma_0$, and our proofs use the strong Markov property at several places.

LEMMA 1. *The process $\{Y_t\}_{t\geq 0}$ is a Markov process on $\Sigma$ with respect to the filtration $\{\mathcal{F}_t^0\}_{t\geq 0}$. Its transition function equals*

$$Q_s(\sigma, \Gamma) = P^\sigma\{Y_s \in \Gamma\}, \qquad s \geq 0, \Gamma \subset \Sigma. \tag{2.13}$$

*Moreover, $t \mapsto Y_t$ is right-continuous if we use the product topology on $\Sigma$.*

To formulate the next lemma we define

$$\alpha_t(z) = P\{S_t^A = -z\} \tag{2.14}$$

and

$$M_s(\sigma) = \sum_{z\in\mathbb{Z}^d} \alpha_s(z) \sum_{\rho\,:\,\sigma(\rho)=(z,A)} 1 = \sum_\rho I[\sigma''(\rho) = A]\alpha_s(\sigma'(\rho)). \tag{2.15}$$

For purposes of comparison it is useful to couple our system with the system in which there are never any $B$-particles and in which all original $A$-particles move forever without interaction. In this system, which we shall denote by $\mathcal{P}^*$ (and which was already mentioned in the heuristic comments in Section 1), an $A$-particle $\rho$ which starts at $z$ will have position $z + \pi_A(t, \rho)$ for all $t$. Thus it coincides with this same particle in the $Y$-process until the minimum of $\theta(\rho)$ and the time at which $\rho$ is moved to its cemetery (if this time is finite). After this time, the increments of $\rho$ in the $Y$-process will be the same as those of $\pi_B(\cdot, \rho)$, or these increments will be 0, while in the $\mathcal{P}^*$ system, the increments of $\rho$ will be the same as those of $\pi_A(\cdot, \rho)$. We write $N^*(x, t)$ for the number of particles at the space–time point $(x, t)$ in the system $\mathcal{P}^*$. $N^*(x, 0)$ is taken equal to $N_A(x, 0-)$, the initial number of $A$-particles at $x$. No initial $B$-particles are introduced in $\mathcal{P}^*$ and all particles have type $A$ forever in $\mathcal{P}^*$. For $x \in \mathbb{Z}^d$, $N^*(x, t)$ is an upper bound for the number of $A$-particles at $(x, t)$ in our original system, because in that system $A$-particles can turn into $B$-particles at some time, after which they are no longer counted in $N_A$. Thus

$$\sum_{\rho\,:\,Y_t(\rho)=(x,A)} 1 \leq N^*(x, t). \tag{2.16}$$

One more piece of notation: We shall write $P^\sigma$ for the measure governing the process $\{Y_t\}$ given that it starts with $Y_0 = \sigma$. This is the unique measure on the space of right-continuous paths into $\Sigma$ with finite-dimensional distributions given by (2.13). These finite-dimensional distributions are determined by

$$\begin{aligned} &P\{Y_{t_i} \in \Gamma_i, 1 \le i \le k | Y_0 = \sigma\} \\ &\qquad = \int_{\sigma_1 \in \Gamma_1} \cdots \int_{\sigma_k \in \Gamma_k} Q_{t_1}(\sigma, d\sigma_1) \cdots Q_{t_k - t_{k-1}}(\sigma_{k-1}, d\sigma_k) \end{aligned} \tag{2.17}$$

for $\Gamma_i \subset \Sigma$ and $0 < t_1 < t_2 < \cdots < t_k$. $E^\sigma$ denotes expectation with respect to $P^\sigma$.

Lemma 2. *Fix the initial state $\sigma \in \Sigma_0$. Then almost surely $[P^\sigma]$ the following properties hold*:

$$\min\{\widehat{\tau}, \tau_\infty\} = \infty; \tag{2.18}$$

$$M_s(Y_t) < \infty \qquad \textit{for all } s, t \ge 0; \tag{2.19}$$

*and for all $z \in \mathbb{Z}, t < \infty$*

$$(\textit{number of particles which visit } z \textit{ during } [0, t]) < \infty. \tag{2.20}$$

Proposition 3. *For each $\sigma \in \Sigma_0$ one has*

$$P^\sigma\{Y_t \in \Sigma_0 \textit{ for all } t \ge 0\} = 1. \tag{2.21}$$

*Also, a.s. $[P^\sigma]$,*

$$\begin{aligned} &\textit{for all } t, s \ge 0 \\ &\quad P^{Y_t}\{Y_s'(\rho) = \partial(\rho) \textit{ for some } \rho\} \le P^{Y_t}\{\min\{\widehat{\tau}, \tau_\infty\} \le s\} = 0. \end{aligned} \tag{2.22}$$

*Moreover, if $\sigma \in \Sigma_0$ and $\mathcal{E}$ is a finite union of sets of the form*

$$\begin{aligned} &\{Y_{s_j}(\rho_j) = (z_j, \eta_j), 1 \le j \le n\} \\ &\qquad = \{\pi(s_j, \rho_j) = z_j, \eta(s_j, \rho_j) = \eta_j, 1 \le j \le n\} \end{aligned} \tag{2.23}$$

*for some fixed $z_j \in \mathbb{Z}^d, \eta_j \in \{A, B\}$, $0 \le s_j < \infty$, then*

$$t \mapsto P^{Y_t}\{\mathcal{E}\} \qquad \textit{is right-continuous a.s. } [P^\sigma]. \tag{2.24}$$

*The process $\{Y_t\}$ starting at $\sigma \in \Sigma_0$ has the strong Markov property with respect to the filtration $\{\mathcal{F}_t\}_{t \ge 0}$, where*

$$\mathcal{F}_t := \bigcap_{h > 0} \mathcal{F}^0_{t+h}. \tag{2.25}$$

We have claimed that the $Y$-process on $\Sigma_0$ is a nice Markov process, but before we can accept it as a version of a process as described in the abstract we have to show that $\Sigma_0$ is not empty. In the next proposition we shall show the even stronger property that $\sigma$ lies in $\Sigma_0$ a.s. if $\sigma$ is chosen by putting $N_A(z,0-)$ $A$-particles at $z$, with the $N_A(z,0-)$ i.i.d. mean-$\mu_A$ Poisson variables, and by adding in total a finite number of $B$-particles. *From now on $P$ without superscript will be used for the measure governing the $Y$-process with such an initial measure.* This notation does not indicate the value of $\mu_A$, nor the location of the $B$-particles introduced at time 0, but these quantities have no significant influence anyway. Expectation with respect to $P$ will be denoted by $E$ without superscript. Note that the description of our system in the abstract forces all particles at any given space–time point to be of the same type. Thus if we put $B$-particles at $z_1,\dots,z_k$ at time 0, then we instantaneously have to change the $A$-particles there into $B$-particles.

The proof of the next proposition is basically a Peierls argument. We associate to each $B$-particle present at time $t$ and with a switching time before $\tau_\infty$ a different "genealogical path" which describes how the $B$-particle arose from the $B$-particles at time 0 by various coincidences between $A$- and $B$-particles, and then more or less count all the genealogical paths to show that the expected number of genealogical paths at each time $t<\infty$ is finite.

PROPOSITION 4. *For any choice of the location of the finite number of initial $B$-particles we have $\sigma\in\Sigma_0$ a.s. $[P]$. Equivalently*

$$\int P\{Y_0\in d\sigma\}P^\sigma\{\min\{\widehat{\tau},\tau_\infty\}=\infty\}=1. \tag{2.26}$$

PROOF. It is a trivial calculation to show that $EM_t<\infty$. We also note the following simple properties of the $\alpha_t$: for $z,z'\in\mathbb{Z}^d$

$$\alpha_{t+u}(z)\ge e^{-D_Au}\alpha_t(z) \tag{2.27}$$

and

$$\alpha_{t+s}(z)\ge\alpha_t(z')\alpha_s(z-z'). \tag{2.28}$$

We now claim that (2.20) holds a.s. $[P]$ on $\{t<\min\{\widehat{\tau},\tau_\infty\}\}$. [Note that Lemma 2 claims that (2.20) holds a.s. with respect to another measure, so that we cannot simply deduce our claim from Lemma 2.] To prove our claim we note that there are only finitely many $B$-particles in the system at any time $t<\min\{\widehat{\tau},\tau_\infty\}$, by the definitions of $\widehat{\tau}$ and $\tau_\infty$. It therefore suffices to prove (2.20) with the number of $A$-particles which visit instead of the number of all particles which visit. In turn, by virtue of (2.16), it suffices to show that for each fixed $(z,t)$

$$E\{(\text{number of particles in } \mathcal{P}^* \text{ which visit } z \text{ during } [0,t])\} < \infty \quad \text{a.s. } [P]. \tag{2.29}$$

To see that (2.29) indeed holds we note that, by a decomposition with respect to the starting point of the particles,

$$\begin{aligned} &E\{(\text{number of particles in } \mathcal{P}^* \text{ which visit } z \text{ during } [0,t])\} \\ &\qquad \le \sum_{y \in \mathbb{Z}^d} \mu_A P\{y + \pi_A(s,\rho) = z \text{ for some } s \le t\}. \end{aligned} \tag{2.30}$$

But, if $\rho$ starts at $y$ with type $A$, then

$$\begin{aligned} &\int_0^{t+1} \alpha_s(y-z)\,ds \\ &\qquad = E\{\text{amount of time spent by } \rho \text{ at } z \text{ during } [0,t+1] \text{ in } \mathcal{P}^*\} \\ &\qquad \ge P\{y + \pi_A(s,\rho) \text{ reaches } z \text{ at some } s \le t \\ &\qquad\qquad \text{and stays at } z \text{ for at least one unit of time}\} \\ &\qquad \ge e^{-D_A} P\{y + \pi_A(s,\rho) = z \text{ for some } s \le t\}. \end{aligned} \tag{2.31}$$

Thus (2.30), combined with (2.27) and (2.28), shows that

$$\begin{aligned} &E\{(\text{number of paricles in } \mathcal{P}^* \text{ which visit } z \text{ during } [0,t])\} \\ &\qquad \le e^{D_A} \sum_{y \in \mathbb{Z}^d} \mu_A \int_0^{t+1} \alpha_s(y-z)\,ds \\ &\qquad \le e^{D_A} \int_0^{t+1} e^{(t+1)D_A} \sum_{y \in \mathbb{Z}^d} \mu_A \alpha_{t+1}(y-z)\,ds \\ &\qquad \le e^{(t+2)D_A} [\alpha_1(z)]^{-1} \int_0^{t+1} \sum_{y \in \mathbb{Z}^d} \mu_A \alpha_{t+2}(y)\,ds \\ &\qquad = (t+1) e^{(t+2)D_A} [\alpha_1(z)]^{-1} E M_{t+2} < \infty. \end{aligned} \tag{2.32}$$

Of course (2.29) and our claim for (2.20) follow from this.

Now the fact that (2.20) holds a.s. $[P]$ tells us that a.s. $[P]$ only finitely many $B$-particles are created at any $\tau_k$ and therefore, a.s. $[P], \tau_{k+1} > \tau_k$ for all $k$ and one cannot have $\widehat{\tau} < \tau_\infty$. Therefore, to prove (2.26), we only have to prove that

$$P\{\tau_\infty = \infty\} = \int P\{Y_0 \in d\sigma\} P^\sigma\{\tau_\infty = \infty\} = 1. \tag{2.33}$$

At any time $t$, we shall associate to any particle $\rho$ that has type $B$, and turned into a $B$-particle strictly before $\tau_\infty$, a unique genealogical path. One component of the genealogical path is a space–time path $\widetilde{\zeta}(\cdot,\rho)$ on $[0,t]$ which keeps track of the space–time paths traversed by the "ancestors" of $\rho$, that is, by $B$-particles which "transmitted" the $B$-type to $\rho$. The genealogical path also contains some additional information about the identity of the transmitting particles. The space–time part $\widetilde{\zeta}$ of the genealogical path is constructed as follows. If $\rho$ starts at $z$ and has type $B$ already at time 0, then its genealogical path is just the space–time path followed by $\rho$, restricted to $[0,t]$, that is, $\widetilde{\zeta}(s,\rho)=\pi(s,\rho)=\pi(0,\rho)+\pi_B(s,\rho)=z+\pi_B(s,\rho), 0\le s\le t$. If $\rho$ initially has type $A$, then $\rho$ first turned into a $B$-particle at its switching time $\theta(\rho)$, which necessarily is less than or equal to $t$. Moreover, we assumed $\rho$ became of type $B$ before $\tau_\infty$, that is, $\theta(\rho)<\tau_\infty$. The path component of the genealogical path of $\rho$ will then be $\widetilde{\zeta}(s,\rho)=\pi(s,\rho)$ for $\theta(\rho)\le s\le t$. Note $\theta(\rho)=\tau_k$ for some $k$. At this time either the particle $\rho$ jumped from some site $y$ to a site which contained some $B$-particle $\rho'$, or there was a $B$-particle $\rho'$ which jumped from some position $y$ onto the position of $\rho$ at time $\theta(\rho)$. In the former case $\rho'$ may not be unique, but we make some choice for $\rho'$ among the $B$-particles at the site to which $\rho$ jumps at time $\theta(\rho)$. We now follow the particle $\rho'$ backward in time till time $\theta(\rho')$ when it first turned into a $B$-particle, and take $\widetilde{\zeta}(s,\rho)=\pi(s,\rho')$ for $\theta(\rho')\le s<\theta(\rho)$. If $\theta(\rho')=0$, then we have defined the genealogical path of $\rho$ on the whole interval $[0,t]$ and we are done. If $\theta(\rho')>0$, then $\theta(\rho')=\tau_{k'}$ for some $k'<k$ and $\rho'$ coincided with some other $B$-particle $\rho''$ at $\theta(\rho')$. We then follow $\rho''$ backward in time, and so on, till we arrive at time 0. If we now traverse $\widetilde{\zeta}$ in the natural direction from 0 to $t$, then we see that this path starts with following the path of some initial $B$-particle, $\rho_{j_0}$, till some time $s_1$. At time $s_1$ either $\rho_{j_0}$ jumps to a point where there is an $A$-particle $\rho_{j_1}$, or some $A$-particle $\rho_{j_1}$ jumps at time $s_1$ to the position of $\rho_{j_0}$ at time $s_1$. Thus $\rho_{j_1}$ turns into a $B$-particle at time $s_1$, so that $s_1=\theta(\rho_{j_1})$. The path $\widetilde{\zeta}$ then follows the path of $\rho_{j_1}$ till some time $s_2$, at which $\rho_{j_1}$ coincides with an $A$-particle $\rho_{j_2}$, which turns into a $B$-particle due to this coincidence. This continues until some time $s_\ell$, at which the $A$-particle $\rho$ turns into a $B$-particle. Thus $s_\ell$ equals what we called $\theta(\rho)$ before. $\widetilde{\zeta}$ then equals $\pi(\cdot,\rho)$ on $[s_\ell,t]$. We shall take $\rho_{j_\ell}$ equal to $\rho$.

We shall want to keep track of some further data in the genealogical path. It will be convenient at this stage to label the initial $A$-particles by their initial position and their number in some arbitrary ordering of the initial particles at that site. Thus $\langle z,m\rangle$ will be used to denote the $m$th particle which started at $z$. We shall say that the particle $\langle z,m\rangle$ *exists* if and only if there are at least $m$ particles at $z$ at time 0. The particle $\rho_{j_i}$ appearing in the genealogical path in the preceding paragraph will also be denoted by $\langle z_i,m_i\rangle$. We denote by $\eta_i$ the type of the particle which jumps at time $s_i$. This particle can be $\rho_{j_{i-1}}$ or $\rho_{j_i}$. $\eta_i$ takes one of the values $A$ or $B$. We further denote, for $1\le i\le\ell$, by $y_i$ the position from which the particle jumps at time $s_i$, and by $x_i$ the position to which the particle jumps at time $s_i$. $x_{\ell+1}$ will be the position of $\rho = \rho_{j_\ell}$ at time $t$. The full genealogical path associated to $\rho$ now consists of $\widetilde{\zeta}(\cdot, \rho)$ plus the $(x_i, s_i, \eta_i, y_i, \rho_{j_0}, \rho_{j_i} = \langle z_i, m_i \rangle)$, that is, the positions and times of the jumps at which there is a changeover from one particle to another, as well as which particles jump at these times and which particles continue along $\widetilde{\zeta}(\cdot, \langle z_i, m_i \rangle)$. We obtain the genealogical paths of all $B$-particles at time $t$ by this forward construction and taking all possible values of $\widetilde{\zeta}, \ell$ and $(x_i, s_i, \eta_i, y_i, \rho_{j_0}, \rho_{j_i}), 1 \le i \le \ell$. Since each genealogical path is the genealogical path of just one particle, namely $\rho_{j_\ell}$, the number of $B$-particles at time $t$ is at most equal to the number of genealogical paths obtained in this forward construction. A crucial observation is that the $\rho_{j_i}, 0 \le i \le \ell$, *have to be distinct.* Indeed, in the construction of the genealogical path, $\rho_{j_i}$ is a particle whose type changes from $A$ to $B$ at time $s_i$ (with $s_0 = 0$), and any particle can change from type $A$ to type $B$ only once. We also note that $\rho_{j_i}$ becomes a $B$-particle at time $s_i$ and then must move from $x_i$ to $x_{i+1}$ during $[s_i, s_{i+1})$ if $\eta_{i+1} = A$, or must move from $x_i$ to $y_{i+1}$ during $[s_i, s_{i+1})$ if $\eta_{i+1} = B$.

We claim that it suffices for (2.33) to prove that for any $t$

$$E\{\text{total number of genealogical paths defined on } [0,t]\} < \infty. \tag{2.34}$$

Indeed, this will imply that the number of $B$-particles which arises before $t$ is almost surely finite. Since infinitely many $B$-particles have been created by time $\tau_\infty$ this will also give

$$P\{\tau_\infty < t\} = 0 \qquad \text{for any } t. \tag{2.35}$$

Thus (2.34) is indeed sufficient for (2.33). For the time being we shall estimate

$$\begin{aligned} &E\{\text{number of genealogical paths associated to} \\ &\qquad \text{some } B\text{-particle that is in the set } E \text{ at time } t\}, \end{aligned} \tag{2.36}$$

for any subset $E$ of $\mathbb{Z}^d$. Only near the end of this proof shall we take $E = \mathbb{Z}^d$ to get (2.33).

We bound the expectation in (2.36) by decomposing with respect to $\ell$ and the data $(x_i, s_i, \eta_i, y_i, \rho_{j_0}, \rho_{j_i} = \langle z_i, m_i \rangle), 1 \le i \le \ell$. Of course we cannot directly decompose with respect to the $s_i$, but have to follow the usual procedure which specifies only that the jump occurs in some interval $J(k) = J_n(k) := (k/n, (k+1)/n]$ and then let $n$ go to infinity. To this end we introduce the following indicator functions (with $\mathbf{z}$ used as an abbreviation for an $\ell$-tuple $z_1, \dots, z_\ell$, and similarly for $\mathbf{m}, \mathbf{k}$): if $\eta_i = A$ and $2 \le i \le \ell$, then

$$\begin{aligned} &I_{i,A}(\mathbf{k}, \mathbf{z}, \mathbf{m}) \\ &\quad = I[\langle z_{i-1}, m_{i-1} \rangle \text{ is at } x_{i-1} \text{ at time } (k_{i-1}+1)/n \end{aligned}$$

$$\begin{aligned}\text{(2.37)} \qquad & \text{and moves from there to } x_i \text{ during } [(k_{i-1}+1)/n, k_i/n)] \\ & \times I[\langle z_i, m_i\rangle \text{ is an } A\text{-particle at } y_i \text{ at time } (k_i/n)-, \\ & \qquad \text{it jumps to } x_i \text{ during } J(k_i) \text{ and becomes the } B\text{-particle } \rho_{j_{i+1}}],\end{aligned}$$

whereas for $\eta_i = B, 2 \le i \le \ell$,

$$\begin{aligned} & I_{i,B}(\mathbf{k},\mathbf{z},\mathbf{m}) \\ & \quad = I[\langle z_{i-1}, m_{i-1}\rangle \text{ is at } x_{i-1} \text{ at time } (k_{i-1}+1)/n \\ \text{(2.38)} \qquad & \qquad \text{and moves from there to } y_i \text{ during } [(k_{i-1}+1)/n, k_i/n)] \\ & \qquad \times I[\langle z_i, m_i\rangle \text{ is an } A\text{-particle at } x_i \text{ at time } (k_i/n)- \text{ and} \\ & \qquad\quad \text{during } J(k_i) \text{ the } B\text{-particle } \langle z_{i-1}, m_{i-1}\rangle \text{ jumps from } y_i \text{ to } x_i];\end{aligned}$$

$$\text{(2.39)} \qquad I_{\ell+1} = I[\rho_{j_\ell} \text{ moves to a position in } E \text{ during } [(k_{\ell+1})/n, t]].$$

For $i = 1$ the definitions of $I_{1,A}$ and $I_{1,B}$ need small changes, which amount to interpreting $\langle z_0, m_0\rangle$ as $\rho_{j_0}, k_0$ as $-1$ and $x_0$ as the initial position of $\rho_{j_0}$. For instance, $I_{1,B}$ is defined as

$$\begin{aligned} & I_{1,B}(\mathbf{k},\mathbf{z},\mathbf{m}) \\ & \quad = I[\rho_{j_0} \text{ moves from } x_0 \text{ to } y_1 \text{ during } [0, k_1/n)] \\ \text{(2.40)} \qquad & \qquad \times I[\langle z_1, m_1\rangle \text{ is an } A\text{-particle at } x_1 \text{ at time } (k_1/n)- \text{ and} \\ & \qquad\quad \text{during } J(k_i) \text{ the } B\text{-particle } \rho_{j_0} \text{ jumps from } y_1 \text{ to } x_1].\end{aligned}$$

We leave the corresponding definition of $I_{1,A}$ to the reader. Finally we define

$$\begin{aligned} & H_n = H_n(\mathbf{k},\mathbf{z},\mathbf{m}) \\ \text{(2.41)} \qquad & \quad = I[\text{the particles } \rho_{j_q}, 1 \le q \le \ell, \text{ together} \\ & \qquad \text{have at most one jump during } J_n(k_i), 1 \le i \le \ell].\end{aligned}$$

We shall use $\prod^{(\eta)}$ to denote the product over the indices $i \in [1, \ell]$ with $\eta_i = \eta$. Also $\sum^{(\ell)}$ is the sum over all ordered $\ell$-tuples $\rho_{j_1}, \dots, \rho_{j_\ell}$ of initial $A$-particles which are distinct, and distinct from $\rho_{j_0}$. Finally, $\sum_{x_1,\dots,x_{\ell+1}}$ will be short for the sum over $x_1, \dots, x_\ell \in \mathbb{Z}^d$, and over $x_{\ell+1} \in E$.

We claim that for fixed $\rho_{j_0}$ and $\mathbf{z},\mathbf{m}$ there are almost surely no common jump times in the paths of any pair of particles from $\{\rho_{j_0}, \rho_{j_i} = \langle z_i, m_i\rangle : 1 \le i \le \ell\}$. This follows from the fact that for any $\rho$, $\pi(\cdot,\rho)$ can have a jump at a time $s$ only if $\pi_A(\cdot,\rho)$ or $\pi_B(\cdot,\rho)$ has a jump at $s$ [by virtue of (2.2), (2.4)]. Our claim then follows because all the pairs of paths $\pi_A(\cdot,\rho_{j_i}), \pi_B(\cdot,\rho_{j_i})$ for different $i$ are independent. It follows that, for given $\mathbf{z}, \mathbf{m}$, we have almost surely $\inf_{k_i \le nt} H_n(\mathbf{k}, \mathbf{z}, \mathbf{m}) \to 1$ as $n \to \infty$. Consequently,

$$\sum_{\rho_0} \sum_{x_1,\dots,x_{\ell+1}} \sum_{\eta_1,\dots,\eta_\ell} \sum_{y_1,\dots,y_\ell} \sum_{0<k_1<\cdots<k_\ell<nt} H_n I_{\ell+1} \times \prod{}^{(A)} I_{i,A}(\mathbf{k}, \mathbf{z}, \mathbf{m}) \prod{}^{(B)} I_{i,B}(\mathbf{k}, \mathbf{z}, \mathbf{m})$$

almost surely converges as $n \to \infty$ to the number of genealogical paths with exactly $\ell$ changeovers among given pairs $\langle z_{i-1}, m_{i-1} \rangle$ and $\langle z_i, m_i \rangle$, and with final position in $E$ at time $t$. Thus, by Fatou's lemma,

$$\begin{aligned} &E\{\text{number of genealogical paths with } \ell \text{ changeover} \\ &\qquad \text{times associated to some } B\text{-particle which is in } E \text{ at time } t\} \\ &\quad \le \liminf_{n\to\infty} E\Bigg\{ \sum_{\rho_{j_0}} \sum_{x_1,\dots,x_{\ell+1}} \sum_{\eta_1,\dots,\eta_\ell} \sum_{y_1,\dots,y_\ell} \sum_{0<k_1<\cdots<k_\ell<nt} \sum{}^{(\ell)} H_n I_{\ell+1} \\ &\qquad \times \prod{}^{(A)} I_{i,A}(\mathbf{k}, \mathbf{z}, \mathbf{m}) \prod{}^{(B)} I_{i,B}(\mathbf{k}, \mathbf{z}, \mathbf{m}) \Bigg\}. \end{aligned} \tag{2.42}$$

Note also that our particles perform simple random walks, so that the sum over $y_i$ in (2.42) can be restricted to the neighbors of $x_i$. The sum over $\rho_{j_0}$ runs over the finite number of initial $B$-particles. For simplicity we restrict ourselves in the remainder of this proof to the case in which there is only one initial $B$-particle, and that it starts from position $x_0$. We can then drop the sum over $\rho_{j_0}$.

We wish to establish some independence between the required jumps and the required movement of $B$-particles in the indicator functions in the right-hand side of (2.42). For this we shall again make use of the particle system $\mathcal{P}^*$, which we coupled to our true particle system just before Lemma 2. We shall use $\mathcal{P}_0$ to denote the true particle system and use $N_A(x,t)$ for the number of $A$ particles at the space–time point $(x,t)$ in this true system. Recall that we coupled $\mathcal{P}^*$ and $\mathcal{P}_0$ in such a way that $N^*(x,0) = N_A(x,0-)$ for all $x$. Thus in the present situation the $N^*(x,0)$ are i.i.d. mean-$\mu_A$ Poisson variables. According to our construction $N_A(x,t) = 0$ for $x \in \mathbb{Z}^d, t \ge \tau_\infty$, and

$$N_A(x,t) \le N^*(x,t) \qquad \text{for all } x \in \mathbb{Z}^d, t \ge 0. \tag{2.43}$$

It also follows that if $\eta_i = A$, then

$$\begin{aligned} &I_{i,A}(\mathbf{k}, \mathbf{z}, \mathbf{m}) \\ &\quad \le I[\text{in } \mathcal{P}_0, \langle z_{i-1}, m_{i-1} \rangle \text{ moves} \\ &\qquad \text{from } x_{i-1} \text{ to } x_i \text{ during } [(k_{i-1}+1)/n, k_i/n)] \\ &\qquad \times I[\text{in } \mathcal{P}^*, \langle z_i, m_i \rangle \text{ is at } y_i \text{ at time } (k_i/n)- \text{ and} \end{aligned}$$

$$
\begin{aligned}
&\qquad\qquad\text{jumps to } x_i \text{ during } J(k_i)] \\
&\le I[\langle z_i, m_i\rangle \text{ exists}] I[\pi_A(k_i/n, \langle z_i, m_i\rangle) = y_i - z_i \\
&\qquad\qquad \text{and } \langle z_i, m_i\rangle \text{ jumps to } x_i \text{ during } J(k_i)] \\
&\quad \times I[\pi_B(k_i/n, \langle z_{i-1}, m_{i-1}\rangle) \\
&\qquad - \pi_B((k_{i-1}+1)/n, \langle z_{i-1}, m_{i-1}\rangle) = x_i - x_{i-1}] \\
&=: K_{i,A}(\mathbf{k}, \mathbf{z}, \mathbf{m}) L_{i,A}(\mathbf{k}, \mathbf{z}, \mathbf{m})
\end{aligned}
\tag{2.44}
$$

with $K_{i,A}$ standing for $I[\langle z_i, m_i\rangle \text{ exists}]$ times the indicator function involving $\pi_A(\cdot, \langle z_i, m_i\rangle)$, while $L_{i,A}$ stands for the indicator function involving $\pi_B(\cdot, \langle z_{i-1}, m_{i-1}\rangle)$ in the right-hand side. For $i = 1$ we interpret $k_0$ as $-1$. Similarly, if $\eta_i = B$, then

$$
\begin{aligned}
&I_{i,B}(\mathbf{k}, \mathbf{z}, \mathbf{m}) \\
&\quad \le I[\langle z_i, m_i\rangle \text{ exists}] I[\pi_A(k_i/n, \langle z_i, m_i\rangle) = x_i - z_i] \\
&\qquad \times I[\pi_B(k_i/n, \langle z_{i-1}, m_{i-1}\rangle) \\
&\qquad\quad - \pi_B((k_{i-1}+1)/n, \langle z_{i-1}, m_{i-1}\rangle) = y_i - x_{i-1} \\
&\qquad\qquad\qquad \text{and jumps to } x_i - x_{i-1} \text{ during } J(k_i)] \\
&\quad =: K_{i,B}(\mathbf{k}, \mathbf{z}, \mathbf{m}) L_{i,B}(\mathbf{k}, \mathbf{z}, \mathbf{m}).
\end{aligned}
\tag{2.45}
$$

For $i = 1$ we again take $k_0 = -1$. Finally,

$$
\begin{aligned}
I_{\ell+1} \le I[&\pi_B(t, \langle z_\ell, m_\ell\rangle) \\
&- \pi_B((k_{\ell+1})/n, \langle z_\ell, m_\ell\rangle) \in E - x_\ell] =: L_{\ell+1}.
\end{aligned}
\tag{2.46}
$$

We may therefore replace $I_{i,A}, I_{i,B}, I_{\ell+1}$ in the right-hand side of (2.42) by the appropriate right-hand sides in (2.44)–(2.46). Consider now an $i$ with $\eta_i = \eta_{i+1} = A$. In this case $\langle z_i, m_i\rangle$ has to exist and to move from $z_i$ to $x_i$ during $[0, (k_i + 1)/n]$ following $\pi_A(\cdot, \langle z_i, m_i\rangle)$, and then from $x_i$ to $x_{i+1}$ during $((k_i + 1)/n, k_{i+1}/n]$ following $\pi_B(\cdot, \langle z_i, m_i\rangle)$. These are requirements on the increments of $\pi_A(\cdot, \langle z_i, m_i\rangle)$ and $\pi_B(\cdot, \langle z_i, m_i\rangle)$ during disjoint time intervals and are therefore independent. The first requirement appears in one of the factors $K_{i,A}$, while the second requirement occurs in one of the factors $L_{i,A}$. A similar situation prevails for the other three possible values of the pair $(\eta_i, \eta_{i+1})$. Because the paths $\pi_A, \pi_B$ have independent increments, the requirements which appear in a $K$-factor and in an $L$-factor are independent for each particle separately. Since further the pairs $(\pi_A, \pi_B)$ for different particles are completely independent, we find that the left-hand side of (2.42) is bounded by

$$\begin{aligned}
&\liminf_{n\to\infty} \sum_{x_1,\dots,x_{\ell+1}} \sum_{\eta_1,\dots,\eta_\ell} \sum_{y_1,\dots,y_\ell} \sum_{0<k_1<\cdots<k_\ell<nt} \sum\nolimits^{(\ell)} \\
&\qquad E\Big\{H_n^A \prod\nolimits^{(A)} K_{i,A}(\mathbf{k},\mathbf{z},\mathbf{m}) \prod\nolimits^{(B)} K_{i,B}(\mathbf{k},\mathbf{z},\mathbf{m})\Big\} \\
&\qquad \times E\Big\{H_n^B L_{\ell+1} \prod\nolimits^{(A)} L_{i,A}(\mathbf{k},\mathbf{z},\mathbf{m}) \prod\nolimits^{(B)} L_{i,B}(\mathbf{k},\mathbf{z},\mathbf{m})\Big\},
\end{aligned} \tag{2.47}$$

where

$$H_n^\eta = I[\text{each particle } \langle z_j, m_j\rangle \text{ with } \eta_j = \eta \text{ has exactly one jump in } J(k_j)].$$

With $S^\eta$ as in the beginning of this section we can write

$$\begin{aligned}
&E\Big\{H_n^B L_{\ell+1} \prod\nolimits^{(A)} L_{i,A}(\mathbf{k},\mathbf{z},\mathbf{m}) \prod\nolimits^{(B)} L_{i,B}(\mathbf{k},\mathbf{z},\mathbf{m})\Big\} \\
&\quad \le \prod\nolimits^{(A)} P\{S^B_{(k_i-k_{i-1}-1)/n} = x_i - x_{i-1}\} \\
&\qquad \times \prod\nolimits^{(B)} P\{S^B_{(k_i-k_{i-1}-1)/n} = y_i - x_{i-1} \\
&\qquad\qquad \text{and has exactly one jump during } J(k_i - k_{i-1} - 1) \\
&\qquad\qquad\qquad \text{and this goes from } y_i - x_{i-1} \text{ to } x_i - x_{i-1}\} \\
&\qquad \times P\{S^B_{t-(k_\ell+1)/n} \in E - x_\ell\}.
\end{aligned} \tag{2.48}$$

To simplify our formulae somewhat we now use that, as in (2.27),

$$\begin{aligned}
&P\{S^B_{(k_i-k_{i-1})/n} = x_i - x_{i-1}\} \\
&\quad \ge P\{S^B_{(k_i-k_{i-1}-1)/n} = x_i - x_{i-1}\} \\
&\qquad \times P\{S^B_\cdot \text{ remains constant during } [(k_i-k_{i-1}-1)/n, (k_i-k_{i-1})/n]\} \\
&\quad = e^{-D_B/n} P\{S^B_{(k_i-k_{i-1}-1)/n} = x_i - x_{i-1}\}.
\end{aligned}$$

We write $\nu = \nu(\eta) = \nu(\eta, \ell)$ for the number of $1 \le i \le \ell$ with $\eta_i = A$. Then the last inequality combined with (2.48) shows that

$$\begin{aligned}
&E\Big\{H_n^B L_{\ell+1} \prod\nolimits^{(A)} L_{i,A}(\mathbf{k},\mathbf{z},\mathbf{m}) \prod\nolimits^{(B)} L_{i,B}(\mathbf{k},\mathbf{z},\mathbf{m})\Big\} \\
&\quad \le e^{\nu D_B/n} P\{S^B_t \in E - x_0, S^B_{k_i/n} = x_i - x_0 \text{ for } \eta_i = A; \\
&\qquad\qquad S^B_{(k_i-1)/n} = y_i - x_0 \text{ and } S^B_\cdot \text{ jumps} \\
&\qquad\qquad \text{from } y_i - x_0 \text{ to } x_i - x_0 \text{ during } J(k_i - 1) \text{ for } \eta_i = B\}.
\end{aligned} \tag{2.49}$$

The right-hand side is independent of $\mathbf{z}, \mathbf{m}$, so the expectation of the $L$ factors in the right-hand side of (2.47) can be replaced by (2.49) and taken outside the sum $\sum^{(\ell)}$.

Next we deal with $\sum^{(\ell)}$ of the expectation of the factors $K_i$. We claim that

$$\sum_{m_1,\dots,m_\ell} E\Big\{H_n^A \prod{}^{(A)} K_{i,A}(\mathbf{k},\mathbf{z},\mathbf{m}) \prod{}^{(B)} K_{i,B}(\mathbf{k},\mathbf{z},\mathbf{m})\Big\}$$

$$\le [\mu_A]^\ell \left[\frac{D_A}{2dn}\right]^\nu \prod{}^{(A)} P\{S^A_{k_i/n} = y_i - z_i\} \times \prod{}^{(B)} P\{S^A_{k_i/n} = x_i - z_i\}, \tag{2.50}$$

provided $\sum_{m_1,\dots,m_\ell}$ runs only over those $\ell$-tuples with $m_i \ge 1$ for which $\langle z_i, m_i\rangle$, $1 \le i \le \ell$, are distinct. To prove this, first fix $\mathbf{z}, \mathbf{m}$ such that all $\langle z_i, m_i\rangle, 1 \le i \le \ell$, are distinct $A$-particles. For such $\ell$-tuples, the paths $\pi_A(\cdot, \langle z_i, m_i\rangle)$ are independent, and

$$P\{\pi_A(k_i/n, \langle z_i, m_i\rangle) = y_i - z_i \text{ and } \langle z_i, m_i\rangle \text{ jumps to } x_i \text{ during } J(k_i)\}$$
$$\le \frac{D_A}{2dn} P\{S^A_{k_i/n} = y_i - z_i\},$$

while

$$P\{\pi_A(k_i/n, \langle z_i, m_i\rangle) = x_i - z_i\} = P\{S^A_{k_i/n} = x_i - z_i\}.$$

Consequently, the left-hand side of (2.50) is at most

$$\sum_{m_1,\dots,m_\ell} E\left\{\prod_{i=1}^{\ell} I[\langle z_i, m_i\rangle \text{ exists}]\right\}$$
$$\times \left[\frac{D_A}{2dn}\right]^\nu \prod{}^{(A)} P\{S^A_{k_i/n} = y_i - z_i\} \prod{}^{(B)} P\{S^A_{k_i/n} = x_i - z_i\}.$$

Therefore it suffices for (2.50) to show that

$$E\left\{\sum_{m_1,\dots,m_\ell} \prod_{i=1}^{\ell} I[\langle z_i, m_i\rangle \text{ exists}]\right\} \le [\mu_A]^\ell. \tag{2.51}$$

To prove this last inequality, we partition the $z_i$ into maximal classes of equal $z$'s. More precisely, let $a_1, \dots, a_p \in \mathbb{Z}^d$ be distinct, and let $T_1, \dots, T_p$ be a partition of $\{1, \dots, \ell\}$ and let $z_i = a_j$ precisely for $i \in T_j$. Finally, let $T_j$ have exactly $q_j$ elements. If we write $[N]_k$ for $N(N-1)\cdots(N-k+1)$, then

$$\sum_{m_1,\dots,m_\ell} \prod_{i=1}^{\ell} I[\langle z_i, m_i\rangle \text{ exists}] = \prod_{j=1}^{p} [N_A(a_j, 0)]_{q_j}. \tag{2.52}$$

Inequality (2.51) now follows by taking the expectation in (2.52). [In fact, since we assumed that the $N_A$ have a Poisson distribution, (2.51) holds with equality. We point out here that (2.51) also holds if $N_A(z,0) \le \mu_A$ with probability 1, rather than distributed like a mean-$\mu_A$ Poisson variable.]

As pointed out, (2.51) proves (2.50). If we sum (2.50) over the $z_i$ and use (2.49), we obtain

$$\begin{aligned}
&\sum{}^{(\ell)} E\Big\{H_n^A \prod{}^{(A)} K_{i,A}(\mathbf{k},\mathbf{z},\mathbf{m}) \prod{}^{(B)} K_{i,B}(\mathbf{k},\mathbf{z},\mathbf{m})\Big\}\\
&\quad \times E\Big\{H_n^B L_{\ell+1} \prod{}^{(A)} L_{i,A}(\mathbf{k},\mathbf{z},\mathbf{m}) \prod{}^{(B)} L_{i,B}(\mathbf{k},\mathbf{z},\mathbf{m})\Big\}\\
&\quad \le [\mu_A]^\ell \left[\frac{D_A e^{D_B/n}}{2dn}\right]^\nu \\
&\qquad \times P\{S_t^B \in E - x_0, S^B_{k_i/n} = x_i - x_0 \text{ for } \eta_i = A;\\
&\qquad\quad S^B_{(k_i-1)/n} = y_i - x_0 \text{ and } S^B_\cdot \text{ jumps from } y_i - x_0 \text{ to } x_i - x_0\\
&\qquad\qquad \text{during } J(k_i - 1) \text{ for } \eta_i = B\}.
\end{aligned} \tag{2.53}$$

We now fix the set of indices for which $\eta_i = B$. Let this set be $\mathcal{D} = \{i_1 < i_2 < \cdots < i_\kappa\} \subset \{1, \dots, \ell\}$. We also fix the $k_{i_j}$ for $1 \le j \le \kappa$. Note that $\mathcal{D} = \varnothing$, or equivalently, $\kappa = 0$ is possible. Further set $i_0 = 0, i_{\kappa+1} = \ell + 1, k_0 = -1, k_{\ell+1} = \lfloor nt \rfloor$. Finally note that

$$\nu = \ell - \kappa = \sum_{j=0}^{\kappa} [i_{j+1} - i_j - 1], \tag{2.54}$$

and that for all integers $a \le b$, and $r \ge 0$

$$\sum_{a < k_{p+1} < k_{p+2} < \cdots < k_{p+r} \le b} 1 = \binom{b-a}{r} \le \frac{(b-a)^r}{r!} \tag{2.55}$$

(the sum here is over $k_{p+1}, \dots, k_{p+r}$). We now sum (2.53) first over all $x_i, y_i$ with $i \notin \mathcal{D}$. The sum of the right-hand side of (2.53) over these $x_i, y_i$ equals

$$\begin{aligned}
&[\mu_A]^\ell \left[\frac{D_A e^{D_B/n}}{n}\right]^\nu \\
&\quad \times P\{S_t^B \in E - x_0, S^B_{(k_i-1)/n} = y_i - x_0 \text{ and } S^B_\cdot \\
&\qquad \text{jumps from } y_i - x_0 \text{ to } x_i - x_0 \text{ during } J(k_i - 1) \text{ for } i \in \mathcal{D}\}.
\end{aligned} \tag{2.56}$$

Next sum over the $k_j$ with $j \ge 1$, but $j \notin \mathcal{D}$. By means of (2.55) we see that the sum over the $k_j$ with $k_{i_s} < k_j < k_{i_{s+1}}$ contributes a factor no larger than

$$\frac{(k_{i_{s+1}} - k_{i_s})^{i_{s+1} - i_s - 1}}{(i_{s+1} - i_s - 1)!}.$$

In this way we obtain that the contribution to (2.42) of the terms with $\eta_i = B$ exactly for $i \in \mathcal{D} = \{i_1 < i_2 < \cdots < i_\kappa\}$, with $\ell, \mathcal{D}$ and $k_{i_j}$ for $i_j \in \mathcal{D}$ fixed (before taking the liminf over $n$), is at most

$$\sum_{x_i, y_i, i \in \mathcal{D}} [\mu_A]^\kappa \left[\frac{D_A \mu_A e^{D_B/n}}{n}\right]^{\ell-\kappa} \frac{(k_{i_1})^{i_1-1}}{(i_1-1)!} \prod_{j=1}^{\kappa} \frac{(k_{i_{j+1}} - k_{i_j})^{i_{j+1}-i_j-1}}{(i_{j+1}-i_j-1)!}$$
$$\times P\{S_t^B \in E - x_0 \text{ and } S_\cdot^B \text{ jumps}$$
$$\text{from } y_i - x_0 \text{ to } x_i - x_0 \text{ during } J(k_i - 1) \text{ for } i \in \mathcal{D}\}$$
$$= [\mu_A]^\kappa \left[\frac{D_A \mu_A e^{D_B/n}}{n}\right]^{\ell-\kappa} \frac{(k_{i_1})^{i_1-1}}{(i_1-1)!} \prod_{j=1}^{\kappa} \frac{(k_{i_{j+1}} - k_{i_j})^{i_{j+1}-i_j-1}}{(i_{j+1}-i_j-1)!}$$
$$\times P\{S_t^B \in E - x_0 \text{ and } S_\cdot^B \text{ has}$$
$$\text{a jump during } J(k_i - 1) \text{ for } i \in \mathcal{D}\}$$

$$\begin{aligned} &= [\mu_A]^\kappa \frac{1}{(i_1-1)!}\left(\frac{D_A \mu_A k_{i_1} e^{D_B/n}}{n}\right)^{i_1-1} \\ &\quad \times \prod_{j=1}^{\kappa}\left[\frac{1}{(i_{j+1}-i_j-1)!}\left(\frac{D_A \mu_A (k_{i_{j+1}} - k_{i_j}) e^{D_B/n}}{n}\right)^{i_{j+1}-i_j-1}\right] \\ &\quad \times P\{S_t^B \in E - x_0 \text{ and } S_\cdot^B \text{ has} \\ &\qquad \text{a jump during } J(k_i - 1) \text{ for } i \in \mathcal{D}\}. \end{aligned} \tag{2.57}$$

We now sum (2.42) also over $\ell \geq i_\kappa$ and use Fatou's lemma to bring the lim inf outside the sum over $\ell$. We also rename $k_{i_j}$ as $r_j$. Since $i_{\kappa+1} = \ell + 1$ and $r_{\kappa+1} = \lfloor nt \rfloor$, this yields

$$E\{\text{number of genealogical paths associated}$$
$$\text{to some } B\text{-particle which is in } E \text{ at time } t\}$$

$$\begin{aligned} &\leq \liminf_{n\to\infty} \sum_{\kappa \geq 0} [\mu_A]^\kappa \\ &\quad \times \sum_{\mathcal{D} = \{i_1 < \cdots < i_\kappa\}} \sum_{0 < r_1 < \cdots < r_\kappa < nt} \exp\left[\frac{D_A \mu_A e^{D_B/n}}{n}(\lfloor nt \rfloor - r_\kappa)\right] \\ &\quad \times \frac{1}{(i_1-1)!}\left(\frac{D_A \mu_A r_1 e^{D_B/n}}{n}\right)^{i_1-1} \\ &\quad \times \prod_{j=1}^{\kappa-1}\left[\frac{1}{(i_{j+1}-i_j-1)!}\left(\frac{D_A \mu_A (r_{j+1} - r_j) e^{D_B/n}}{n}\right)^{i_{j+1}-i_j-1}\right] \\ &\quad \times P\{S_t^B \in E - x_0 \text{ and } S_\cdot^B \text{ has} \end{aligned} \tag{2.58}$$

$$\text{a jump during } J(r_j-1) \text{ for } 1\le j\le\kappa\}.$$

We next carry out the sum over $\mathcal{D}=\{i_1<\cdots<i_\kappa\}$. This transforms the right-hand side of (2.58) into

$$\begin{aligned}
&\liminf_{n\to\infty}\sum_{\kappa\ge 0}[\mu_A]^\kappa\sum_{0<r_1<\cdots<r_\kappa<nt}\exp\biggl[\frac{D_A\mu_A r_1 e^{D_B/n}}{n}\\
&\qquad\qquad+\sum_{j=1}^{\kappa}\frac{D_A\mu_A(r_{j+1}-r_j)e^{D_B/n}}{n}\biggr]\\
&\qquad\times P\{S_t^B\in E-x_0 \text{ and } S_\cdot^B \text{ has}\\
&\qquad\qquad \text{a jump during } J(r_j-1) \text{ for } 1\le j\le\kappa\}\\
&=\liminf_{n\to\infty}\sum_{\kappa\ge 0}[\mu_A]^\kappa\sum_{0<r_1<\cdots<r_\kappa<nt}\exp[D_A\mu_A t]\\
&\qquad\times P\{S_t^B\in E-x_0 \text{ and } S_\cdot^B \text{ has}\\
&\qquad\qquad \text{a jump during } J(r_j-1) \text{ for } 1\le j\le\kappa\}.
\end{aligned}\tag{2.59}$$

At this point we finally specialize to $E=\mathbb{Z}^d$. With this choice the right-hand side of (2.59) is at most

$$\begin{aligned}
&\exp[D_A\mu_A t]\liminf_{n\to\infty}\sum_{\kappa\ge 0}[\mu_A]^\kappa\sum_{0<r_1<\cdots<r_\kappa<nt}\biggl[\frac{D_B}{n}\biggr]^\kappa\\
&\le\exp[D_A\mu_A t]\liminf_{n\to\infty}\sum_{\kappa\ge 0}\frac{1}{\kappa!}\biggl[\frac{D_B\mu_A nt}{n}\biggr]^\kappa \qquad \text{[by (2.55)]}\\
&=\exp[(D_A+D_B)\mu_A t]<\infty.
\end{aligned}\tag{2.60}$$

This proves (2.34) and the proposition in the case when we start with one $B$-particle. If we start with $N_B$ $B$-particles at $x_{0,1},x_{0,2},\dots,x_{0,N_B}$, respectively, then we only have to replace the probability in (2.59) by

$$\begin{aligned}
&\sum_{m=1}^{N_B}P\{S_t^B\in E-x_{0,m} \text{ and } S_\cdot^B \text{ has}\\
&\qquad \text{a jump during } J(r_j-1) \text{ for } 1\le j\le\kappa\}.
\end{aligned}\tag{2.61}$$

(The $x_{0,m}$ do not have to be distinct here.) □

REMARK 2. A check of the proof shows that (2.51) is the only property of the initial distribution which is used. In particular, (2.33) also holds for any initial distribution for which $N_A(z,0)$ is a.s. bounded by a constant. A special case of this last situation is also treated in [11], Lemma 3.2. The proof also works for initial distributions which are stochastically below a Poisson distribution. Also, by the argument given in the next section for Theorem 1, we obtain by specializing to $E = [\mathcal{C}(C_1 t)]^c$, that (2.51) is sufficient to conclude that (1.3) holds.

## 3. A linear upper bound for $B(t)$.

In this section we give the

Proof of Theorem 1. The arguments preceding (2.42) show that it is enough to show that for $E =$ the complement of $\mathcal{C}(C_1 t)$, the left-hand side of (2.42) is bounded by $2N_B \exp(-t)$. In turn, it suffices to prove that the right-hand side of (2.59) [with the last factor replaced by (2.61)] is bounded by $2N_B \exp(-t)$ if we take $E = [\mathcal{C}(C_1 t)]^c$. In order to show this we split the sum over $\kappa$ in (2.59) into two pieces. The first sum is over $\kappa \geq K_1 t$ and the second over $\kappa < K_1 t$, where $K_1$ is chosen so large that the first piece is bounded by [cf. (2.60)]

$$\begin{aligned}
&\limsup_{n\to\infty} N_B \sum_{\kappa \geq K_1 t} \sum_{0<r_1<\cdots<r_\kappa<nt} [\mu_A]^\kappa \exp\left[\frac{D_A \mu_A r_1 e^{D_B/n}}{n}\right. \\
&\qquad\qquad \left. + \sum_{j=1}^{\kappa} \frac{D_A\mu_A(r_{j+1}-r_j)e^{D_B/n}}{n}\right] \qquad (3.1)\\
&\qquad \times P\{S^B_\cdot \text{ has a jump during } J(r_j - 1) \text{ for } 1 \leq j \leq \kappa\} \\
&\leq N_B \exp[D_A\mu_A t] \limsup_{n\to\infty} \sum_{\kappa\geq K_1 t} \frac{1}{\kappa!}\left[\frac{D_B \mu_A n t}{n}\right]^\kappa \\
&\leq N_B e^{-t} \qquad (\text{for } t \geq 1).
\end{aligned}$$

Note that this estimate is uniform in $C_1$.

To estimate the second sum, over $\kappa < K_1 t$ (for fixed $K_1$), we note that the increments of $S^B$ over disjoint intervals are independent. Thus the sum of the increments of $S^B$ over

$$[0, t] \setminus \bigcup_{j=1}^{\kappa} J(r_j - 1)$$

has the same distribution as

$$S^{B,0}_{r_1 - 1} + \sum_{j=1}^{\kappa-1} S^{B,j}_{(r_{j+1}-r_j-1)/n} + S^{B,\kappa+1}_{t-r_\kappa/n},$$

where the $S^{B,j}$ are independent copies of $S^B$. In turn, this sum has the same distribution as $S^B_{t-\kappa/n}$, and is independent of the increments of $S^B$ over the $J(r_j-1)$. In addition, given that $S^B$ has a jump in $J(r_j-1)$, the conditional distribution of $S^B_{r_j/n} - S^B_{(r_j-1)/n}$ is the distribution of

$$\sum_{m=1}^{\psi} Z_m,$$

where $Z_1, Z_2, \dots$ are independent random variables, each with the distribution of a generic jump of $S^B$, and $\psi$ is independent of the $Z_i$, and $\psi$ has the conditional distribution of a mean-$D_B/n$ Poisson variable, given that this variable is at least 1. In our case $P\{Z_i = \pm e_j\} = 1/(2d)$, so that $\|\sum_{m=1}^{\psi} Z_m\| \le \psi$, and conditionally on the event $\{S^B$ has a jump in $J(r_j - 1), 1 \le j \le \kappa\}$, $\|S^B_t\|$ is stochastically smaller than

$$\|S^B_{t-\kappa/n}\| + \psi_1 + \cdots + \psi_\kappa,$$

with the $\psi_i$ independent copies of $\psi$, which are also independent of $S^B$. It is now a standard large deviation estimate that for fixed $K_1$ and $x_{0,m}$, and sufficiently large $C_1$ (independent of the $x_{0,m}$, though), and all sufficiently large $t$ and $\kappa < K_1 t$

$$\begin{aligned} P\{&S^B_t \notin \mathcal{C}(C_1 t) - x_{0,m}, S^B_\cdot \text{ has a jump in } J(r_j-1), 1 \le j \le \kappa\} \\ &\le \left[\frac{D_B}{n}\right]^\kappa P\left\{\|S^B_{t-\kappa/n}\| + \sum_{1\le j<K_1 t} \psi_j \ge C_1 t/2\right\} \\ &\le \left[\frac{D_B}{n}\right]^\kappa \exp[-(D_A\mu_A + D_B\mu_A + 1)t]. \end{aligned}$$

We leave the details of this to the reader (cf. (2.40) in [7]). For such a choice of $C_1$ it follows that the sum of the terms with $\kappa < K_1 t$ in (2.59) [with the replacement of (2.61)] is at most

$$\begin{aligned} \liminf_{n\to\infty} &\sum_{1\le m\le N_B} \sum_{0\le\kappa<K_1 t} \sum_{0<r_1<\cdots<r_\kappa<nt} [\mu_A]^\kappa \exp[D_A\mu_A t] \\ &\times P\{S^B_t \notin \mathcal{C}(C_1 t) - x_{0,m} \\ &\qquad \text{and } S^B_\cdot \text{ has a jump during } J(r_j - 1) \text{ for } 1 \le j \le \kappa\} \\ &\le N_B \exp[D_A\mu_A t] \liminf_{n\to\infty} \sum_{0\le\kappa<K_1 t} \frac{1}{\kappa!}\left(\frac{D_B\mu_A n t}{n}\right)^\kappa \\ &\qquad \times \exp[-(D_A\mu_A + D_B\mu_A + 1)t] \\ &\le N_B e^{-t}. \end{aligned}$$

For $K_1, C_1$ as above and $E = [\mathcal{C}(C_1 t)]^c$ we find that the expectation (2.36) is bounded by $2N_B \exp(-t)$ for all large $t$, so that (1.3) holds. Equation (1.4) now follows from the Borel–Cantelli lemma and the fact that $B(t)$ is increasing in $t$. $\square$

**4. A linear lower bound for $B(t)$ when $D_A = D_B$.** In this section we shall prove Theorem 2. We remind the reader that $P$ without superscript stands the measure governing the $Y$-process when the initial $N_A(x, 0-)$ are i.i.d. mean-$\mu_A$ Poisson variables and a finite number of $B$-particles are added at time 0. Throughout this section we assume that the $A$- and $B$-particles perform random walks with the same distribution, that is,

$$D_A = D_B. \tag{4.1}$$

We shall write $D$ for the common value of $D_A$ and $D_B$. As explained in Section 2 we then take $\pi_A(\cdot, \rho) \equiv \pi_B(\cdot, \rho)$. We then have that the position at time $s$ of a particle $\rho$ which starts in $z$ is $z + \pi_A(s, \rho)$ for all $s$. However, the type of $\rho$ will change from $A$ to $B$ at $\theta(\rho)$, the first instant when $\rho$ coincides with a $B$-particle. These paths $\pi_A(\cdot, \rho)$ for different $\rho$ are independent and so, as far as the positions of the particles are concerned, there is no interaction. Thus the system of particles which start out as $A$-particles (i.e., all particles but the finitely many initial $B$-particles) is the same as the system $\mathcal{P}^*$ described in Section 2 just before (2.16), as far as positions of particles are concerned. In agreement with Section 2 we write $N^*(x, t)$ for the number of particles at the space–time point $(x, t)$ which started out as an $A$-particle (but whose type may have changed to $B$ by time $t$). Of course, this assumes, as before, that $\mathcal{P}^*$ is coupled with the true system such that $N^*(x, 0) = N_A(x, 0-)$. In the present setup this means that the $N^*(x, 0), x \in \mathbb{Z}^d$, are i.i.d., mean-$\mu_A$ Poisson variables. The system $\mathcal{P}^*$ is then stationary in time for $t \geq 0$. It is convenient to extend the system $\mathcal{P}^*$ to a stationary system defined for all times $t \in \mathbb{R}$, including negative ones. For our system of noninteracting random walkers this can easily be done by extending the path $t \mapsto \pi_A(t, \rho)$ for each particle $\rho$ present at time zero to all $t$ in such a way that $\{-\pi_A(-t, \rho)\}_{t \geq 0}$ has the same distribution as $\{\pi_A(t, \rho)\}_{t \geq 0}$, and in such a way that the paths $\{\pi_A(\cdot, \rho)\}$, with $\rho$ varying over all particles present at time 0, are completely independent. We shall still use the notation $\mathcal{P}^*$ for the extended system. The configurations $\{N^*(x, t), x \in \mathbb{Z}^d\}$ are stationary in time in $\mathcal{P}^*$.

We remind the reader that we gave an outline of the proof of Theorem 2 under the heading "Some heuristics" in the Introduction. We now fill in the details of the proof. This proof is almost a "mirror image" of the proof of Theorem in [7]. The difference is that in [7] we wanted not too many particles near certain space–time paths, and here we want not too few particles near these space–time paths.

We repeat most of the definitions of the Introduction, but now in the form needed for a general dimension $d \geq 1$. The constants $C_0$ and $\gamma_r$ are chosen as follows: $\gamma_0 > 0$ is a constant which satisfies

$$0 < \gamma_0 \prod_{j=1}^{\infty} [1 - 2^{-j/4}]^{-1} \leq \tfrac{1}{2}. \tag{4.2}$$

We take

$$\gamma_1 = \gamma_0, \qquad \gamma_{r+1} = \gamma_0 \prod_{j=1}^{r} \left[1 - \frac{1}{C_0^{j/4}}\right]^{-1}, \qquad r > 0. \tag{4.3}$$

Further, $C_0 \geq 2$ is an integer which is so large that for all $r \geq 1$,

$$\begin{aligned} C_0^{-r/2} - \left(1 - \frac{C_4(r \log C_0)^d}{C_0^r}\right) \\ \times (1 - e^{-C_0^{-r/2}})[1 - C_0^{-r/4}]^{-1} \leq -\frac{1}{2} C_0^{-3r/4}, \end{aligned} \tag{4.4}$$

as well as

$$3^{d+1} C_0^{6(d+1)(r+1)} \exp[-\tfrac{1}{2}\gamma_0 \mu_A C_0^{(d-3/4)r}] \leq 1, \qquad r \geq 1. \tag{4.5}$$

Here $C_4$ is the constant of Lemma 5 below. Since $C_4$ will not depend on $C_0$, we can indeed fulfill (4.4) and (4.5) by taking $C_0$ large. Then (4.2), together with $C_0 \geq 2$, implies (1.11) with

$$\gamma_\infty := \lim_{r \to \infty} \gamma_r = \gamma_0 \prod_{j=1}^{\infty} \left[1 - \frac{1}{C_0^{j/4}}\right]^{-1}.$$

We take $\Delta_r = C_0^{6r}$ as before, and for $\mathbf{i} = (i(1), \ldots, i(d)) \in \mathbb{Z}^d$ we define (see Figure 1)

$$\mathcal{B}_r(\mathbf{i}, k) = \prod_{s=1}^{d} [i(s)\Delta_r, (i(s)+1)\Delta_r) \times [k\Delta_r, (k+1)\Delta_r)$$

and

$$\widetilde{\mathcal{B}}_r(\mathbf{i}, k) \risingdotseq \prod_{s=1}^{d} [(i(s)-3)\Delta_r, (i(s)+4)\Delta_r) \times [(k-1)\Delta_r, (k+1)\Delta_r).$$

For $x = (x(1), \ldots, x(d)) \in \mathbb{Z}^d$ we further take

$$U_r(x, v) = \sum_{y \in \mathcal{Q}_r(v)} N^*(y, v) \qquad \text{with } \mathcal{Q}_r(x) = \prod_{s=1}^{d} [x(s), x(s) + C_0^r). \tag{4.6}$$

Note that the edge size of the cube $\mathcal{Q}_r$ is only $C_0^r$ and not $\Delta_r$. Further we define

$$V_r(\mathbf{i}) = \prod_{s=1}^{d} [(i(s)-3)\Delta_r, (i(s)+4)\Delta_r),$$

and the *pedestal* of $\mathcal{B}_r(\mathbf{i},k)$ which is defined as

$$\mathcal{V}_r(\mathbf{i},k) = V_r(\mathbf{i}) \times \{(k-1)\Delta_r\} = \prod_{s=1}^{d} [(i(s)-3)\Delta_r, (i(s)+4)\Delta_r) \times \{(k-1)\Delta_r\}.$$

The $r$-block $\mathcal{B}_r(\mathbf{i},k)$ is called *bad* if

$$U_r(x,v) < \gamma_r \mu_A C_0^{dr} \qquad \text{for some } (x,v) \text{ with integer } v \text{ for which } \mathcal{Q}_r(x) \times \{v\} \text{ is contained in } \widetilde{\mathcal{B}}_r(\mathbf{i},k). \tag{4.7}$$

The pedestal $\mathcal{V}_r(\mathbf{i},k)$ of $\mathcal{B}_r(\mathbf{i},k)$ is called bad if

$$U_r(x,(k-1)\Delta_r) < \gamma_r \mu_A C_0^{dr} \qquad \text{for some } x \text{ with } \mathcal{Q}_r(x) \subset V_r(\mathbf{i}). \tag{4.8}$$

A block or pedestal is called *good* if it is not bad. Note that in contrast to [7], the good blocks and pedestals have $U(x,v)$ large.

Still more definitions are needed. As in Section 1 $\widehat{\pi}(\{s_i,x_i\})$ will denote the space–time path for which $\widehat{\pi}(s) = x_i$ for $s_i \le s < s_{i+1}$. Then exactly as in Section 1,

$$\Xi(\ell,t) := \{\widehat{\pi}(\{s_i,x_i\}_{0\le i\le \ell}) \text{ with } 0 = s_0 < s_1 < \cdots < s_\ell < t \text{ and } x_i \in \mathcal{C}(t\log t)\}, \tag{4.9}$$

$$\phi_r(\widehat{\pi}) = \text{number of bad } r\text{-blocks which intersect the space–time path } \widehat{\pi}, \tag{4.10}$$

$$\Phi_r(\ell) = \sup_{\widehat{\pi}\in\Xi(\ell,t)} \phi_r(\widehat{\pi}), \tag{4.11}$$

$$\psi_{r+1}(\widehat{\pi}) = \text{number of } (r+1)\text{-blocks which intersect the space–time path } \widehat{\pi} \text{ and which have a good pedestal but contain a bad } r\text{-block} \tag{4.12}$$

and

$$\Psi_r(\ell) = \sup_{\widehat{\pi}\in\Xi(\ell,t)} \psi_r(\widehat{\pi}) \tag{4.13}$$

(we suppress the dependence on $t$ in these quantities).

Exactly as in the argument for the one-dimensional case in the Introduction, or as in Lemma 8 of [7], we now have for any $\widehat{\pi} \in \Xi(\ell,t)$

$$\phi_r(\widehat{\pi}) \le C_0^{6(d+1)} \Phi_{r+1}(\ell) + C_0^{6(d+1)} \psi_{r+1}(\widehat{\pi}) \tag{4.14}$$

and

$$\Phi_r(\ell) \le C_0^{6(d+1)} \Phi_{r+1}(\ell) + C_0^{6(d+1)} \Psi_{r+1}(\ell). \tag{4.15}$$

Now choose some large constant $K$ and take $R = R(t)$ to be the integer for which

$$C_0^R \geq [K_4 \log t]^{1/d} > C_0^{R-1}. \tag{4.16}$$

(This differs slightly from [7], which had a 1 instead of the $K_4$ here.) As in [7], Lemmas 5 and 9, simple Poisson distribution estimates now show that we can take $K_4 = K_4(K, d, \mu_A)$ so large that

$$\begin{aligned} &P\{\Phi_r(\ell) > 0 \text{ for some } r \geq R \text{ and some } \ell \geq 0\} \\ &\quad \leq P\{\text{for some } r \geq R \text{ and } \ell \geq 0 \\ &\qquad\quad \text{a bad } r\text{-block intersects some } \widehat{\pi} \in \Xi(\ell, t)\} \\ &\quad \leq \sum_{r \geq R} P\{U_r(x, v) < \gamma_r \mu_A C_0^{dr} \\ &\qquad\qquad \text{for some } (x, v) \text{ with integer } v \in [-\Delta_r, t + \Delta_r) \\ &\qquad\qquad\quad \text{for which } \mathcal{Q}_r(x) \text{ intersects } \mathcal{C}(t \log t + 3\Delta_r)\} \\ &\quad \leq \sum_{r \geq R} P\{U_r(x, v) < \tfrac{1}{2} \mu_A C_0^{dr} \\ &\qquad\qquad \text{for some } (x, v) \text{ with integer } v \in [-\Delta_r, t + \Delta_r) \\ &\qquad\qquad\quad \text{for which } \mathcal{Q}_r(x) \text{ intersects } \mathcal{C}(t \log t + 3\Delta_r)\} \\ &\quad \leq t^{-K} \text{ for all large } t. \end{aligned} \tag{4.17}$$

(We used $\gamma_r \leq 1/2$ for the one but last inequality.) Note also that the required value of $K_4$ depends on $K, d$ and $\mu_A$ only.

We shall also need the following analogue of Lemma 6 in [7] (note that this time the inequality goes in the opposite direction from (5.22) in [7]). $\{S_u\}_{u \geq 0}$ is short for what we formerly denoted by $\{S_u^A\}_{u \geq 0}$, that is, a continuous-time simple random walk with jump rate $D$ (starting at $\mathbf{0}$).

LEMMA 5. *There exists a constant $C_4 = C_4(d, D)$, which is independent of $C_0$, such that for $r \geq 1$, if $\mathcal{V}_{r+1}(\mathbf{i}, k)$ is good, and $u$ an integer with $\Delta_{r+1} - \Delta_r \leq u \leq 2\Delta_{r+1}$, then for*

$$y \in \prod_{s=1}^{d} [(i(s) - 1)\Delta_{r+1}, (i(s) + 2)\Delta_{r+1}), \tag{4.18}$$

*it holds that*

$$\begin{aligned} &\sum_{z \in V_{r+1}(\mathbf{i})} N^*(z, (k-1)\Delta_{r+1}) P\{z + S_u \in \mathcal{Q}_r(y)\} \\ &\qquad \geq \gamma_{r+1} \mu_A C_0^{dr} \left[1 - \frac{C_4 (r \log C_0)^d}{C_0^r}\right]. \end{aligned} \tag{4.19}$$

PROOF. Let $r \geq 1$ be fixed. In addition to the blocks $\prod_{s=1}^d [i(s)\Delta_{r+1}, (i(s)+1)\Delta_{r+1})$ which have edge length $\Delta_{r+1} = C_0^{6(r+1)}$, we also need the blocks

$$\mathcal{M}(\mathbf{j}) := \prod_{s=1}^d [j(s)C_0^{r+1}, (j(s)+1)C_0^{r+1}).$$

In our previous notation $\mathcal{M}(\mathbf{j}) = \mathcal{Q}_{r+1}(x)$ with $x(s) = j(s)C_0^{r+1}$. These blocks have edge length only $C_0^{r+1}$, and the set $V_{r+1}(\mathbf{i})$ is a disjoint union of $7^d C_0^{5d(r+1)}$ of these smaller blocks. Let $\Lambda = \Lambda(\mathbf{i}, r+1)$ be the set of $\mathbf{j} \in \mathbb{Z}^d$ with

$$\mathcal{M}(\mathbf{j}) \subset V_{r+1}(\mathbf{i}).$$

Also, for each $\mathbf{j} \in \Lambda$ let $z_{\mathbf{j}} \in \mathcal{M}(\mathbf{j})$ be such that

$$P\{z_{\mathbf{j}} + S_u \in \mathcal{Q}_r(y)\} = \min_{z \in \mathcal{M}(\mathbf{j})} P\{z + S_u \in \mathcal{Q}_r(y)\}.$$

Then the left-hand side of (4.19) equals

$$\begin{aligned} &\sum_{\mathbf{j}\in\Lambda} \sum_{z\in\mathcal{M}(\mathbf{j})} N^*(z,(k-1)\Delta_{r+1}) P\{z + S_u \in \mathcal{Q}_r(y)\} \\ &\qquad \geq \sum_{\mathbf{j}\in\Lambda} \sum_{z\in\mathcal{M}(\mathbf{j})} N^*(z,(k-1)\Delta_{r+1}) P\{z_{\mathbf{j}} + S_u \in \mathcal{Q}_r(y)\}. \end{aligned} \tag{4.20}$$

Since $\mathcal{V}_{r+1}(\mathbf{i}, k)$ is assumed to be good, we have

$$\begin{aligned} \sum_{z\in\mathcal{M}(\mathbf{j})} N^*(z,(k-1)\Delta_{r+1}) &= U_{r+1}(\mathbf{j}C_0^{r+1}, (k-1)\Delta_{r+1}) \\ &\geq \gamma_{r+1}\mu_A C_0^{d(r+1)} = \sum_{z\in\mathcal{M}(\mathbf{j})} \gamma_{r+1}\mu_A. \end{aligned}$$

We can therefore continue (4.20) to obtain that the left-hand side of (4.19) is at least

$$\begin{aligned} &\sum_{\mathbf{j}\in\Lambda} \sum_{z\in\mathcal{M}(\mathbf{j})} \gamma_{r+1}\mu_A P\{z_{\mathbf{j}} + S_u \in \mathcal{Q}_r(y)\} \\ &\qquad \geq \sum_{\mathbf{j}\in\Lambda} \sum_{z\in\mathcal{M}(\mathbf{j})} \gamma_{r+1}\mu_A P\{z + S_u \in \mathcal{Q}_r(y)\} \\ &\qquad\quad - \sum_{\mathbf{j}\in\Lambda} \sum_{z\in\mathcal{M}(\mathbf{j})} \gamma_{r+1}\mu_A |P\{z_{\mathbf{j}} + S_u \in \mathcal{Q}_r(y)\} - P\{z + S_u \in \mathcal{Q}_r(y)\}|. \end{aligned} \tag{4.21}$$

Now, by virtue of (4.18), the first multiple sum in the right-hand side of (4.21) is at least

$$\sum_{z\,:\,z-y\in[-2\Delta_{r+1},2\Delta_{r+1})^d} \gamma_{r+1}\mu_A \sum_{w\in\mathcal{Q}_r(y-z)} P\{S_u = w\}$$

$$\begin{aligned}
&\geq \sum_{w\in[-\Delta_{r+1},\Delta_{r+1})^d} P\{S_u=w\} \sum_{z\in\mathcal{Q}_r(y-w)} \gamma_{r+1}\mu_A \\
&= \sum_{w\in[-\Delta_{r+1},\Delta_{r+1})^d} P\{S_u=w\}\gamma_{r+1}\mu_A C_0^{dr} \qquad (4.22)\\
&= \gamma_{r+1}\mu_A C_0^{dr}[1-P\{S_u\notin[-\Delta_{r+1},\Delta_{r+1})^d\}]\\
&\geq \gamma_{r+1}\mu_A C_0^{dr}[1-K_5\exp[-K_6\Delta_{r+1}]]
\end{aligned}$$

for some constants $K_5, K_6$, depending on $d, D_A$ only. In the last inequality we used simple large deviation estimates for $S_u$ (see, e.g., (2.40) in [7]) and the fact that $u \leq 2\Delta_{r+1}$.

On the other hand, we have for any $z \in \mathcal{M}(\mathbf{j})$ that

$$\begin{aligned}
&|P\{z_{\mathbf{j}}+S_u\in\mathcal{Q}_r(y)\}-P\{z+S_u\in\mathcal{Q}_r(y)\}|\\
&\quad\leq \sum_{w\in\mathcal{Q}_r(y)} |P\{z_{\mathbf{j}}+S_u=w\}-P\{z+S_u=w\}|\\
&\quad\leq \sum_{v\in\mathcal{Q}_r(y-z)} \sup_{w:\,\|w-v\|\leq C_0^{r+1}} |P\{S_u=v\}-P\{S_u=w\}|.
\end{aligned}$$

It follows that the second multiple sum in the right-hand side of (4.21) is bounded in absolute value by

$$\begin{aligned}
&\sum_z \gamma_{r+1}\mu_A \sum_{v\in\mathcal{Q}_r(y-z)} \sup_{w:\,\|w-v\|\leq C_0^{r+1}} |P\{S_u=v\}-P\{S_u=w\}|\\
&\quad\leq \gamma_{r+1}\mu_A \sum_{v\in\mathbb{Z}^d}\sum_{z\in\mathcal{Q}_r(y-v)} \sup_{w:\,\|w-v\|\leq C_0^{r+1}} |P\{S_u=v\}-P\{S_u=w\}|\\
&\quad= \gamma_{r+1}\mu_A C_0^{dr} \sum_{v\in\mathbb{Z}^d} \sup_{w:\,\|w-v\|\leq C_0^{r+1}} |P\{S_u=v\}-P\{S_u=w\}|.
\end{aligned}$$

The right-hand side here has been estimated in (6.37) and in (5.26) and the following lines in [7]. The result is that the right-hand side here is bounded by

$$K_7\gamma_{r+1}\mu_A C_0^{dr-2r-2}[r\log C_0]^d$$

for some constant $K_7$ which does not depend on $C_0$. Combining this with the estimates (4.21) and (4.22) we obtain (4.19). □

We define some $\sigma$-fields analogously to [7] (but with some differences):

$\mathcal{H}_{r+1}(\mathbf{i},k) :=$ $\sigma$-field generated by the paths of all particles through time $(k-1)\Delta_{r+1}$ and the paths through time $(k+1)\Delta_{r+1}-1$ of the particles which are (4.23)

$$\text{outside } V_{r+1}(\mathbf{i}) \text{ at time } (k-1)\Delta_{r+1},$$

$$\mathcal{K}_{r+1} := \sigma\text{-field generated by } \{N^*(x,(k-1)\Delta_{r+1}) : x \in V_{r+1}(\mathbf{i})\}.$$

Note that

$$\mathcal{K}_{r+1} \subset \mathcal{H}_{r+1}(\mathbf{i},k),$$

because if one knows all paths through time $(k-1)\Delta_{r+1}$, then one also knows how many particles there are at each $x$ at time $(k-1)\Delta_{r+1}$. In other words, all $N^*(x,(k-1)\Delta_{r+1}), x \in \mathbb{Z}^d$, are $\mathcal{H}_{r+1}(\mathbf{i},k)$-measurable.

We also need certain events $\mathcal{A}(\mathbf{i},k)$ which are somewhat larger than the event $\{\mathcal{B}_{r+1}(\mathbf{i},k)$ contains some bad $\mathcal{B}_r(\mathbf{j},q)\}$. For given $(\mathbf{i},k)$ and any $(y,v)$ with $v \geq (k-1)\Delta_{r+1}$ we define

$$W_r(y,v) = \text{ number of particles in the system } \mathcal{P}^* \text{ in } \mathcal{Q}_r(y) \times \{v\}$$
$$\text{which were in } V_{r+1}(\mathbf{i}) \text{ at time } (k-1)\Delta_{r+1}.$$

A block $\mathcal{B}_r(\mathbf{j},q) \subset \mathcal{B}_{r+1}(\mathbf{i},k)$ will be called *inferior* if $W_r(y,v) < \gamma_r \mu_A C_0^{dr}$ for some $(y,v)$ for which $v$ is an integer and $\mathcal{Q}_r(y) \times \{v\}$ is contained in $\widetilde{\mathcal{B}}_r(\mathbf{j},q)$.

It is apparent from the definitions that

$$W_r(y,v) \leq U_r(y,v), \tag{4.24}$$

since we count only particles which passed through $\mathcal{V}_{r+1}(\mathbf{i},k)$ in $W_r(y,v)$, whereas $U_r(y,v)$ also counts particles which do not satisfy this requirement. It follows from this that a bad block is also inferior. Finally, we define the event

$$\mathcal{A}(\mathbf{i},k) = \mathcal{A}(\mathbf{i},k,r) = \{\mathcal{B}_{r+1}(\mathbf{i},k) \text{ contains some inferior } r\text{-block } \mathcal{B}_r(\mathbf{j},q)\}.$$

One now has the following analogue of Lemma 7 and part of the proof of Lemma 8 in [7].

LEMMA 6. *Let*

$$\rho_{r+1} = 3^{d+1} C_0^{6(d+1)(r+1)} \exp[-\tfrac{1}{2}\gamma_r \mu_A C_0^{(d-3/4)r}], \qquad r \geq 1. \tag{4.25}$$

*Then for* $r \geq 1$, *on the event* $\{\mathcal{V}_{r+1}(\mathbf{i},k)$ *is good*$\}$,

$$P\{\mathcal{A}(\mathbf{i},k) | \mathcal{H}_{r+1}(\mathbf{i},k)\} = P\{\mathcal{A}(\mathbf{i},k) | \mathcal{K}_{r+1}(\mathbf{i},k)\} \leq \rho_{r+1}. \tag{4.26}$$

*Moreover, for fixed* $a(s) \in \{0,1,\dots,11\}$ *and* $b = 0$ *or* $1$, *the collection of pairs* $(\mathbf{i},k), i(s) \equiv a(s) \bmod 12, 1 \leq s \leq d, k \equiv b \bmod 2$, *for which* $\mathcal{V}_{r+1}(\mathbf{i},k)$ *is good, but* $\mathcal{A}(\mathbf{i},k)$ *occurs, is stochastically smaller than an independent percolation system in which each site* $(\mathbf{i},k), i(s) \equiv a(s) \bmod 12, 1 \leq s \leq d, k \equiv b \bmod 2$, *is open with probability* $\rho_{r+1}$.

PROOF. By the Markov property of the particle system $\mathcal{P}^*$, the conditional distribution of the particles during $[(k-1)\Delta_{r+1}, \infty)$, given the behavior of all particles during $(-\infty, (k-1)\Delta_{r+1}]$, is the same as the conditional distribution given the positions of all particles at time $(k-1)\Delta_{r+1}$. Moreover, given the positions of the particles at time $(k-1)\Delta_{r+1}$, the future paths of all particles are conditionally independent. In particular, given the particles at time $(k-1)\Delta_{r+1}$, the paths after time $(k-1)\Delta_{r+1}$ of the particles in $\mathcal{V}_{r+1}(\mathbf{i}, k)$ are conditionally independent of the future paths of all particles outside $\mathcal{V}_{r+1}(\mathbf{i}, k)$ at time $(k-1)\Delta_{r+1}$. By definition the event $\mathcal{A}(\mathbf{i}, k)$ depends only on the particles in $\mathcal{V}_{r+1}(\mathbf{i}, k)$ at time $(k-1)\Delta_{r+1}$ and the increments of their paths after time $(k-1)\Delta_{r+1}$. It follows from these comments that $P\{\mathcal{A}(\mathbf{i}, k)|\mathcal{H}_{r+1}(\mathbf{i}, k)\}$ is a function of the particles in $\mathcal{V}_{r+1}(\mathbf{i}, k)$ only. In fact, since $\mathcal{A}(\mathbf{i}, k)$ depends on particle counts only, $P\{\mathcal{A}(\mathbf{i}, k)|\mathcal{H}_{r+1}(\mathbf{i}, k)\}$ depends only on the $N^*(x, v)$ with $(x, v) \in \mathcal{V}_{r+1}(\mathbf{i}, k)$ (see the explicit computation in the next paragraph). Thus the left-hand side of (4.26) is $\mathcal{K}_{r+1}(\mathbf{i}, k)$-measurable, and equals the middle member of (4.26).

We next prove the inequality in (4.26). If $\mathcal{A}(\mathbf{i}, k)$ occurs, then $W_r(y, v) < \gamma_r \mu_A C_0^{dr}$ for some integer $v$ and

$$
\begin{aligned}
(y, v) \in & \bigcup_{B_r(\mathbf{j},q) \subset B_{r+1}(\mathbf{i},k)} \widetilde{B}_r(\mathbf{j}, q) \\
& \subset \prod [i(s)\Delta_{r+1} - 3\Delta_r, \\
& \qquad (i(s)+1)\Delta_{r+1} + 3\Delta_r) \times [k\Delta_{r+1} - \Delta_r, (k+1)\Delta_{r+1}) \\
& \subset \prod [(i(s)-1)\Delta_{r+1}, \\
& \qquad (i(s)+2)\Delta_{r+1} - \Delta_r) \times [k\Delta_{r+1} - \Delta_r, (k+1)\Delta_{r+1}).
\end{aligned} \tag{4.27}
$$

Now consider any $(y, v)$ satisfying (4.27) and let the particles in $\mathcal{V}_{r+1}(\mathbf{i}, k)$ be given such that $\mathcal{V}_{r+1}(\mathbf{i}, k)$ is good. Conditionally on this, the distribution of $W_r(y, v)$ is the distribution of

$$
\sum_{z \in V_{r+1}(\mathbf{i})} \sum_{q=1}^{N^*(z,(k-1)\Delta_{r+1})} I[z + S_u^{z,q} \in \mathcal{Q}_r(y)],
$$

where the $\{S^{z,q}\}$ are independent copies of the random walk $\{S\}$ and $u = v - (k-1)\Delta_{r+1} \in [\Delta_{r+1} - \Delta_r, 2\Delta_{r+1}]$ (see the proof of Lemma 7, and in particular the lines following (5.37) in [7]). Therefore, $P\{W_r(y, v) < \gamma_r \mu_A C_0^{dr}\}$ is the probability of fewer than $\gamma_r \mu_A C_0^{dr}$ successes in

$$
\sum_{z \in V_{r+1}(\mathbf{i})} N^*(z, (k-1)\Delta_{r+1})
$$

trials, $N^*(z, (k-1)\Delta_{r+1})$ of which have success probability

$$
p(y - z, u) := P\{z + S_u \in \mathcal{Q}_r(y)\}.
$$

Very much as in (5.38), (5.39) of [7] we therefore have for $\theta \ge 0$,

$$\begin{aligned} &E\{\exp[-\theta W_r(y,v)]|\mathcal{K}_{r+1}\} \\ &\quad = \prod_{z \in V_{r+1}(\mathbf{i})} [1 - p(y-z,u) + p(y-z,u)e^{-\theta}]^{N^*(z,(k-1)\Delta_{r+1})} \\ &\quad \le \exp\left[ - \sum_{z \in V_{r+1}(\mathbf{i})} N^*(z,(k-1)\Delta_{r+1}) p(y-z,u)(1-e^{-\theta}) \right]. \end{aligned}$$

For $\theta = C_0^{-r/2}$ this gives, by virtue of (4.3), (4.4) and Lemma 5, that on the event $\{\mathcal{V}_{r+1}(\mathbf{i},k)$ is good$\}$ and for $(y,v)$ satisfying (4.27),

$$\begin{aligned} &P\{W_r(y,v) < \gamma_r \mu_A C_0^{dr} | \mathcal{K}_{r+1}\} \\ &\quad \le \exp\Bigg[\theta\gamma_r\mu_A C_0^{dr} \\ &\qquad - \sum_{z \in V_{r+1}(\mathbf{i})} N^*(z,(k-1)\Delta_{r+1}) p(y-z,u)(1-e^{-\theta})\Bigg] \\ &\quad \le \exp[-\tfrac{1}{2}\gamma_r\mu_A C_0^{(d-3/4)r}]. \end{aligned} \tag{4.28}$$

The inequality in (4.26) now follows from the fact that $P\{\mathcal{A}(\mathbf{i},k)\}$ is bounded by the sum over $(y,v)$ with integral $v$ and satisfying (4.27).

The last statement of the lemma concerning the stochastic ordering between the collection of pairs $(\mathbf{i},k)$ for which $\mathcal{V}_{r+1}(\mathbf{i},k)$ is good and $\mathcal{A}(\mathbf{i},k)$ occurs, and an independent percolation system, now follows in exactly the same way as in the proof of (5.43) in [7]. ☐

The next lemma is basically a copy of parts of Lemma 8 and Lemma 11 in [7]. Note that now $R = R(t)$ is defined in (4.16).

LEMMA 7. *Inequalities (4.14) and (4.15) hold. Moreover, there exist some constants $C_5 = C_5(d,\mu_A)$, $\kappa_0 = \kappa(d,\mu_A)$ and $t_0 = t_0(d,\mu_A)$ (independent of $r,\ell$), such that for $1 \le r \le R(t) - 1, \kappa \ge \kappa_0, t \ge t_0$ and any $\ell \ge 0$*

$$\begin{aligned} &P\left\{\Psi_{r+1}(\ell) \ge \frac{\kappa(t+\ell)}{\Delta_{r+1}}[\rho_{r+1}]^{1/(d+1)}\right\} \\ &\quad \le \exp\left[-(t+\ell)C_5\kappa\exp\left[-\frac{1}{2(d+1)}\gamma_r\mu_A C_0^{(d-3/4)r}\right]\right]. \end{aligned} \tag{4.29}$$

PROOF. We already observed that (4.14) and (4.15) hold, for the same reasons as in Lemma 8 of [7].

Inequality (4.29) follows by a percolation argument which is given in the proof of (6.28) and Lemma 8 of [7]; see also the proof of Theorem 9 in [9]. This time we take an integer $\nu$ such that

$$[\rho_{r+1}]^{-1/(d+1)} \le \nu \le 2[\rho_{r+1}]^{-1/(d+1)} \tag{4.30}$$

and define

$$\begin{aligned}\mathcal{D}(\mathbf{m}, q) = \prod_{s=1}^{d} &[\nu m(s)\Delta_{r+1}, \nu(m(s)+1)\Delta_{r+1}) \\ &\times [q\nu\Delta_{r+1}, (q+1)\nu\Delta_{r+1}).\end{aligned} \tag{4.31}$$

$\mathbf{m}$ here is short for $(m(1), \dots, m(d))$; the $\mathbf{m}$ and $\nu$ here have nothing to do with the $m_i$ and $\nu$ in the proof of Proposition 4. Note that $\rho_{r+1} \le 1$ by (4.5) and (1.11), so that (4.30) can be satisfied. Each $\mathcal{D}(\mathbf{m}, q)$ is the disjoint union of $\nu^{d+1}(r+1)$-blocks. Moreover, as shown in (6.30) of [7], for $\ell \ge 0$ at most

$$\lambda(\ell) := 3^d \left( \frac{t+\ell}{\nu\Delta_{r+1}} + 2 \right) \tag{4.32}$$

blocks $\mathcal{D}(\mathbf{m}, q)$ can intersect a space–time path $\widehat{\pi} \in \Xi(\ell, t)$ (with jump times $s_1 < \cdots < s_\ell < t$ and positions $x_1, \dots, x_\ell$). Now fix $a(1), \dots, a(d) \in \{0, \dots, 11\}$, $b \in \{0, 1\}$ and define for any space–time path $\widehat{\pi}$

$$\begin{aligned}\psi_{r+1}(\widehat{\pi}, \mathbf{a}, k) = &\text{ number of } (r+1)\text{-blocks } \mathcal{B}_{r+1}(\mathbf{i}, k) \\ &\text{with } i(s) \equiv a(s) \bmod 12,\ k \equiv b \bmod 2, \\ &\text{which intersect the space–time path } \widehat{\pi} \text{ and} \\ &\text{which have a good pedestal but contain a bad } r\text{-block.}\end{aligned} \tag{4.33}$$

Define further

$$\Psi_{r+1}(\ell, \mathbf{a}, b) = \sup_{\widehat{\pi} \in \Xi(\ell, t)} \psi_{r+1}(\widehat{\pi}, \mathbf{a}, b).$$

Then

$$\Psi_{r+1}(\ell) \le \sum_{(\mathbf{a}, b)} \Psi_{r+1}(\ell, \mathbf{a}, b). \tag{4.34}$$

As in [7], let $Z(\mathbf{i}, k)$ be independent random variables with

$$P\{Z(\mathbf{i}, k) = 1\} = 1 - P\{Z(\mathbf{i}, k) = 0\} = \rho_{r+1}.$$

Then, as in (6.31) of [7],

$$P\left\{ \Psi_{r+1}(\ell, \mathbf{a}, b) \ge 2^{-1}(12)^{-d} \frac{\kappa(t+\ell)}{\Delta_{r+1}} [\rho_{r+1}]^{1/(d+1)} \right\}$$

$$\leq \sum_{\mathcal{D}(\mathbf{m}_0,0),\dots,\mathcal{D}(\mathbf{m}_{\lambda-1},\lambda-1)} P\Bigg\{\bigcup_{q=0}^{\lambda-1} \mathcal{D}(\mathbf{m}_q,q) \text{ contains at least}$$

$$2^{-1}(12)^{-d}\frac{\kappa(t+\ell)}{\Delta_{r+1}}[\rho_{r+1}]^{1/(d+1)} \tag{4.35}$$

$$(r+1)\text{-blocks } \mathcal{B}_{r+1}(\mathbf{i},k) \text{ with } Z(\mathbf{i},k)=1,$$

$$\text{and } i(s)\equiv a(s) \bmod 12, k\equiv b \bmod 2\Bigg\}.$$

Here $(\mathcal{D}(\mathbf{m}_0,0),\dots,\mathcal{D}(\mathbf{m}_{\lambda-1},\lambda-1))$ runs over the possible collections of blocks $\mathcal{D}$ which intersect a space–time path $\widehat{\pi}\in\Xi(\ell,t)$. For some constant $K_8$ which depends on $d$ only, there are at most

$$[2t\log t+1]^d\exp[K_8\lambda] \tag{4.36}$$

collections of this form. If we fix such a collection $\mathcal{D}(\mathbf{m}_0,0),\dots,\mathcal{D}(\mathbf{m}_{\lambda-1},\lambda-1)$, then the probability that

$$\bigcup_{q=0}^{\lambda-1}\mathcal{D}(\mathbf{m}_q,q)$$

contains at least $2^{-1}(12)^{-d}\kappa(t+\ell)\Delta_{r+1}^{-1}[\rho_{r+1}]^{1/(d+1)}$ $(r+1)$-blocks $\mathcal{B}_{(r+1)}(\mathbf{i},k)$ with $Z(\mathbf{i},k)=1$ and $(\mathbf{i},k)\equiv(\mathbf{a},b)$, is bounded by

$$P\bigg\{T\geq 2^{-1}(12)^{-d}\frac{\kappa(t+\ell)}{\Delta_{r+1}}[\rho_{r+1}]^{1/(d+1)}\bigg\}, \tag{4.37}$$

where $T$ has a binomial distribution corresponding to $\lambda\nu^{d+1}$ trials with success probability $\rho_{r+1}$. As in (6.33) or (5.52) in [7] one obtains from Bernstein's inequality [together with (4.30) and (4.32)] that this probability is at most

$$K_9\exp\bigg[-K_{10}\frac{\kappa(t+\ell)}{\Delta_{r+1}}[\rho_{r+1}]^{1/(d+1)}\bigg]$$

for $1\leq r\leq R(t)-1,\kappa\geq$ some $\kappa_0$, $t\geq$ some $t_0$, and constants $K_9,K_{10}$, depending on $d$ and $\mu_A$ only. Inequality (4.29) now follows from (4.35), (4.36) and (4.34). □

PROPOSITION 8. *For any choice of $K$ and $\varepsilon_0>0$, there exist constants $r_0,t_1$ such that for all $t\geq t_1$,*

$$P\{\Phi_r(\ell)\geq\varepsilon_0 C_0^{-6r}(t+\ell) \text{ for some } r\geq r_0,\ell\geq 0\}\leq\frac{2}{t^K}. \tag{4.38}$$

Proof. Consider a sample point for which

$$\Phi_r(\ell)=0 \qquad \text{for all } r\geq R(t) \text{ and } \ell\geq 0, \tag{4.39}$$

and for which

$$\Phi_r(\ell)\leq C_0^{6(d+1)}\Phi_{r+1}(\ell)+C_0^{6(d+1)}\frac{\kappa_0(t+\ell)}{\Delta_{r+1}}[\rho_{r+1}]^{1/(d+1)} \tag{4.40}$$

for all $t\geq t_0, 1\leq r\leq R-1, \ell\geq 0$. For such a sample point one also has for $t\geq t_0, r_0\leq r\leq R-1, \ell\geq 0$,

$$\begin{aligned}\Phi_r(\ell)&\leq C_0^{6(d+1)}\frac{\kappa_0(t+\ell)}{\Delta_{r+1}}[\rho_{r+1}]^{1/(d+1)}+C_0^{6(d+1)}\Phi_{r+1}(\ell)\\&\leq C_0^{6(d+1)}3\kappa_0(t+\ell)\exp\left[-\frac{\gamma_0\mu_A}{2(d+1)}C_0^{(d-3/4)r}\right]+C_0^{6(d+1)}\Phi_{r+1}(\ell)\\&\leq C_0^{6(d+1)}3\kappa_0(t+\ell)\exp\left[-\frac{\gamma_0\mu_A}{2(d+1)}C_0^{(d-3/4)r}\right]\\&\quad+C_0^{12(d+1)}3\kappa_0(t+\ell)\exp\left[-\frac{\gamma_0\mu_A}{2(d+1)}C_0^{(d-3/4)(r+1)}\right]\\&\quad+C_0^{12(d+1)}\Phi_{r+2}(\ell)\\&\leq\cdots\leq\sum_{j=1}^{R-r}C_0^{6j(d+1)}3\kappa_0(t+\ell)\exp\left[-\frac{\gamma_0\mu_A}{2(d+1)}C_0^{(d-3/4)(r+j-1)}\right]\\&\quad+C_0^{6(d+1)(R-r)}\Phi_R(\ell)\\&\leq 6\kappa_0(t+\ell)C_0^{6(d+1)}\exp\left[-\frac{\gamma_0\mu_A}{2(d+1)}C_0^{(d-3/4)r}\right]\\&\leq\varepsilon_0C_0^{-6r}(t+\ell),\end{aligned}\tag{4.41}$$

provided $r_0$ is sufficiently large. The required value for $r_0$ is independent of $t,\ell$.

By (4.15)–(4.17), and (4.29), relations (4.39) and (4.40) hold outside a set of probability

$$t^{-K}+\sum_{r=1}^{R-1}\sum_{\ell\geq 0}\exp\left[-(t+\ell)C_5\kappa_0\exp\left[-\frac{1}{2(d+1)}\gamma_r\mu_AC_0^{(d-3/4)r}\right]\right]\leq 2t^{-K},$$

provided $t\geq$ some $t_1\geq t_0$. This proves the proposition. □

Proposition 8 is the main technical estimate for the proof of Theorem 2. We now start on this proof proper. The strategy will be to show that (with overwhelming probability) there exists a (random) path $u\mapsto\lambda(u,x)$ along which a $B$-particle moves with a "drift" toward a fixed point $x$, at least at the times $u$ when there are at least two particles at $\lambda(u, x)$. Lemma 9 below expresses this more formally in terms of quantities $I_1, I_{\geq 2}, \Gamma_1$ and $\Gamma_{\geq 2}$ which will be defined in (4.42) and (4.43) below; at the times when $I_{\geq 2}(u) = 1$ [and hence $I_1(u) = 0$], $\|\lambda(u, x) - x\|_2$ has a drift $D\Gamma_{\geq 2}(u)$, which will be shown to be negative in (4.68).

We now give the details. For $x \in \mathbb{Z}^d$ construct a path $\lambda(\cdot) = \lambda(\cdot, x) \in \mathbb{Z}^d$ by the rules (i)–(v) below:

(i) $\lambda(0, x)$ is the location of some initial $B$-particle, say $\lambda(0, x) = z_0$;

(ii) for all times $s$ there is a distinguished $B$-particle, $\widehat{\rho}(s)$ say, at $\lambda(s, x)$; at time 0 we designate any of the $B$-particles at $z_0$ as $\widehat{\rho}(0)$;

(iii) $s \mapsto \lambda(s, x)$ can jump only at times when $\widehat{\rho}(s-)$ jumps away from $\lambda(s, x)$, and $\lambda(\cdot, x)$ is constant between such jumps;

(iv) if $\widehat{\rho}(s-)$ jumps from $\lambda(s-, x) = w$ to $w'$ at some time $s$, and if this was the only particle at $w$ at time $s-$, then $\lambda(\cdot, x)$ also jumps to $w'$ at time $s$ [so that $\lambda(s, x) = w'$] and $\widehat{\rho}(s) = \widehat{\rho}(s-)$, the particle which jumped at time $s$;

(v) if $\widehat{\rho}(s-)$ jumps from $\lambda(s-, x) = w$ to $w'$ at some time $s$ such that there is at least one other particle $\rho'$ at $w$ at time $s-$, then $\lambda(\cdot, x)$ jumps to $w'$ at time $s$ if and only if $\|w' - x\|_2 < \|w - x\|_2$, and in this case again $\widehat{\rho}(s) = \widehat{\rho}(s-)$; if, however, $\|w' - x\|_2 \geq \|w - x\|_2$, then $\lambda(\cdot, x)$ does not jump at time $s$ and we take $\widehat{\rho}(s) = \rho'$.

In general, these rules do not determine $\lambda(\cdot)$ uniquely, because there may be more than one possible choice for $\rho'$ in rule (v). However, we can use any a priori rule to make $\lambda$ unique. We can, for instance, choose $\rho'$ as the particle $\rho_k$ with the minimal $k$ among all possible ones.

We shall say that the distinguished particle *attempts a jump* at time $s$ when the distinguished particle at time $s-$, that is, $\widehat{\rho}(s-)$, jumps away from $\lambda(s)$. Due to the fact that we may then declare another particle to be the distinguished one at time $s$ [see rule (v)], an attempted jump of the distinguished particle at time $s$ does not necessarily make $\lambda(s) \neq \lambda(s-)$.

We start with the distinguished particle being of type $B$, and by rules (iii)–(v), right after each attempted jump of the distinguished particle, there still is a distinguished $B$-particle at the location of $\lambda$. From this it is easy to check recursively, from one attempted jump of the distinguished particle to the next, that $\lambda(\cdot, x)$ automatically satisfies (ii). We merely have to note that if the distinguished particle has type $B$ just before an attempted jump, then all particles which coincide with the distinguished particle at that time also have type $B$.

Note that $\lambda(\cdot, x)$ can change only at the time at which the distinguished particle $\widehat{\rho}$ jumps [by rule (iii)]. The next lemma shows that the jumps which occur at a time $t$ when there are more than one particle at $\lambda(t)$ cause a certain drift of $\lambda(\cdot)$ toward $x$. We remind the reader that $\mathcal{F}_t^0 = \sigma$-field generated by $\{Y_s : s \le t\}$ and that

$$\mathcal{F}_t = \bigcap_{h>0} \mathcal{F}_{t+h}^0.$$

We remind the reader that $P$ without superscript is discussed just before Proposition 4. We take $e_{d+i} = -e_i$ for $1 \le i \le d$, and define

$$\begin{aligned} I_1(u) &= I[N_B(\lambda(u,x),u) = 1] \\ &= I[\widehat{\rho}(u) \text{ is the only particle present at } (\lambda(u,x),u)], \\ I_{\ge 2}(u) &= I[N_B(\lambda(u,x),u) \ge 2], \end{aligned} \tag{4.42}$$

$$\begin{aligned} \Gamma_1(u) &= \frac{1}{2d} \sum_{i=1}^{2d} [\|\lambda(u,x) + e_i - x\|_2 - \|\lambda(u,x) - x\|_2], \\ \Gamma_{\ge 2}(u) &= \frac{1}{2d} \sum{}^* [\|\lambda(u,x) + e_i - x\|_2 - \|\lambda(u,x) - x\|_2], \end{aligned} \tag{4.43}$$

where $\sum^*$ is the sum over those $i \in \{1, \dots, 2d\}$ for which

$$\|\lambda(u,x) + e_i - x\|_2 - \|\lambda(u,x) - x\|_2 < 0.$$

LEMMA 9.

$$\begin{aligned} M(t) = M(t,x) &:= \|\lambda(t,x) - x\|_2 \\ &\qquad - D \int_0^t [I_1(u)\Gamma_1(u) + I_{\ge 2}(u)\Gamma_{\ge 2}(u)]\, du \end{aligned} \tag{4.44}$$

*is a right-continuous $\{\mathcal{F}_t\}$-martingale under the measure $P$.*

The proof of this lemma is standard and we shall skip it here. The reader can find a proof in [8], Lemma 10.

We now want to use known exponential bounds for large deviations of martingales with suitable bounds on their increments. The following lemma is a special case of estimates for discrete-time (super) martingales with bounded jumps, such as can be found in [10], pages 154–155 (see also the estimation of $\lambda$ on page 334 in [6]).

LEMMA 10. *Assume that $\{\mathcal{G}_n\}_{n\ge 0}$ is an increasing sequence of $\sigma$-fields and that $D_n, n \ge 1$, are random variables which satisfy for all $n \ge 1$*

$$D_n \text{ is } \mathcal{G}_n\text{-measurable}, \tag{4.45}$$

$$|D_n| \le c, \tag{4.46}$$

*for some constant* $0 \le c < \infty$, *and*

$$E\{D_n | \mathcal{G}_{n-1}\} = 0. \tag{4.47}$$

*Define* $V_0 = 0$ *and*

$$V_n = \sum_{i=1}^{n} D_i, \qquad A_n = \sum_{i=1}^{n} E\{D_n^2 | \mathcal{G}_{n-1}\},$$

*for* $n \ge 1$. *Then* $\{V_n\}_{n\ge 0}$ *is a* $\{\mathcal{G}_n\}$*-martingale and there exists a constant* $K_3$, *depending on* $c$ *only, such that*

$$\begin{aligned} P\{|V_n| \ge a + bA_n \text{ for some } n \ge 0\} \le 2\exp[-K_3 ab], \\ a \ge 0, 0 \le b \le 1. \end{aligned} \tag{4.48}$$

To deduce estimates for $M(t)$ from this lemma, we define $\sigma_0 = 0$ and for $k \ge 0$

$$\begin{aligned} \sigma_{k+1} = \min[\sigma_k + 1, \inf\{t > \sigma_k : &\text{the distinguished particle} \\ &\widehat{\rho} \text{ attempts a jump at time } t\}]. \end{aligned}$$

We further take $D_n = M(\sigma_n) - M(\sigma_{n-1}), n \ge 1$, and $\mathcal{G}_n = \mathcal{F}_{\sigma_n}$. We then have $V_0 = 0$ and

$$\begin{aligned} V_n &= M(\sigma_n, x) - M(0, x) \\ &= M(\sigma_n, x) - \|\lambda(0, x) - x\|_2 \\ &= M(\sigma_n, x) - \|z_0 - x\|_2, \end{aligned}$$

with $M$ given by (4.44). It is immediate from the definitions that

$$\sup_{\sigma_n \le s \le \sigma_{n+1}} |M(s) - M(\sigma_n)| \le 1 + D. \tag{4.49}$$

Thus, (4.45)–(4.47) are satisfied with $c = 1 + D$. Moreover, $A_n \le c^2 n$. Consequently,

$$\begin{aligned} &P\{|M(\sigma_n) - M(0)| \ge a + bn \text{ for some } n \ge 0\} \\ &\quad \le P\left\{|M(\sigma_n) - M(0)| \ge a + \frac{b}{c^2} A_n \text{ for some } n \ge 0\right\} \\ &\quad \le 2\exp[-K_3 abc^{-2}] \qquad \text{if } b \le c^2. \end{aligned} \tag{4.50}$$

The attempted jump times of $\widehat{\rho}$ are distributed like the jump times of a rate-$D$ Poisson process, so that

$$P\{\sigma_{\lfloor 2Dt \rfloor} \le t\} \le \sum_{k \ge \lfloor 2Dt \rfloor} e^{-Dt} \frac{[Dt]^k}{k!} \le K_4 \exp[-K_5 t]. \tag{4.51}$$

Combined with (4.49) and (4.50) this shows that for $a \geq 2 + 2D$, $0 \leq b \leq 1$ and some constant $K_6 = K_6(D) > 0$

$$
\begin{aligned}
&P\Big\{\sup_{s \leq t} |M(s) - M(0)| \geq a + bt\Big\} \\
&\quad \leq P\{\sigma_{\lfloor 2Dt \rfloor} \leq t\} \\
&\qquad + P\Big\{|M(\sigma_n) - M(0)| \geq a - 1 - D + \frac{b}{(1+2D)} n \\
&\qquad\qquad \text{for some } n \leq 2Dt\Big\} \\
&\quad \leq K_4 \exp[-K_5 t] + 2\exp[-K_6 ab].
\end{aligned}
\tag{4.52}
$$

In particular, if we take $K_6 \leq K_5$ (as we may), then we obtain for $a = bt, 0 \leq b \leq 1$ and $t \geq t_2 := 2(1+D)/b$, that

$$
P\Big\{\sup_{s \leq t} |M(s) - M(0)| \geq 2bt\Big\} \leq (2 + K_4) \exp[-K_6 b^2 t]. \tag{4.53}
$$

Next we must find a lower bound for $\int_0^t I_{\geq 2}(u)\, du$. Before we can do this we need a preparatory lemma. For $L \geq 2$ we define

$$
\beta(L, d) = \begin{cases} 1, & \text{if } d = 1, \\ [\log L]^{-1}, & \text{if } d = 2, \\ L^{2-d}, & \text{if } d \geq 3, \end{cases}
$$

$$
\begin{aligned}
\mathcal{E}_n = \{&\text{there is some particle } \rho' \neq \widehat{\rho}(3L^2(n-1)) \text{ in} \\
&\lambda(3L^2(n-1), x) + [-L, L]^d \text{ at time } 3L^2(n-1)\}
\end{aligned}
\tag{4.54}
$$

and

$$
\begin{aligned}
J_n = I[&\widehat{\rho}(u) \text{ coincides with another particle} \\
&\text{at some time } u \in (3L^2(n-1), L^2(3n-1)]].
\end{aligned}
$$

Lemma 11. *There exists a constant $K_7 > 0$, depending on $d$ only, such that for all $4 \leq L^2 \leq t/4$*

$$
E\{J_n | \mathcal{F}_{3L^2(n-1)}\} \geq K_7 \beta(L, d) \qquad \text{on the event } \mathcal{E}_n. \tag{4.55}
$$

Proof. Fix an integer $n$ and for brevity write $m$ for $3L^2(n-1)$. Recall that the position of $\widehat{\rho}(m)$ is $\lambda(m)$. Write $\rho''$ for the particle which is the distinguished particle $\widehat{\rho}(m)$ at time $m$. Of course $\rho''$ does not have to be the distinguished particle anymore at some later time $u$. However, $\rho''$ can fail to be the distinguished particle at time $u > m$ only if rule (v) is invoked at some time in $(m, u]$.

Now $\rho'$ and $\rho''$ continue after time $m$ to perform random walks $\{S'\}$ and $\{S''\}$ which are independent of each other and all other particles. Since, on the event $\mathcal{E}_n$, these two particles have a distance at most $L\sqrt{d}$ from each other at time $m$, we can use standard random walk estimates to find a lower bound for

$$\begin{aligned} &P\{\rho' \text{ and } \rho'' \text{ coincide at some time } u \in (m, m+L^2] | \mathcal{F}_m\} \\ &\quad = P\{S'_u - S''_u = z \text{ for some } u \le L^2\}, \end{aligned}$$

where $z = -\pi(m, \rho') + \pi(m, \rho'') = -\pi(m, \rho') + \pi(m, \widehat{\rho}(m))$. Indeed,

$$\begin{aligned} &\int_{u \le L^2} P\{S'_u - S''_u = z\}\, du \\ &\quad = E\{\text{amount of time during } [0, L^2] \text{ with } S'_u - S''_u = z\} \\ &\quad = \int_{s \le L^2} P\{\text{smallest } u \text{ with } S'_u - S''_u = z \text{ lies in } ds\} \\ &\qquad\qquad \times E\{\text{amount of time during } [0, L^2 - s] \text{ with } S'_u - S''_u = \mathbf{0}\} \\ &\quad \le P\{S'_u - S''_u = z \text{ for some } u \le L^2\} \int_{u \le L^2} P\{S'_u - S''_u = \mathbf{0}\}\, du. \end{aligned}$$

The integrals in the extreme left- and right-hand sides here can be estimated by means of the local central limit theorem to obtain that

$$P\{\rho' \text{ and } \rho'' \text{ coincide at some time } u \in (m, m+L^2] | \mathcal{F}_m\} \ge K_7 \beta(L, d) \tag{4.56}$$

(see Theorem 2.2 in [1] or Lemmas 5.1 and 5.2 in [4] for similar arguments). If $\rho''$ is still the distinguished particle at the time $u$ when $\rho'$ and $\rho''$ coincide, then there are at least the two particles $\rho'$ and $\rho''$ at $\lambda(u)$ at time $u$ so that $I_{\ge 2}(u) = 1$ and $J_n = 1$ in this case.

As pointed out, $\rho''$ does not have to equal the distinguished particle $\widehat{\rho}(u)$ at a time $u$, but this can happen only if for some $u' \in (m, u]$ we use rule (v). This means that there must have been some time $u' \in (m, u]$ at which $\widehat{\rho}(u')$ coincided with another particle. It follows that $J_n = 1$ also in this case. □

We next derive a lower bound for

$$Z(t) = Z(t, x) := \int_0^t I_{\ge 2}(u)\, du$$

in terms of

$$V(t, L) = V(t, L, x)$$

$$\begin{aligned} &:= \sum_{1\le n\le L^{-2}t/3} I[\text{there is some particle other than } \widehat{\rho}(3L^2(n-1)) \\ &\qquad \text{inside } \lambda(3L^2(n-1),x)+[-L,L]^d \\ &\qquad\qquad \text{at time } 3L^2(n-1)]. \end{aligned} \tag{4.57}$$

LEMMA 12. *For $0<\varepsilon\le 1$ and $u\le L^2\le t/4$*

$$\begin{aligned} &P\{Z(t)\le \varepsilon\beta(L,d)L^{-2}t\} \\ &\quad \le P\left\{V(t,L)\le \frac{2\varepsilon}{K_7}e^{2D}L^{-2}t\right\} \\ &\qquad + 2\exp\left[-\frac{K_3}{3}\varepsilon^2\beta^2(L,d)L^{-2}t\right]. \end{aligned} \tag{4.58}$$

PROOF. We define

$$\begin{aligned} \mathcal{G}_n &= \mathcal{F}_{3L^2n}, \\ \widetilde{J}_n &= \min\left\{1, \int_{3L^2(n-1)}^{3L^2n} I_{\ge 2}(u)\,du\right\}, \\ D_n &= \widetilde{J}_n - E\{\widetilde{J}_n|\mathcal{G}_{n-1}\}. \end{aligned}$$

Note that $0\le \widetilde{J}_n\le 1$, so that $|D_n|\le 1$, $E\{D_n^2|\mathcal{G}_{n-1}\}\le 1$ and

$$A_{\lfloor L^{-2}t/3\rfloor} := \sum_{1\le n\le L^{-2}t/3} E\{D_n^2|\mathcal{G}_{n-1}\}\le L^{-2}t/3.$$

Therefore, (4.48) with $c=1$,

$$a=\varepsilon\beta(L,d)L^{-2}t/3 \quad\text{and}\quad b=\varepsilon\beta(L,d)$$

yields

$$\begin{aligned} &P\left\{\left|\sum_{1\le n\le L^{-2}t/3}[\widetilde{J}_n - E\{\widetilde{J}_n|\mathcal{G}_{n-1}\}]\right|\ge \varepsilon\beta(L,d)L^{-2}t\right\} \\ &\quad \le 2\exp\left[-\frac{K_3}{3}\varepsilon^2\beta^2(L,d)L^{-2}t\right]. \end{aligned}$$

In particular,

$$\begin{aligned} &P\{Z(t)\le \varepsilon\beta(L,d)L^{-2}t\} \\ &\quad \le P\left\{\sum_{1\le n\le L^{-2}t/3}\widetilde{J}_n\le \varepsilon\beta(L,d)L^{-2}t\right\} \end{aligned}$$

$$\le P\Biggl\{\sum_{1\le n\le L^{-2}t/3} E\{\widetilde{J}_n|\mathcal{G}_{n-1}\} \le 2\varepsilon\beta(L,d)L^{-2}t\Biggr\}$$

$$+ P\Biggl\{\Biggl|\sum_{1\le n\le L^{-2}t/3} [\widetilde{J}_n - E\{\widetilde{J}_n|\mathcal{G}_{n-1}\}]\Biggr| \ge \varepsilon\beta(L,d)L^{-2}t\Biggr\} \tag{4.59}$$

$$\le P\Biggl\{\sum_{1\le n\le L^{-2}t/3} E\{\widetilde{J}_n|\mathcal{G}_{n-1}\} \le 2\varepsilon\beta(L,d)L^{-2}t\Biggr\}$$

$$+ 2\exp\biggl[-\frac{K_3}{3}\varepsilon^2\beta^2(L,d)L^{-2}t\biggr].$$

Finally, we observe that on the event $\mathcal{E}_n$ [see (4.54)]

$$\begin{aligned} E\{\widetilde{J}_n|\mathcal{G}_{n-1}\} \ge P\{&\widehat{\rho}(u) \text{ coincides with another particle } \rho' \text{ at some time}\\ &u \in (3L^2(n-1), L^2(3n-1)] \text{ and the positions of } \widehat{\rho}(u')\\ &\text{and } \rho' \text{ and } \lambda(u',x) \text{ all stay together for } u' \in [u,u+1]|\mathcal{G}_{n-1}\}\\ &\ge \exp[-2D]P\{J_n = 1|\mathcal{G}_{n-1}\}\\ &\ge \exp[-2D]K_7\beta(L,d) \qquad [\text{by } (4.55)]. \end{aligned}$$

The lemma now follows from

$$\begin{aligned} \sum_{1\le n\le L^{-2}t/3} E\{\widetilde{J}_n|\mathcal{G}_{n-1}\} &\ge \exp[-2D]K_7\beta(L,d) \sum_{1\le n\le L^{-2}t/3} I[\mathcal{E}_n]\\ &= \exp[-2D]K_7\beta(L,d)V(t,L). \end{aligned}$$

□

The next lemma gives an upper bound for $\int_0^t [I_1(u)\Gamma_1(u) + I_{\ge 2}(u)\Gamma_{\ge 2}(u)]\,du$ in terms of $Z(t)$.

LEMMA 13. *There exist constants $0 < K_8, K_9, K_{10} < \infty$, which depend on $d$ and $D$ only, such that for all $z > 0$*

$$\begin{aligned} &\int_0^t [I_1(u)\Gamma_1(u) + I_{\ge 2}(u)\Gamma_{\ge 2}(u)]\,du\\ &\quad \le \frac{K_8\beta(L,d)}{zL^2}t + K_8\int_0^t I\biggl[\|\lambda(u) - x\|_2 \le \frac{zL^2}{\beta(L,d)}\biggr]\,du - K_9 Z(t). \end{aligned} \tag{4.60}$$

*Consequently, for*

$$0 < \varepsilon \le 1, \qquad \|x - z_0\|_2 \le \frac{K_9\varepsilon D\beta(L,d)}{4L^2}t, \qquad z = \frac{4K_8}{K_9\varepsilon}, \tag{4.61}$$

$$L \ge L_0 := \biggl[\frac{\varepsilon D K_9}{8}\biggr]^{1/2} \vee 3, \tag{4.62}$$

*and for* $t \geq L^2[\varepsilon\beta(L,d)]^{-1}t_3$ *for some* $t_3 = t_3(D)$ *it holds that*

$$
\begin{aligned}
&P\biggl\{\int_0^t I\biggl[\|\lambda(u)-x\|_2 \leq \frac{zL^2}{\beta(L,d)}\biggr]\,du \leq \frac{K_9\varepsilon\beta(L,d)}{4K_8L^2}t\biggr\} \\
&\qquad \leq (2+K_4)\exp[-K_{10}\varepsilon^2\beta^2(L,d)L^{-4}t] \\
&\qquad\quad + P\biggl\{V(t,L) \leq \frac{2\varepsilon}{K_7}e^{2D}L^{-2}t\biggr\} + 2\exp\biggl[-\frac{K_3}{3}\varepsilon^2\beta^2(L,d)L^{-2}t\biggr].
\end{aligned}
\tag{4.63}
$$

Proof. We shall show by simple calculus that there exist some constants $0 < K_8, K_9 < \infty$ which depend on $d$ only, such that for $\lambda, x \in \mathbb{Z}^d$

$$
\frac{1}{2d}\sum_{i=1}^{2d}[\|\lambda + e_i - x\|_2 - \|\lambda - x\|_2] \leq \frac{K_8}{\|\lambda - x\|_2 + 1}, \tag{4.64}
$$

and, with $\sum^*$ as in (4.43),

$$
\frac{1}{2d}\sum{}^*[\|\lambda + e_i - x\|_2 - \|\lambda - x\|_2] \leq -K_9 + \frac{K_8}{\|\lambda - x\|_2 + 1}. \tag{4.65}
$$

Moreover, the left-hand sides of (4.64) and (4.65) are at most 1 in absolute value.

Before we prove these inequalities we show that they imply the lemma. Indeed, it follows from (4.64), (4.65) and the definitions (4.43) that the left-hand side of (4.60) is at most

$$
\begin{aligned}
&\frac{K_8\beta(L,d)}{zL^2}\int_0^t (I_1(u) + I_{\geq 2}(u))I\biggl[\|\lambda(u)-x\|_2 > \frac{zL^2}{\beta(L,d)}\biggr]\,du \\
&\quad + K_8\int_0^t (I_1(u) + I_{\geq 2}(u))I\biggl[\|\lambda(u)-x\|_2 \leq \frac{zL^2}{\beta(L,d)}\biggr]\,du - K_9\int_0^t I_{\geq 2}(u)\,du \\
&\quad \leq \frac{K_8\beta(L,d)}{zL^2}t + K_8\int_0^t I\biggl[\|\lambda(u)-x\|_2 \leq \frac{zL^2}{\beta(L,d)}\biggr]\,du - K_9Z(t).
\end{aligned}
$$

This proves (4.60).

To prove (4.63) we take $b = \varepsilon DK_9\beta(L,d)/(8L^2)$ in (4.53). For $L \geq L_0$ this $b$ satisfies $b \leq 1$. Then we obtain, by means of (4.60), that outside a set of probability at most $(2+K_4)\exp[-K_6b^2t]$ it holds that

$$
\begin{aligned}
0 &\leq \|\lambda(t) - x\|_2 \\
&\leq \|\lambda(0) - x\|_2 + 2bt + \frac{DK_8\beta(L,d)}{zL^2}t \\
&\quad + DK_8\int_0^t I\biggl[\|\lambda(u)-x\|_2 \leq \frac{zL^2}{\beta(L,d)}\biggr]\,du - DK_9Z(t)
\end{aligned}
$$

for $t \geq t_2 = L^2[\varepsilon\beta(L,d)]^{-1}t_3$ for some $t_3 = t_3(D)$. By substitution of the chosen values of $x, b$ and $z$ this yields

$$\int_0^t I\left[\|\lambda(u) - x\|_2 \leq \frac{zL^2}{\beta(L,d)}\right] du \geq \frac{K_9}{K_8} Z(t) - \frac{3K_9\varepsilon\beta(L,d)}{4K_8L^2} t. \tag{4.66}$$

If we exclude a further set of probability at most equal to the right-hand side of (4.58), then the right-hand side of (4.66) exceeds $K_9\varepsilon\beta(L,d)[4K_8L^2]^{-1}t$. Thus (4.63) also follows from (4.64) and (4.65).

We turn to the proof of (4.64) and (4.65). The sentence following (4.65) is trivial. We can therefore adjust $K_8$ so that (4.64) and (4.65) are valid on any given finite set of values for $\|\lambda - x\|_2$. In particular, we may restrict ourselves to proving (4.64), (4.65) for $\|\lambda - x\|_2 \geq 2$. Now the Taylor expansion

$$\begin{aligned} \|a+b\|_2 &= \sqrt{\|a\|_2^2 + 2a\cdot b + \|b\|_2^2} \\ &= \|a\|_2 + \frac{2a\cdot b + \|b\|_2^2}{2\|a\|_2} + O\left(\frac{\|a\|_2^2\|b\|_2^2 + \|b\|_2^4}{\|a\|_2^3}\right) \end{aligned}$$

shows that the left-hand side of (4.64) equals

$$\frac{1}{2d}\sum_{i=1}^{2d} \frac{(\lambda - x)\cdot e_i}{\|\lambda - x\|_2} + \frac{H(\lambda - x)}{\|\lambda - x\|_2} = \frac{H(\lambda - x)}{\|\lambda - x\|_2} \tag{4.67}$$

for some function $H$ which is bounded on $\{\lambda \neq x\} = \{\|\lambda - x\|_2 \geq 1\}$ (recall that $\lambda, x \in \mathbb{Z}^d$). Thus (4.64) holds.

For (4.65) we write $\lambda - x = \sum_{i=1}^d n_i e_i$, with integer coefficients $n_i$ (since $\lambda - x \in \mathbb{Z}^d$). Then for a given $i \in \{1, \dots, d\}$ there are three possibilities: $n_i > 0, n_i < 0, n_i = 0$. If $n_i > 0$, and hence $n_i \geq 1$, then $d + i$ is contained in $\sum^*$, but not $i$. Thus in this case, $2d$ times the contribution of the term with $d + i$ to the left-hand side of (4.65) is

$$\begin{aligned} &\left[\sum_{k\neq i} n_k^2 + (n_i - 1)^2\right]^{1/2} - \left[\sum_{k\neq i} n_k^2 + n_i^2\right]^{1/2} \\ &\quad = \left\{\left[\sum_{k\neq i} n_k^2 + (n_i - 1)^2\right]^{1/2} + \left[\sum_{k\neq i} n_k^2 + n_i^2\right]^{1/2}\right\}^{-1} (-2n_i + 1) \\ &\quad \leq -\tfrac{1}{2}\left[\sum_{k\neq i} n_k^2 + n_i^2\right]^{-1/2} n_i. \end{aligned}$$

[In the last inequality we used that $2n_i - 1 \geq n_i$ and that the term in the denominator with $n_i^2$ exceeds that with $(n_i - 1)^2$.] If $n_i < 0$, then only $i$ is contained in $\sum^*$, but not $d + i$, and the preceding estimates hold with $n_i$ replaced by $-n_i$. Finally, if $n_i = 0$, then neither $i$ nor $d+i$ is contained in $\sum^*$. Thus the left-hand side of (4.65) is at most

$$-\frac{1}{4d}\left[\sum_{k=1}^d n_k^2\right]^{-1/2}\sum_{i=1}^d |n_i| \le -\frac{1}{4d}I[\lambda - x \neq 0].$$

□

Proof of Theorem 2. We now have everything in place to prove Theorem 2. Fix $K > 0$ and a large $t$. We first use that the distinguished particle attempts to jump at the constant rate $D$. (It may, however, lose its distinguished character due to a jump.) Therefore, it holds for any $x$ that

$$\begin{aligned} &P\{\lambda(\cdot, x) \text{ has more than } 2Dt \text{ jumps during } [0,t]\} \\ &\quad \le P\{\widehat{\rho}(\cdot) \text{ attempts to jump more than } 2Dt \text{ times during } [0,t]\} \\ &\quad \le K_4 \exp[-K_5 t] \qquad [\text{see } (4.51)]. \end{aligned}$$

Note that if $\lambda(\cdot, x)$ has no more than $2Dt$ jumps during $[0,t]$, then also $\|\lambda(s,x) - \lambda(0,x)\| = \|\lambda(s,x) - z_0\|_\infty \le 2Dt$ for $s \le t$. Thus, for sufficiently large $t$ [see (4.9) for $\Xi$]

$$\begin{aligned} &P\left\{\{\lambda(s,x)\}_{s\le t} \notin \bigcup_{\ell \le 2Dt} \Xi(\ell, t) \text{ for some } x \text{ with } \|x\| \le t\right\} \\ &\qquad \le K_{11} t^d \exp[-K_5 t]. \end{aligned} \tag{4.68}$$

Proposition 8 (with $K$ replaced by $K+d$) now tells us that outside a further set of probability at most $2/t^{K+d}$, we have for $r \ge r_0, \ell \le 2Dt$ and $t \ge t_1$ that $\Phi_r(\ell) \le \varepsilon_0 C_0^{-6r}(1+2D)t$. Therefore, if we write $\widetilde{\lambda}(s,x)$ for the space–time point $\lambda(s,x) \times \{s\}$, then

$$\begin{aligned} &\text{for each } \|x\| \le t \text{ the path } \{\widetilde{\lambda}(s,x)\}_{s\le t} \text{ intersects} \\ &\text{at most } C_0^{-6r}\varepsilon_0(1+2D)t \text{ bad } r\text{-blocks} \end{aligned} \tag{4.69}$$

[see (4.7) for the definition of a bad block]. We choose $\varepsilon_0 > 0$ such that

$$\varepsilon_0(1+2D) < \frac{1}{7} \quad \text{and} \quad \frac{2\varepsilon_0}{K_7} e^{2D} < \frac{1}{6}, \tag{4.70}$$

and then we fix $r$ at some value $r_1 \ge r_0$ such that

$$\gamma_{r_1}\mu_A C_0^{dr_1} \ge 2, \qquad C_0^{6r_1} \ge \frac{3DK_9}{8} \tag{4.71}$$

[(4.71) is possible because $C_0 \ge 2$ and $\gamma_{r_1} \ge \gamma_0 > 0$]. We claim that with this choice, (4.69) implies (for large $t$) that

$$\begin{aligned} &\text{for each } \|x\| \le t \text{ and corresponding path } \{\lambda(s,x)\}_{s\le t} \text{ there are at least} \\ &(1/6)C_0^{-6r_1} t \text{ integers } 0 \le n \le C_0^{-6r_1} t/3 - 1 \text{ for which there exists a} \\ &\text{particle } \rho' \neq \widehat{\rho} \text{ inside } \lambda(3C_0^{6r_1} n, x) + [-C_0^{3r_1}, C_0^{3r_1}] \text{ at time } 3C_0^{6r_1} n. \end{aligned} \tag{4.72}$$

To see this recall that $\Delta_r = C_0^{6r}$, and note that each point $\widetilde{\lambda}(kC_0^{6r_1}, x)$ belongs to a unique $r_1$-block $\mathcal{B}_{r_1}(\mathbf{i}, k)$, and for different $k$, these blocks are disjoint. Thus for each $x$, $\{\lambda(s,x)\}_{s\le t}$ intersects at least $\lfloor C_0^{-6r_1} t\rfloor$ distinct $r_1$-blocks. If (4.69) holds, then, by (4.70), at most $(1/7)C_0^{-6r_1}t$ of these blocks are bad, so that for large $t$ there are at least $(6/7)C_0^{-6r_1}t - 4$ values of $0 \le k \le \lfloor C_0^{-6r_1}t\rfloor - 3$ such that $\widetilde{\lambda}(kC_0^{6r_1}, x)$ belongs to a good $r_1$-block. At least $(1/6)C_0^{-6r_1}t$ of these will have $k$ divisible by 3, say $k = 3n$, with $n \le \lfloor C_0^{-6r_1}t\rfloor/3 - 1$. If $\widetilde{\lambda}(3C_0^{6r_1}n, x)$ belongs to a good $r_1$-block, then by definition

$$\begin{aligned} U_{r_1}(\lambda(3C_0^{6r_1}n, x), 3C_0^{6r_1}n) &= \sum_{y\in\mathcal{Q}_{r_1}(\lambda^*)} N^*(y, 3C_0^{6r_1}n) \\ &\ge \gamma_{r_1}\mu_A C_0^{dr_1} \ge 2 \end{aligned}$$

[see (4.71)], where we have temporarily written $\lambda^*$ for $\lambda(3C_0^{6r_1}n, x)$. In particular, there have to be two particles in $\prod_{s=1}^d [\lambda^*(s), \lambda^*(s) + C_0^r)$ at time $3C_0^{6r_1}n$, and one of these must be different from the distinguished particle at $\lambda(3C_0^{6r_1}n)$. This justifies our claim (4.72).

In the notation of (4.57) the preceding paragraph shows that (4.69)–(4.71) imply that for all $x$ with $\|x\| \le t$ and

$$L = C_0^{3r_1},$$

it holds that

$$\begin{aligned} V(t, L, x) &\ge \frac{1}{6} C_0^{-6r_1} t > \frac{2\varepsilon_0}{K_7} e^{2D} C_0^{-6r_1} t \\ &= \frac{2\varepsilon_0}{K_7} e^{2D} L^{-2} t. \end{aligned}$$

Thus the bound (4.68) and the lines following it prove (for large $t$)

$$\begin{aligned} P\left\{V(t, L, x) < \frac{2\varepsilon_0}{K_7} e^{2D} L^{-2} t\right\} &\le K_{11} t^d \exp[-K_5 t] + 2t^{-K-d} \\ &\le 3t^{-K-d}. \end{aligned}$$

Equation (4.63) with $\varepsilon = \varepsilon_0$ and $z$ as in (4.61) then shows that for large $t$ and for all $x$ with

$$\|x - z_0\|_2 \le \frac{K_9 \varepsilon_0 D \beta(L,d)}{4L^2} t \quad \text{and} \quad \|x\| \le t, \tag{4.73}$$

$$\begin{aligned} &P\left\{\int_0^t I\left[\|\lambda(u,x) - x\|_2 \le \frac{4K_8 L^2}{K_9 \varepsilon_0 \beta(L,d)}\right] du \le \frac{K_9 \varepsilon_0 \beta(L,d)}{4K_8 L^2} t\right\} \\ &\quad \le 4t^{-K-d}. \end{aligned} \tag{4.74}$$

Now fix $x$ and assume that

$$\int_0^t I\left[\|\lambda(u,x)-x\|_2 \le \frac{4K_8L^2}{K_9\varepsilon_0\beta(L,d)}\right]du > \frac{K_9\varepsilon_0\beta(L,d)}{4K_8L^2}t. \tag{4.75}$$

This implies trivially that $\|\lambda(u,x)-x\|_2 \le 4K_8L^2/(K_9\varepsilon_0\beta(L,d))$ for some $u \le t$. However, we want more. For Theorem 2 we want that for each $x \in \mathcal{C}(C_2t)$ there exists a $u \le t$ at which $x$ is visited by a $B$-particle. To show that such a $u$ exists with high probability, we define the stopping times

$$\begin{aligned} u_0 = 0, \qquad u_{i+1} &= u_{i+1}(x) \\ &= \inf\left\{u \ge u_i + L^4\beta^{-2}(L,d) : \|\lambda(u,x)-x\|_2 \le \frac{4K_8L^2}{K_9\varepsilon_0\beta(L,d)}\right\}. \end{aligned}$$

Equation (4.75) implies that $u_i + L^4\beta^{-2}(L,d) \le t$ for at least

$$\chi := \left\lceil \frac{K_9\varepsilon_0\beta^3(L,d)}{4K_8L^6}t\right\rceil - 2$$

values of $i \ge 1$ with $u_i \le t$. By definition $u_{i+1}-u_i \ge L^4\beta^{-2}(L,d)$ and by the right-continuity of $\lambda(\cdot)$

$$\|\lambda(u_i,x)-x\|_2 \le \frac{4K_8L^2}{K_9\varepsilon_0\beta(L,d)} =: K_{12}\frac{L^2}{\beta(L,d)}. \tag{4.76}$$

Define

$$\widetilde{L} = \frac{L^2}{\beta(L,d)}.$$

Then the same random walk estimates as for (4.56) give that

$$\begin{aligned} &P\{\text{the particle which is the distinguished particle} \\ &\qquad \text{at time } u_i \text{ visits } x \text{ at some time in } (u_i,u_{i+1}]|\mathcal{F}_{u_i}\} \\ &\quad \asymp K_{13}\beta(\widetilde{L},d) \asymp \begin{cases} K_{13}, & \text{if } d=1, \\ K_{13}[4\log L]^{-1}, & \text{if } d=2, \\ K_{13}L^{d(d-2)}, & \text{if } d \ge 3. \end{cases} \end{aligned} \tag{4.77}$$

Note that we are estimating here the probability that the distinguished particle of time $u_i$ visits $x$ at some time $u$, rather than that $\lambda(u,x)=x$. But conditionally on $\mathcal{F}_{u_i}$, the $B$-particle which is the distinguished one at time $u_i$ performs a random walk which is a copy of $S$, and this fact is the basis for the estimate (4.77). Finally we apply Lemma 10 once more. We take $\mathcal{G}_n = \mathcal{F}_{u_n}$,

$$\begin{aligned} H_n := I[&\text{the particle which is the distinguished} \\ &\quad \text{particle at time } u_{n-1} \text{ visits } x \text{ at some time in } (u_{n-1},u_n]] \end{aligned}$$

and

$$D_n = H_n - E\{H_n|\mathcal{G}_{n-1}\}.$$

Since $H_n$ takes on only the values 0 or 1, it is easy to see from (4.77) that on the event (4.75)

$$\sum_{n=1}^{\chi} E\{H_n|\mathcal{G}_{n-1}\} \ge \chi K_{13}\beta(\widetilde{L},d) \ge K_{14}\beta(\widetilde{L},d)\frac{\beta^3(L,d)}{L^6}t$$

and

$$A_\chi = \sum_{n=1}^{\chi} E\{D_n^2|\mathcal{G}_{n-1}\} \le K_{15}\beta(\widetilde{L},d)\frac{\beta^3(L,d)}{L^6}t.$$

It is then easy to deduce from Lemma 10, with $a = (1/2)K_{14}\beta(\widetilde{L},d)\frac{\beta^3(L,d)}{L^6}t, b = K_{14}/(2K_{15}) \wedge 1$ and $c = 1$, that for $x$ satisfying (4.73) and $t$ large

$$\begin{aligned}
&P\{x \text{ is not visited by a } B\text{-particle by time } t\}\\
&\quad \le P\{(4.75) \text{ does not occur}\} + P\left\{(4.75) \text{ occurs, but } \sum_{i=1}^{\chi} H_n = 0\right\}\\
&\quad \le 4t^{-K-d}\\
&\qquad + P\Bigg\{(4.75) \text{ occurs, but}\\
&\qquad\qquad \sum_{i=1}^{\chi}[H_n - E\{H_n|\mathcal{G}_{n-1}\}] \le -K_{14}\beta(\widetilde{L},d)\frac{\beta^3(L,d)}{L^6}t\Bigg\}\\
&\quad \le 4t^{-K-d} + P\left\{\sum_{n=1}^{\chi}[H_n - E\{H_n|\mathcal{G}_{n-1}\}] \le -a - bA_\chi\right\}\\
&\quad \le 4t^{-K-d} + 2\exp\left[-K_{16}\beta(\widetilde{L},d)\frac{\beta^3(L,d)}{L^6}t\right]\\
&\quad \le 5t^{-K-d}.
\end{aligned} \tag{4.78}$$

Theorem 2 with

$$C_2 = \frac{K_9\varepsilon_0 D\beta(L,d)}{8L^2\sqrt{d}} \wedge 1$$

now follows by summing (4.78) over all $x$ which satisfy (4.73). $\square$

**5. Proof of Theorem 3.** A basic step for the proof is a monotonicity property which is proven via a coupling argument. We formulate it as a separate lemma.

Lemma 14. *Assume $D_A = D_B$ and let $\sigma^{(2)} \in \Sigma_0$. Assume further that $\sigma^{(1)}$ lies below $\sigma^{(2)}$ in the following sense:*

$$\text{for any site } z \in \mathbb{Z}^d\text{, all particles present in } \sigma^{(1)} \text{ at } z \text{ are also present in } \sigma^{(2)} \text{ at } z \tag{5.1}$$

*and*

$$\text{at any site } z \text{ at which the particles in } \sigma^{(2)} \text{ have type } A\text{, the particles also have type } A \text{ in } \sigma^{(1)}. \tag{5.2}$$

*Let $\pi_A(\cdot,\rho) = \pi_B(\cdot,\rho)$ be the random walk paths associated to the various particles and assume that the Markov processes $\{Y_t^{(1)}\}$ and $\{Y_t^{(2)}\}$ are constructed by means of the same set of paths $\pi_A(\cdot,\rho) = \pi_B(\cdot,\rho)$ and starting with state $\sigma^{(1)}$ and $\sigma^{(2)}$, respectively (as defined in Section 2). Then, almost surely, $\{Y_t^{(1)}\}$ and $\{Y_t^{(2)}\}$ satisfy (5.1) and (5.2) for all $t$ with $\sigma^{(i)}$ replaced by $Y_t^{(i)}, i = 1, 2$. In particular, $\sigma^{(1)} \in \Sigma_0$ and (2.20) holds almost surely for $\{Y_t^{(1)}\}$.*

Proof. Couple the processes $\{Y_t^{(1)}\}$ and $\{Y_t^{(2)}\}$ as in the statement of the lemma. Specifically, first choose independent paths $s \mapsto \pi_A(s,\rho)$ for all particles $\rho$ present in $\sigma^{(2)}$ and construct $\{Y_t^{(2)}\}$ with the help of these paths [as in (2.6) and (2.7), with $\pi_B(s,\rho) = \pi_A(s,\rho)$ for all $s,\rho$]. We then assign to each particle $\rho$ present in $\sigma^{(1)}$ the same path $s \mapsto \pi_A(s,\rho)$ as assigned to $\rho$ in $\sigma^{(2)}$. By (5.1) this assigns a path to each particle present in $\sigma^{(1)}$. We then construct $\{Y_t^{(1)}\}$ on the basis of these paths. Note that a.s. no particles are ever moved to a cemetery point in the $\{Y_t^{(2)}\}$-system, since $\sigma^{(2)} \in \Sigma_0$. The position at time $t$ of a particle $\rho$ starting at $z$ is then $z + \pi_A(t,\rho)$, in whichever of the systems the particle is present. It is immediate from this that (5.1) with $\sigma^{(i)}$ replaced by $Y_t^{(i)}, i = 1, 2$, is valid.

To show (5.2) with $\sigma^{(i)}$ replaced by $Y_t^{(i)}, i = 1, 2$, we first note that (2.20) a.s. holds for $\{Y_t^{(2)}\}$, because $\sigma^{(2)} \in \Sigma_0$. Then, by (5.1), a.s. (2.20) holds in both systems. Now let $\tau_0^{(i)} = 0$ and for $k \geq 1$, let $\tau_k^{(i)}$ be the $k$th time at which a new particle changes from type $A$ to type $B$ in $\{Y_t^{(i)}\}$. More formally, as in (2.3), $\tau_{k+1}^{(i)} = \inf\{t > \tau_k^{(i)}:$ a $B$-particle coincides with an $A$-particle at time $t$ in $\{Y_t^{(i)}\}\}$. We shall show by induction on $k \geq 0$ that at each time $\tau_k^{(1)}$ the property (5.2) with $\sigma^{(i)}$ replaced by $Y^{(i)}_{\tau_k^{(1)}}$ still holds, and that there are only finitely many $B$-particles in both systems at time $\tau_k^{(1)}$. We may restrict ourselves to sample points for which $\min(\widehat{\tau}, \tau_\infty) = \infty$ in the $\{Y_t^{(2)}\}$-system, because $\sigma^{(2)} \in \Sigma_0$ (see Lemma 2). Assume then that at time $\tau_k^{(1)}$, (5.2) with $\sigma^{(i)}$ replaced by $Y^{(i)}_{\tau_k^{(1)}}$ still holds. Since the second system has only finitely many $B$-particles at time $\tau_k^{(1)}$, this, together with (5.2) with $\sigma^{(i)}$ replaced by $Y^{(i)}_{\tau_k^{(1)}}$, shows that also the first system has only finitely many $B$-particles at time $\tau_k^{(1)}$. Moreover, $\tau_{k+1}^{(1)}$ is the first time after $\tau_k^{(1)}$ at which some $B$-particle $\rho$ coincides with some $A$-particle $\rho'$ in the first system. From the right-continuity of the paths $\pi_A(\cdot, \zeta)$ for all $B$-particles $\zeta$, plus (2.20), it then follows that $\tau_{k+1}^{(1)} > \tau_k^{(1)}$, a.s. By the induction hypothesis, $\rho$ must also have type $B$ at time $\tau_k^{(1)}$ in the second system. Therefore also $\rho'$ turns into a $B$-particle in the second system no later than $\tau_{k+1}^{(1)}$. ($\rho'$ may already have turned to type $B$ before $\tau_{k+1}^{(1)}$ in the second system, in which case no type change occurs for $\rho'$ in the second system at $\tau_{k+1}^{(1)}$.) In any case, any particle which turns to type $B$ at time $\tau_{k+1}^{(1)}$ in the first system also has type $B$ at or before time $\tau_{k+1}^{(1)}$ in the second system, so that (5.2) with $\sigma^{(i)}$ replaced by $Y^{(i)}_{\tau_{k+1}^{(1)}}$ still holds. This completes the inductive step.

Since we already know that there are only finitely many $B$-particles at each $\tau_k^{(1)}$ in the second system, we conclude that a.s. this holds in both systems at each $\tau_k^{(1)}$. As we remarked above this shows that a.s. $\tau_{k+1}^{(1)} > \tau_k^{(1)}$ for all $k$, so that $\widehat{\tau} < \tau_\infty$ has probability 0 in both systems. Also, at each $\tau_k^{(1)}$, the number of $B$-particles in the first system is at most equal to the number of $B$-particles in the second system. Since there are at least $k$ $B$-particles in the first system at time $\tau_k^{(1)}$, this shows that

$$P^{\sigma^{(1)}}\{\tau_\infty < \infty\} \le P^{\sigma^{(2)}}\{\tau_\infty < \infty\} = 0.$$

Thus, $\sigma^{(1)} \in \Sigma_0$. $\square$

Proof of Theorem 3. Fix $K$. Note that if a particle has type $B$ at some time $s \le t$ and is outside the cube $\mathcal{C}(C_1 t)$ at that time, then by symmetry of the random walk $\{S\}$, the particle has a conditional probability, given $\mathcal{F}_s$, at least 1/2 of being outside $\mathcal{C}(C_1 t)$ at time $t$. Therefore

$$\begin{aligned}&E\{\text{number of particles outside } \mathcal{C}(C_1 t) \text{ at some time } s \le t\}\\&\quad \le 2E\{\text{number of particles outside } \mathcal{C}(C_1 t) \text{ at time } t\}.\end{aligned}$$

Thus, by (1.3)

$$P\{\text{a site outside } \mathcal{C}(C_1 t) \text{ is visited by a } B\text{-particle during } [0,t]\} \le t^{-K-1} \tag{5.3}$$

for $t \ge$ some $t_0$. We may therefore restrict ourselves for (1.7) to space–time points $(z,t)$ with $z \in \mathcal{C}(C_1 t)$. Now fix a $(z,t)$ which satisfies this condition and assume $z$ is first visited by a $B$-particle at time $s \le t - [K_1 t \log t]^{1/2}$. Outside a set of probability $Dt^{-K-d-1}$ this $B$-particle stayed at $z$ for at least $t^{-K-d-1}$ units of time, so that it is still at $z$ at a time $s$ of the form $kt^{-K-d-1} \le t - [K_1 t \log t]^{1/2} + t^{-K-d-1}$. Further, the probability that there is an $A$-particle at $(z,t)$ which was at a point $y$ with $\|y-z\| > [K_2 t \log t]^{1/2}$ at one of the times $kt^{-K-d-1}$ is bounded by

$$\sum_{\substack{s=kt^{-K-d-1} \\ \le t-[K_1 t\log t]^{1/2}+1}} \sum_{y:\, \|y-z\|>[K_2 t \log t]^{1/2}} EN_A(y,s) P\{S_{t-s} = z-y\}. \tag{5.4}$$

We now remind the reader of the particle system $\mathcal{P}^*$ which we introduced just before (2.16). In this process interactions between particles are ignored. Thus, in $\mathcal{P}^*$, each $A$-particle $\rho$ with initial position $\pi(0,\rho)$ continues to follow the path $t \mapsto \pi(0,\rho) + \pi_A(t,\rho)$ even after its switching time $\theta(\rho)$. In the present case with $D_A = D_B$, this is also the path which the particle follows in the true system. We write $N^*(z,t)$ for the number of particles at $(z,t)$ in $\mathcal{P}^*$. In our case this is just the total number of particles at $(z,t)$ which are different from the finitely many original $B$-particles. $\{N^*(z,t) : z \in \mathbb{Z}^d, t \ge 0\}$ is stationary in time, and at each $t$ the $N^*(x,t), x \in \mathbb{Z}^d$, are i.i.d. mean-$\mu_A$ Poisson variables. From this description we see that

$$N_A(z,t) \le N^*(z,t) \quad \text{and} \quad N_A(z,t) + N_B(z,t) \ge N^*(z,t), \qquad z \in \mathbb{Z}^d, t \ge 0. \tag{5.5}$$

In particular, $EN_A(y,s) \le EN^*(y,s) = \mu_A$. Therefore, (5.4) is bounded by

$$\sum_{\substack{s=kt^{-K-d-1} \\ \le t-[K_1 t\log t]^{1/2}+1}} \mu_A P\{\|S_{t-s}\| > [K_2 t \log t]^{1/2}\} \le K_3 t^{K+d+2} \exp[-K_4 K_2 \log t] \tag{5.6}$$

(see (2.42) in [7]). Taking into account the number of possibilities for $z$ we find that the probability in the left hand-side of (1.7) is for large $t$ bounded by

$$\frac{1}{t^{K+1}} + K_5 t^d \frac{1}{t^{K+d+1}} + K_6 t^{K+2d+2} \exp[-K_4 K_2 \log t]$$

$$+ \sum_{\substack{s=kt^{-K-d-1} \\ \leq t-[K_1 t \log t]^{1/2}+1}} \sum_{z \in \mathcal{C}(C_1 t)} \sum_{y:\, \|y-z\| \leq [K_2 t \log t]^{1/2}} P\{\text{there is a } B\text{-particle at } z \text{ and an } A\text{-particle at } y \text{ at time } s = kt^{-K-d-1}, \text{ but this } A\text{-particle is not turned into a } B\text{-particle by time } t\}.$$

We shall write $\mathcal{U}(k, z, y)$ for the event in the summand corresponding to $k, z, y$ in the triple sum here. We further write $s$ for $kt^{-K-d-1}$. We choose $K_2$ so large that for $t \geq$ some $t_1$ the sum of the first three terms here is at most $1/(2t^K)$. The probability in the left-hand side of (1.7) is then bounded by

$$\frac{1}{2t^K} + \sum_{\substack{s=kt^{-K-d-1} \\ \leq t-[K_1 t \log t]^{1/2}+1}} \sum_{z \in \mathcal{C}(C_1 t)} \sum_{y:\|y-z\| \leq [K_2 t \log t]^{1/2}} P\{\mathcal{U}(k, z, y)\}. \tag{5.7}$$

The triple sum in (5.7) contains at most $K_7 t^{K+d+2}[C_1 t]^d [K_2 t \log t]^{d/2} \leq K_8 t^{K+3d+2}$ summands. We shall complete the proof of (1.7) by showing that each summand in (5.7) is at most $K_9 t^{-2K-3d-3}$ for some constant $K_9$ which does not depend on $t$.

At an intuitive level it seems clear that there should be a way to estimate these summands which is similar to the one used in the proof of Theorem 2 for estimating the probability that a fixed vertex $x$ has not been reached by a $B$-particle by time $t$. We have to be a bit careful, though, not to bias the relevant distributions by the fact that there is a $B$-particle at the space–time point $(z, s)$ or a certain number of $A$-particles at $(y, s)$. Nevertheless, the proof will follow Theorem 2 closely. As before, we shall write $\mathcal{P}_0$ for our original particle system.

In order to estimate the summand in (5.7) we shall make use of the monotonicity property in Lemma 14. We shall compare the real system with a modified system, which is constructed via two modifications. In each of these modifications we remove particles and change some $B$-particles to $A$-particles, at time $s$. According to the monotonicity property of Lemma 14, any particle which is present in both systems at time $t$, and which has type $A$ in the unmodified system, also must have type $A$ in the modified system. Now fix $y$ and $z$ as well as $s < t$, and make the first modification as follows: at time $s$, at each $x \neq z$ remove all particles of type $B$ which were added at time 0, and reset the type of the remaining particles at $x$ to $A$. At $x = z$, if there are particles at $(z, s)$ give them type $B$ at time $s$. If there is no particle at $(z, s)$, put one $B$-particle at $z$ at time $s$. If there is no $B$-particle at $z$ or no $A$-particle at $y$ at time $s$ before the modification, then $\mathcal{U}(k, z, y)$ does not occur in the original system, so we do not care what the modification does in this case. In the other cases there is at least one $A$-particle at $y$, and a $B$-particle at $z$, so $y \neq z$. In these cases the $A$-particle at $y$ is not removed, and the type of the particles at $z$ is unchanged by the modification, so the monotonicity property gives us that $P\{\mathcal{U}(k,z,y)\}$ can only go up by this first modification. Note that after the first modification, we have $N^*(x,s)$ $A$-particles at $x$, for any $x \neq z$, and $N_B(z,s) \vee 1$ $B$-particles at $z$. These

$$N^*(x,s), x \in \mathbb{Z}^d, \text{ are i.i.d. mean-}\mu_A \text{ Poisson variables.} \tag{5.8}$$

Now there are $N^*(y,s)$ $A$-particles at $(y,s)$ after the first modification. Since all particles at the same space–time point play the same role, each of these $N^*(y,s)$ particles at $(y,s)$ has the same probability of still having type $A$ at time $t$. Order the particles at $(y,s)$ by some arbitrary rule. Then the $k,z,y$ summand in (5.7) is at most

$$\begin{aligned} &E\{N^*(y,s)P\{\text{in the first modified system the first} \\ &\qquad \text{particle at } (y,s) \text{ is still of type } A \text{ at time } t|\mathcal{F}_s\}\}. \end{aligned} \tag{5.9}$$

There is some dependence between $N^*(y,s)$ and the conditional probability factor in this expectation. To handle this we make one further modification at $y$. In this second modification we remove all but the first $A$-particle at $(y,s)$. If there is no $A$-particle at $(y,s)$, then we add an $A$-particle at $y$. Since there is no contribution to (5.9) from the sample points with $N^*(y,s) = 0$, we can ignore how this modification influences the conditional probability in (5.9) on the event $\{N^*(s,y) = 0\}$. On the event $\{N^*(y,s) \geq 1\}$ the second modification cannot decrease the conditional probability in (5.9), once again by Lemma 14. Now, the second modified system does not involve $N^*(y,s)$ anymore. In the second modified system the conditional probability in (5.9) is replaced by

$$\begin{aligned} &P\{\text{in the second modified system, the unique} \\ &\qquad \text{particle at } (y,s) \text{ is still of type } A \text{ at time } t|\mathcal{F}_s\}, \end{aligned} \tag{5.10}$$

which is a function of the $N^*(u,s)$ with $u \neq y$ only. As we already observed, these are independent mean-$\mu_A$ Poisson variables, independent of $N^*(y,s)$. Consequently (5.9) is at most

$$\begin{aligned} &E\{N^*(y,s)\}P\{\text{in the second modified system the unique} \\ &\qquad \text{particle at } (y,s) \text{ is still of type } A \text{ at time } t\}. \end{aligned} \tag{5.11}$$

The statement (5.11) suggests that we introduce two systems $\mathcal{P}^{z,y}$ and $\mathcal{P}^z$ say, which have initial distributions as follows: Let $\widetilde{N}(u), u \in \mathbb{Z}^d$, be a family of independent mean-$\mu_A$ Poisson variables. Then $\mathcal{P}^{z,y}$ starts with $\widetilde{N}(u)$ $A$-particles at $u$ if $u \notin \{z,y\}, \widetilde{N}(z) \vee 1$ $B$-particles at $z$ and one $A$-particle at $y$. No other particles are in the initial state. For $\mathcal{P}^z$ the only change is that at $y$ we put initially $\widetilde{N}(y)$ particles, so that $y$ is treated like all other sites $u \neq z$. The particles then move and change type in $\mathcal{P}^{z,y}$ and in $\mathcal{P}^z$ according to the same rules as in $\mathcal{P}_0$. It follows from (5.8) that the probability factor in (5.11) equals

$$\begin{aligned}
&P\{\text{in the system } \mathcal{P}^{z,y} \text{ the unique} \\
&\quad \text{particle at } (y,0) \text{ is still of type } A \text{ at time } t-s\} \\
&\quad = P\{\text{in the system } \mathcal{P}^{z} \text{ the unique particle at } (y,0) \\
&\qquad \text{is still of type } A \text{ at time } t-s \mid \text{start with one } A\text{-particle at } y\} \\
&\quad \le \frac{1}{e^{-\mu_A}\mu_A} P\{\text{in the system } \mathcal{P}^{z} \text{ the first particle} \\
&\qquad \text{at } (y,0) \text{ is still of type } A \text{ at time } t-s\}.
\end{aligned} \tag{5.12}$$

Here $e^{-\mu_A}\mu_A$ in the right-hand side represents $P\{\widetilde{N}(y) = 1\}$. We further have $E\{N^*(y,s)\} = \mu_A$, and $t - s \ge [K_1 t \log t]^{1/2} - 1$, so that the $k, z, y$ summand in (5.7) is bounded by

$$\begin{aligned}
&e^{\mu_A} P\{\text{the first particle at } (y,0) \text{ is still of type } A \\
&\qquad \text{at time } [K_1 t \log t]^{1/2} - 1 \text{ in the system } \mathcal{P}^z\}.
\end{aligned} \tag{5.13}$$

Our task has now been reduced to estimating the probability factor in (5.13). But the system $\mathcal{P}^z$ is either equal to the original system $\mathcal{P}_0$ with one $B$-particle added at $z$ [in case $\widetilde{N}(z) = 0$] or is like the system $\mathcal{P}_0$, but with all particles at $z$ turned into $B$-particles at time 0, without the addition of an extra $B$-particle. In both these cases we can basically repeat the proof of Theorem 2 to estimate (5.13). We form a path $s \mapsto \lambda(s) = \lambda(s, y, z), s \ge 0$, which starts at $\lambda(0, y, z) = z$ and from there proceeds according to the rules (i)–(v) in the proof of Theorem 2 (after the proof of Proposition 8) with only one change in rule (v). We now want the path to have a drift to the first particle which started at $y$, instead of to a fixed vertex $x$. If we denote this first particle at $y$ by $\phi$, then the position of $\phi$ at a time $s$ equals $\pi(s, \phi) = y + \pi_A(s, \phi)$. Accordingly we change rule (v) to the following:

(v′) if $\widehat{\rho}(s-)$ jumps from $\lambda(s-) = w$ to $w'$ at some time $s$ such that there is at least one other particle $\rho'$ at $w$ at time $s-$, then $\lambda(\cdot)$ jumps to $w'$ at time $s$ if and only if $\|w' - \pi(s,\phi)\|_2 < \|w - \pi(s,\phi)\|_2$, and in this case again $\widehat{\rho}(s) = \widehat{\rho}(s-)$; if, however, $\|w' - \pi(s,\phi)\|_2 \ge \|w - \pi(s,\phi)\|_2$, then $\lambda(\cdot)$ does not jump at time $s$ and we take $\widehat{\rho}(s) = \rho'$.

Under these changed rules $\lambda(\cdot) - \pi(\cdot, \phi)$ has a drift toward zero in the sense that now

$$\widetilde{M}(t) := \|\lambda(t,x) - \pi(t,\phi)\|_2$$

(5.14)
$$- D\int_0^t [I_1(u)\widetilde{\Gamma}_1(u) + I_{\geq 2}(u)\widetilde{\Gamma}_{\geq 2}(u)]du - D\int_0^t \widetilde{\Gamma}_1(u)\,du$$

is an $\{\mathcal{F}_t\}$-martingale, where analogously to (4.43)

$$\widetilde{\Gamma}_1(u) := \frac{1}{2d}\sum_{i=1}^{2d}[\|\lambda(u) + e_i - \pi(u,\phi)\|_2 - \|\lambda(u) - \pi(u,\phi)\|_2],$$

$$\widetilde{\Gamma}_{\geq 2}(u) := \frac{1}{2d}\sum{}^*[\|\lambda(u) + e_i - \pi(u,\phi)\|_2 - \|\lambda(u) - \pi(u,\phi)\|_2],$$

and $\sum^*$ is the sum over those $i \in \{1,\dots,2d\}$ for which

$$\|\lambda(u) + e_i - \pi(u,\phi)\|_2 - \|\lambda(u) - \pi(u,\phi)\|_2 < 0,$$

and $e_{d+i} = -e_i, 1 \leq i \leq d$; $I_1(u)$ and $I_{\geq 2}(u)$ are the same as in (4.42). The extra integral $D\int_{[0,t]}\widetilde{\Gamma}_1(u)\,du$ [which was not present in the $M(\cdot)$ of (4.44)] has to be introduced to compensate for the jumps of $\phi$. However, the proof that $\widetilde{M}$ is a martingale is quite the same as for Lemma 9.

From here on one can follow the proof of Theorem 2. One merely has to replace $\lambda(s) - x$ by $\lambda(s) - \pi(s,\phi)$ at most places and to take into account the jumps in $\widetilde{M}(\cdot)$ due to a jump of $\phi$. For instance, the definition of $\sigma_{k+1}$ right after (4.48) should now be

$$\sigma_{k+1} = \min[\sigma_k + 1, \inf\{t > \sigma_k : t \text{ the distinguished particle attempts a jump at } t \text{ or } \phi \text{ jumps at } t\}].$$

One also has to note that $t$ should be replaced by $[K_1 t\log t]^{1/2} - 1$ in everything that comes after Lemma 10, but this is a trivial change to make. The conclusion is that if $K_1$ is chosen sufficiently large with respect to $K_2$, then each of the summands in (5.7) is bounded by $K_9 t^{-3d-2K-3}$, uniformly in $k, y, z$ in the ranges over which they are summed. As pointed out before this completes the proof of (1.7).

Once we have (1.7), (1.8) easily follows by means of (1.5). Indeed, by (1.5), outside a set of probability at most $t^{-K}$, each point in $\mathcal{C}((C_2/2)t)$ has been visited by a $B$-particle during $[0, t/2]$. We then have by (1.7) that the additional probability of some vertex $z \in \mathcal{C}((C_2/2)t)$ being occupied by an $A$-particle at time $t$ is at most $t^{-K}$. Thus, for large $t$ the left-hand side of (1.8) is at most $2t^{-K}$. □

**Acknowledgments.** Much of the research by H. Kesten for this paper was carried out at the Mittag–Leffler Institute in Djursholm, while he was supported by a Tage Erlander Professorship. H. Kesten thanks the Swedish Research Council for awarding him a Tage Erlander Professorship for 2002 and the Mittag–Leffler Institute in Djursholm for providing him with excellent facilities and for its hospitality.

V. Sidoravicius thanks Cornell University and the Mittag–Leffler Institute for their hospitality and travel support.

Both authors thank Balint Toth, Antal Jarai and J. van den Berg for many helpful discussions. They are also grateful to an Associate Editor and a referee (both anonymous) for a careful reading of the paper and for suggestions for improvements of the exposition.

Department of Mathematics
Malott Hall, Cornell University
Ithaca, New York 14853
USA
e-mail: kesten@math.cornell.edu

Instituto Nacional De Matemática
Pura E Aplicada
Estrada Dona Castorina 110
Rio de Janeiro
Brasil
e-mail: vladas@impa.br